\documentclass[a4paper,11pt]{article}
\usepackage{amsmath}
\usepackage{amssymb}
\numberwithin{equation}{section}
\newtheorem{theo}{Theorem}[section]
\newtheorem{cor}{Corollary}[theo]
\newtheorem{lem}{Lemma}[section]
\newtheorem{prop}{Proposition}[section]

\newcommand{\Proof}[1]{\noindent{\textit{Proof.}}#1 $\quad\Box$
{\vspace{.2cm}}}
\newcommand{\Rem}[2]{\noindent{\bf Remark #1} {#2} {\vspace{.2cm}}}

\newcommand{\lam}{\lambda}
\newcommand{\ep}{\epsilon}
\newcommand{\sig}{\sigma}
\newcommand{\ta}{\theta}

\newcommand{\sigpmj}{\sig^\pm_j}

\newcommand{\sigpj}{\sig^+_j}

\newcommand{\sigmj}{\sig^-_j}

\newcommand{\phipmj}{\phi^\pm_j}
\newcommand{\phipj}{\phi^+_j}
\newcommand{\phimj}{\phi^-_j}
\newcommand{\psipmj}{\psi^\pm_j}

\newcommand{\Abf}{\mathbf{A}}
\newcommand{\BB}{\mathbf{B}}

\newcommand{\DD}{\mathbf{D}}

\newcommand{\JJ}{\mathbf{J}}
\newcommand{\KK}{\mathbf{K}}
\newcommand{\LL}{\mathbf{L}}
\newcommand{\II}{\mathbf{I}}
\newcommand{\MM}{\mathbf{M}}

\newcommand{\PP}{\mathbf{P}}

\newcommand{\VV}{\mathbf{V}}
\newcommand{\ZZ}{\mathbf{Z}}
\newcommand{\R}{\mathbb{R}}
\newcommand{\C}{\mathbb{C}}
\newcommand{\Z}{\mathbb{Z}}
\newcommand{\E}{\mathrm{E}_{\sigma ,\ep}}
\newcommand{\WE}{\widetilde{\mathrm{E}_{\mu ,\ep}}}
\newcommand{\Om}{\Omega}

\newcommand{\cir}{\mathbb{S}^1}

\newcommand{\cl}{\mathcal{L}}
\newcommand{\clep}{\cl_\ep}
\newcommand{\clepstar}{\cl_\ep^\ast}
\newcommand{\clstar}{\cl^\ast}
\newcommand{\specl}{\mathrm{Spec}(\cl)}
\newcommand{\speclep}{\mathrm{Spec}(\clep)}

\newcommand{\specb}{\mathrm{Spec}}
\newcommand{\ind}{\mathrm{Ind}}
\newcommand{\charac}{\mathrm{Char}}
\newcommand{\tr}{\mathrm{tr}}

\newcommand{\dis}{\displaystyle}
\newcommand{\ov}[1]{\overline{#1}}
\newcommand{\pa}{\partial}
\newcommand{\dd}[2]{\dis\frac{\pa #1}{\pa #2}}
\newcommand{\begeq}{\begin{equation}}
\newcommand{\stopeq}{\end{equation}}
\newcommand{\begar}{\begin{array}}
\newcommand{\stopar}{\end{array}}
\newcommand{\ei}[1]{\mathrm{e}^{#1}}

\title{On First and Second Order \\
Planar Elliptic Equations with Degeneracies}

\author{ Abdelhamid Meziani\\
Department of Mathematics\\
Florida International University\\
Miami, Florida 33199\\
\\
email: meziani@fiu.edu}

\date{}

\begin{document}

\maketitle

\noindent{\bf Abstract.} This paper deals with elliptic equations in
the plane with degeneracies. The equations are generated by a
complex vector field that is elliptic everywhere except along a
simple closed curve. Kernels for these equations are constructed.
Properties of solutions, in a neighborhood of the degeneracy curve,
are obtained through integral and series representations. An
application to a second order elliptic equation with a punctual
singularity is given.

\vspace{2cm}


\newpage

\begin{center}
CONTENTS
\end{center}

\begin{itemize}
\item[] Introduction
\item[1] Preliminaries
\item[2] Basic Solutions
\begin{itemize}
\item[2.1] Properties of basic solutions
\item[2.2] The spectral  equation and $\mathrm{Spec}(\cl_0)$
\item[2.3] Existence of basic solutions
\item[2.4] Properties of the fundamental matrix
\item[2.5] The system of equations for the
adjoint operator
\item[2.6] Continuation of a simple spectral value
\item[2.7] Continuation of a double spectral value
\item[2.8] Purely imaginary spectral values
\item[2.9] Main result about basic solutions
\end{itemize}
\item[3] Example
\item[4] Asymptotic Behavior of the Basic Solutions
\begin{itemize}
\item[4.1] Estimate of $\sigma$
\item[4.2] First estimate of $\phi$ and $\psi$
\item[4.3] End of the proof of Theorem 4.1
\end{itemize}
\item[5] The Kernels
\begin{itemize}
\item[5.1] Two lemmas
\item[5.2] Proof of Theorem 5.1
\item[5.3] Modified kernels
\end{itemize}
\item[6] The Homogeneous Equation $\cl u=0$
\begin{itemize}
\item[6.1] Representation of solutions in a cylinder
\item[6.2] The Cauchy Integral Formula
\item[6.3] Consequences
\end{itemize}
\item[7] The Nonhomogeneous Equation $\cl u=F$
\begin{itemize}
\item[7.1] Generalized Cauchy Integral Formula
\item[7.2] The integral operator $T$
\item[7.3] Compactness of the operator $T$
\end{itemize}
\item[8] The Semilinear Equation
\item[9] The Second Order Equation: Reduction
\item[10] The Homogeneous Equation $Pu=0$
\begin{itemize}
\item[10.1] Some properties
\item[10.2] Main result about the homogeneous equation
\item[10.3] A maximum principle
\end{itemize}
\item[11] The Nonhomogeneous Equation $Pu=F$
\item[12] Normalization of a Class of Second Order Equations with
a  Singularity
\item[] References
\end{itemize}

\section*{Introduction}
This paper deals with the properties of solutions of first and
second order equations in the plane. These equations are generated
by a complex vector field $X$ that is elliptic everywhere except
along a simple closed curve $\Sigma\subset \R^2$. The vector field
$X$ is tangent to $\Sigma$ and $X\wedge\ov{X}$ vanishes to first
order along $\Sigma$ (and so $X$ does not satisfy H\"{o}rmander's
bracket condition). Such vector fields have canonical
representatives (see {\cite{Mez-JFA}}). More precisely, there is a
change of coordinates in a tubular neighborhood of $\Sigma$ such
that $X$ is conjugate to a unique vector field $L$ of the form
 \begeq
 L=\lam\dd{}{t}-ir\dd{}{r}
 \stopeq
defined in a neighborhood of the circle $r=0$ in $\R\times\cir$,
where $\lam\in\R^++i\R$ is an invariant of the structure generated
by $X$. We should point out that normalizations for vector fields
$X$ such that  $X\wedge\ov{X}$ vanishes to a constant order $n>1$
along $\Sigma$ are obtained in {\cite{Mez-TAMS}}, but we will
consider here only the case $n=1$. This canonical representation
makes it possible to study the equations generated by $X$ in a
neighborhood of the characteristic curve $\Sigma$. We would like to
mention a very recent paper by F. Treves {\cite{Trev-New}} that uses
this normalization to study hypoellipticity and local solvability of
complex vector fields in the plane near a linear singularity. The
motivation for our work stems
 from the theory of hypoanalytic structures (see {\cite{Trev-Book}} and
 the references therein)
 and from the theory of generalized analytic functions
 (see {\cite{Vek}}).

The equations considered here are of the form $Lu=F(r,t,u)$.
 and $Pu=G(r,t,u,Lu)$, where $P$ is the (real) second
order operator
 \begeq
 P=L\ov{L}+\beta(t)L+\ov{\beta}(t)\ov{L}\, .
 \stopeq
 It should be noted that very little is known, even locally,
 about the structure of the solutions of second order equations generated by
 complex vector fields with degeneracies. The paper {\cite{Gil-Trev}}
 explores the local solvability of a particular case generated by a
 vector field of finite type.

 An application to a class of second order elliptic operators with
 a punctual singularity in $\R^2$ is given. This class consists of
 operators of the form
 \begeq
 D=a_{11}\dd{^2}{x^2}+2a_{12}\dd{^2}{xy}+a_{22}\dd{^2}{y^2}
 +a_1\dd{}{x}+a_2\dd{}{y}\, ,
 \stopeq
 where the coefficients are real-valued, smooth, vanish at 0, and satisfy
 \[
 C_1\,\le\,\frac{a_{11}a_{22}-a_{12}^2}{(x^2+y^2)^2}\,\le \,C_2
 \]
 for some positive constants $C_1\le C_2$. It turns out that each
 such operator $D$ is conjugate in $U\backslash 0$ (where $U$ is
 an open neighborhood of $0\in\R^2$) to a multiple of an
 operator $P$ given by (0.2).

Our approach is based on a thorough study of the operator $\cl$
given by
 \begeq
 \cl u=Lu-A(t)u-B(t)\ov{u}\, .
 \stopeq
 For the equation $\cl u=0$, we introduce particular solutions,
 called here basic solutions. They have the form
 \[
 w(r,t)= r^\sigma \phi(t)+\ov{r^\sigma \psi(t)}\, ,
 \]
 where $\sigma\in\C$ and $\phi(t)$, $\psi(t)$ are $2\pi$-periodic
 and $\C$-valued. Sections 2 and 4 establish the main properties of
 the basic solutions. In particular, we show that for every
 $j\in\Z$, there are (up to real multiples) exactly two
 $\R$-independent basic solutions
 \[
 w_j^\pm(r,t)=r^{\sigpmj} \phipmj(t)+\ov{r^{\sigpmj} \psipmj(t)}
 \]
with winding number $j$. For a given $j$, if
$\sigpj\in\C\backslash\R$, then $\sigmj=\sigpj$; and if
$\sigpj\in\R$ then we have only $\sigmj\le\sigpj$. The basic
solutions play a fundamental role in the structure of the space of
solutions of the equation $\cl u=F$.

In section 6, we show that any solution of $\cl u=0$ in a cylinder
$(0,\ R)\times\cir$ has a Laurent type series expansion in the
$w_j^\pm$'s. From the basic solutions of $\cl$ and those of the
adjoint operator $\clstar$, we construct, in section 5, kernels
$\Om_1$ and $\Om_2$ that allow us to obtain a Cauchy Integral
Formula (section 6)
 \begeq
 u(r,t)=\int_{\pa_0U}\Om_1u\, \frac{d\zeta}{\zeta}+
 \ov{\Om_2u}\, \frac{d\ov{\zeta}}{\ov{\zeta}}
 \stopeq
that represents the solution $u$ of $\cl u=0$ in terms of its values
on the distinguished boundary $\pa_0U=\pa U\backslash\Sigma$.

For the nonhomogeneous equation $\cl u =F$, we construct, in section
7,  an integral operator $T$, given by
 \begeq
 TF=\frac{-1}{2\pi}\iint_U\left(\Om_1F+\ov{\Om_2F}\right)
 \frac{d\rho d\ta}{\rho}\, .
 \stopeq
 This operator produces H\"{o}lder continuous solutions (up to the
 characteristic circle $\Sigma$), when $F$ is in an appropriate
 $L^p$-space. The properties of $T$ allow us to establish, in section 8,
  a similarity principle between the solutions of the homogeneous
 equations $\cl u=0$ and those of a semilinear equation
 $\cl u=F(r,t,u)$

The properties of the (real-valued) solutions of $Pu=G$ are studied
in sections 9 to 11. To each  function $u$ we associate a complex
valued function $w=BLu$, called here the $L$-potential of $u$, and
such that $w$ solves an equation of the form $\cl w =F$. The
properties of the solutions of $Pu=G$ can  thus be understood in
terms of the properties of their $L$-potentials. In particular
series representations and integral representations are obtained for
$u$. A maximum principle for  the solutions of $Pu=0$ holds on the
distinguished boundary $\pa_0U$, if the spectral values $\sigpmj$
satisfy a certain condition. In the last section, we establish the
conjugacy between the operator $D$ and the operator $P$.

\section{Preliminaries}
We start by reducing the main equation $Lu=Au+B\ov{u}$ into a
simpler form. Then, we define a family of operators $\clep$,  their
adjoint $\clepstar$, and prove a Green's formula. The operators
$\clep$ will be extensively used in the next section.

 Let $\lam =a+ib\in \R^++i\R^\ast$ and define the vector field $L$
 by
 \begeq
 L=\dis\lam\dd{}{t}-ir\dd{}{r}\,.
 \stopeq
For $ A\in C^k(\cir,\C)$, with $k\in\Z^+$, set
\[
A_0=\frac{1}{2\pi}\int_0^{2\pi}A(t)dt,\qquad \nu =
1-\mathrm{Im}\frac{A_0}{\lam}+\left[\mathrm{Im}\frac{A_0}{\lam}\right]
\]
where for $x\in\R$, $\ [ x]$ denotes the greatest integer less or
equal than $x$. Hence, $\nu\in [0,\ 1)$. Define the function
\[
m(t)=\exp\left(it+i\left[\mathrm{Im}\frac{A_0}{\lam}\right]t
+\frac{1}{\lam}\int_0^t(A(s)-A_0)ds \right)\, .
\]
Note that $m(t)$ is $2\pi$-periodic. The following lemma is easily
verified.

{\lem Let $A,B\in C^k(\cir,\C)$ and $m(t)$ be as above. If $u(r,t)$
is a solution of the equation
\begin{equation}
Lu=A(t)u+B(t)\ov{u}
\end{equation}
then the function $w(r,t)=\dis\frac{u(r,t)}{m(t)}$ solves the
equation
\begin{equation}
Lw=\lam\left(\mathrm{Re}\frac{A_0}{\lam}-i\nu\right)w+C(t)\ov{w}
\end{equation}
where $C(t)=\dis B(t)\frac{\ov{m}(t)}{m(t)}$.}

\vspace{.2cm}

 In view of this lemma, from now on, we will assume
that $\dis\mathrm{Re}\frac{A_0}{\lam}=0$ and deal with the
simplified equation
 \begeq
Lu=-i\lam\nu u+c(t)\ov{u} \stopeq
 where $\nu\in [0, 1)$ and $c(t)\in C^k(\cir,\C)$.

 Consider the family of vector fields
 \begeq
L_\ep =\lam_\ep\dd{}{t}-ir\dd{}{r}
 \stopeq
 where $\lam_\ep=a+ib\ep $, $\ \ep\in\R$, and the associated
 operators $\clep$ defined by
 \begeq
\clep u(r,t)=\lam_\ep\dd{u}{t}(r,t)-ir\dd{u}{r}(r,t)+i\lam_\ep\nu
u(r,t)-c(t)\ov{u(r,t)}
 \stopeq
For $\C$-valued functions defined on an open set $U\in \R^+\times
\cir$, we define the bilinear form
\[
<f,g>=\mathrm{Re}\left(\iint_Uf(r,t)g(r,t)\frac{drdt}{r}\right)\, .
\]
 For the duality induced by this form, the adjoint
 of $\clep$ is
 \begeq
 \clepstar v(r,t)=-\left(\lam_\ep\dd{v}{t}(r,t)-ir\dd{v}{r}(r,t)-i\lam_\ep\nu
v(r,t)+\ov{c(t)}\,\ov{v(r,t)}\right)
 \stopeq
 The function $z_\ep(r,t)=|r|^{\lam_\ep}\ei{it}$ is a first integral
 of $L_\ep$ in $\R^\ast\times\cir$. That is,
  $L_\ep z_\ep =0$, $dz_\ep\ne 0$.
 Furthermore $z_\ep :\, R^+\times \cir\,\longrightarrow\, \C^\ast$ is a
 diffeomorphism. The following Green's identity will be used
 throughout.

 {\prop Let $U\subset\R^+\times\cir$ be an open set with piecewise smooth
 boundary. Let $u,v\in C^0(\ov{U})$ with $L_\ep u$ and $L_\ep v$
 integrable. Then,
 \begeq
\mathrm{Re}\left(\int_{\pa U}uv\frac{dz_\ep}{z_\ep}\right)
=<u,\clepstar v>-<\clep u,v>.
 \stopeq}

\vspace{.2cm}

 \Proof{ Note that for a differentiable function
$f(r,t)$, we have
\[
df=\frac{i}{2a}\left( -\ov{L_\ep}f\frac{dz_\ep}{z_\ep}+L_\ep
f\frac{d\ov{z_\ep}}{\ov{z_\ep}}\right) \quad\mathrm{and}\quad
\frac{d\ov{z_\ep}}{\ov{z_\ep}}\wedge\frac{dz_\ep}{z_\ep}=\frac{2ia}{r}dr\wedge
dt\, .
\]
Hence,
 \[ \begar{ll}
  \dis \int_{\pa U}uv\frac{dz_\ep}{z_\ep} & =\dis \iint_U
  \frac{i}{2a}\left(uL_\ep v+vL_\ep u \right)
  \frac{d\ov{z_\ep}}{\ov{z_\ep}}\wedge\frac{dz_\ep}{z_\ep}\\
  & =-\dis \iint_U\left(v\clep u -u\clepstar v
  +cv\ov{u}-u\ov{cv}\right)\frac{dr dt}{r}  \stopar
 \]
 By taking the real parts, we get $(1.8)$
 }

\Rem{\bf 1.1} { When $b=0$ so that $\lam =a\in\R^+$. The pushforward
via the first integral $r^a\ei{it}$ reduces the equation $\cl u=F$
into a Cauchy Riemann equation with a singular point of the form
 \begeq
 \dd{W}{\ov{z}} =\frac{a_0}{\ov{z}}W+\frac{B(t)}{\ov{z}}\ov{W}
 +G(z)
 \stopeq
 Properties of the solutions of such  equations are
  thoroughly studied in {\cite{Mez-CV3}}. Many aspects of CR
  equations with punctual singularities have been studied by
 a number of
  authors and we would like to mention in particular the following
 papers {\cite{Beg-Dao}}, {\cite{Tung-1}}, {\cite{Tung-2}},
 {\cite{Usm1}} and {\cite{Usm2}}}.

\Rem{\bf 1.2} { We should point out that the vector fields involved
here satisfy the Nirenberg-Treves Condition (P)  at each point of
the characteristic circle. For vector fields $X$ satisfying
condition (P), there is a rich history for the local solvability of
the $\C$-linear  equation $Xu=F$ (see the books
{\cite{Ber-Cor-Hou}}, {\cite{Trev-Book}} and the references
therein). Our focus here is first, on the semiglobal solvability in
a tubular neighborhood of the characteristic circle, and second, on
the equations containing the term $\ov{u}$ which makes them not
$\C$-linear. }

 \Rem{\bf 1.3} {The operator $\clep$ is invariant under
 the diffeomorphism $\Phi(r,t)=(-r,t)$ from
 $\R^+\times\cir$ to $\R^-\times\cir$. Hence, all the results
 about $\clep$ stated in domains contained in $\R^+\times\cir$
 have their counterparts for domains in $\R^-\times \cir$.
 Throughout this paper, we will be mainly stating results for
 $r\ge 0$.}

\section{Basic Solutions}
 In this section we introduce the notion of basic solutions for
 $\clep$.
 We say that $w$ is a basic solution of $\clep$ if it is a nontrivial
solution of $\clep w=0$, in $\R^+\times \cir$, of the form
 \begeq w(r,t)=r^\sigma
\phi(t)+\ov{r^\sigma\psi(t)},
 \stopeq
 with $\sigma\in\C$ and where $\phi(t)
,\psi(t)\ $ are $2\pi$-periodic functions. These solutions play a
crucial role for the equations generated by $L_\ep$. In a sense,
they play  roles similar to those played by the functions $z^n$ in
classical complex and harmonic analysis.

Consider, as our starting point, the basic solutions of $\cl_0$.
These basic solutions are known, since they can be recovered from
those of equation (1.9) (see Remark 1.1). From $\cl_0$, we obtain
the properties of the basic solution for $\clep$. This is done
through continuity arguments in the study of an associated system of
$2\times 2$ ordinary differential equations in $\C^2$ with periodic
coefficients. By using analytic dependence of the system with
respect to the parameters, the spectral values $\sigma$ of the
monodromy matrix can be tracked down. The main result (Theorem 2.1)
states that for every $j\in\Z$, the operator $\clep$ has exactly two
$\R$-independent basic solutions with winding number $j$.

 \subsection{Properties of basic solutions}
We prove that a basic solution has no vanishing points when $r>0$
and that one of its components $\phi$ or $\psi$ is always
dominating.

 It is immediate, from (1.6), that in order for a function $w(r,t)$,
 given by (2.1), to satisfy $\clep w=0$, the components
 $\phi$ and $\psi$ need to be periodic solutions of the system of
ordinary differential equations
 \begeq\left\{\begar{l}
 \lam_\ep\phi'(t)=i(\sigma-\lam_\ep\nu)\phi(t)+c(t)\psi(t)\\
 \ov{\lam_\ep}\psi'(t)=-i(\sigma-\ov{\lam_\ep}\nu)\psi(t)+\ov{c(t)}\phi(t)\,
 .
\stopar \right.
 \stopeq
 Note that if $\sigma\in\R$, then $w=r^\sigma (\phi(t)+\ov{\psi(t)})$ and
 $f=\phi +\ov{\psi}$ solves the equation
 \begeq
\lam_\ep f'(t)= i(\sigma -\lam_\ep\nu)f(t)+c(t)\ov{f(t)}\, .
 \stopeq
 Now we prove that a basic solution cannot have zeros when $r>0$.

{\prop Let $w(r,t)$, given by $(2.1)$, be a basic solution of
$\clep$. Then
 \[
w(r,t)\ne 0\qquad\forall (r,t)\in\R^+\times\cir\, .\] }

\vspace{.2cm}

\Proof{ If $\sigma\in\R$, we have $w(r,t)=r^\sigma f(t)$ with $f(t)$
satisfying $(2.3)$. If $w(r_0,t_0)=0$ for some $r_0>0$, then
$f(t_0)=0$ and so $f\equiv 0$ by uniqueness of solutions of the
differential equation $(2.3)$. Now, assume that $\sigma
=\alpha+i\beta$ with $\beta \in\R^\ast$. Suppose that $w$ is a basic
solution and $w(r_0,t_0)=0$ for some $(r_0,t_0)\in\R^+\times\cir$.
Consider the sequence of real numbers $\dis r_k=r_0\exp(-k\pi
/|\beta |)$ with $k\in\Z^+$. Then $r_k\longrightarrow 0$ as
$k\longrightarrow\infty$ and $r_k^{2i\beta}=r_0^{2i\beta}$. It
follows
 at once from $w(r_0,t_0)=0$ and $(2.1)$ that
 $w(r_k,t_0)=0$ for every $k\in\Z^+$. Note that  from
 $(2.1)$ we have $|w(r,t)|\le Er^a$, where
 $E=\max(|\phi(t)|+|\psi(t)|)$. Note also that since $\clep$ is
 elliptic in $\R^+\times\cir$, then the zeros of any solution of the
 equation $\clep u=0$ are isolated in  $\R^+\times\cir$.

 The pushforward via the mapping $z=r^{\lam_\ep}\ei{it}$ of the
 equation $\clep w=0$ in $\R^+\times\cir$ is the singular CR equation
 \[
\dd{W}{\ov{z}}=\frac{\lam_\ep\nu\ei{2i\ta}}{2az}W-
\frac{C(z)\ei{2i\ta}}{2iaz}\ov{W}
 \]
 where $W(z)$ and $C(z)$ are the pushforwards of $w(r,t)$ and $c(t)$ and
 where $\ta$ is the argument of $z$.
  We are going to show that $W$ has the form
 $W(z)=H(z)\exp(S(z))$ where $H$ is holomorphic in the punctured disc
  $D^\ast (0,R)$, $S(z)$ continuous in $D^\ast (0,R)$ and satisfies the growth
  condition $\dis |S(z)|\le \log\frac{K}{|z|^p}$ for some positive constants
  $K$ and $p$. For this, consider the function $M(z)$ defined by
  \[
M(z)=\frac{\lam_\ep\nu\ei{2i\ta}}{2a}-
\frac{C(z)\ei{2i\ta}}{2ia}\frac{\ov{W(z)}}{W(z)}
  \]
 for $0<|z|<R$, $W(z)\ne 0$ and by $M(z)=1$ on the set of isolated points
 where $W(z)=0$. This function is bounded and it follows from the classical
 theory of CR equations (see {\cite{Beg-Book}} or {\cite{Vek}}) that
 \[
N(z)=\frac{-1}{\pi}\iint_{D(0,R)}\frac{M(\zeta)}{\zeta -z}d\xi d\eta
 \]
 ($\zeta =\xi+i\eta$) is continuous, satisfies
 $\dis
\dd{N(z)}{\ov{z}}=M(z)$ and
 \[
 |N(z_1)-N(z_2)|\le A||M||_\infty
 |z_1-z_2|\log\frac{2R}{|z_1-z_2|}\quad
 \forall z_1,z_2\in D(0,R)
 \]
 for some positive constant $A$. Define $S$ by
 $S(z)=\dis\frac{N(z)-N(0)}{z}$. We have then, for $z\ne 0$,
 \[\dd{S}{\ov{z}}
=\frac{W_{\ov{z}}(z)}{W(z)}
 \quad\mathrm{and}\quad |S(z)|\le B\log\frac{2R}{|z|}\, .\]
 with $B=A||M||_\infty$
 Let $H(z)=W(z)\exp(-S(z))$. Then $H$ is holomorphic in $0<|z|<R$ and it
 satisfies
 \[
|H(z)|\le |W(z)|\exp(|S(z)|)\le |W(z)|\frac{(2R)^B}{|Z|^B}\le
C_1|z|^s\, .
 \]
 for some constants $C_1$ and $s\in R$. The last inequality follows
 from the estimate $|w|\le Er^\alpha$. This means that the function
 $H$ has at most a pole at $z=0$. Since $w(r_k,t_0)=0$, then
 $H(z_k)=0$ for every $k$ and $z_k=r_k^{\lam_\ep}\ei{it_0}\longrightarrow
 0$. Hence $H\equiv 0$ and $w\equiv 0$ which is a contradiction}

{ \cor{ If $w=r^\sigma\phi(t)+\ov{r^\sigma\psi(t)}$ is a basic
solution
 of $\clep$ with $\sigma =\alpha +i\beta$ and $\beta\ne 0$,
 then for every $t\in\R$,$\ |\phi(t)|\ne |\psi(t)|$ }}

\vspace{.2cm}

 \Proof{ By contradiction, suppose that there is
$t_0\in\R$ such that $|\phi(t_0)|=|\psi(t_0)|$. Let $x_0\in\R$ such
that $\ov{\psi(t_0)}=-\ei{ix_0}\phi(t_0)$.
 Then the positive number $r_0=\exp(x_0/2\beta)$ satisfies
 $r_0^{i\beta}=r_0^{-i\beta}\ei{ix_0}$ and consequently,
 \[
w(r_0,t_0)=r_0^{\alpha}(r_0^{i\beta}\phi(t_0)+\ov{r_0^{i\beta}\psi(t_0)})=0.
 \]
 This contradicts Proposition 2.1. }

This corollary implies that, for a given basic solution
$w=r^\sigma\phi +\ov{r^\sigma\psi}$ with $\sigma\in\C\backslash\R$,
one of the functions $\phi$ or $\psi$ is dominant. That is, $|\phi
(t)|>|\psi(t)|$ or $|\psi(t)|>|\phi(t)|$ for every $t\in\R$. Hence
the winding number of $w$, $\ \ind (w)$ is well defined and we have
$\ind (w)=\ind (\phi)$ if $|\phi|>|\psi|$ and $\ind (w)=\ind
(\ov{\psi})$ otherwise. When $\sigma\in\R$, we have $w=r^\sigma
f(t)$ with $f$ nowhere 0 and so $\ind (w)=\ind (f)$.

For a basic solution $w=r^\sigma\phi +\ov{r^\sigma\psi}$ with
$|\phi|>|\psi|$, we  will refer to $\sigma$ as the exponent of $w$
(or a spectral value of $\clep$) and define the character of $w$ by
\[
\charac (w)=(\sigma ,\ind (w))
\]
We will denote by $\speclep$ the set of exponents of basic
solutions. That is
 \begeq
\speclep =\{ \sigma\in\C ;\, \exists w,\ \charac (w)=(\sigma ,\ind
(w))\}
 \stopeq

\Rem{\bf 2.1} {When $\sigma\in\C\backslash\R$ and $w=r^\sigma\phi(t)
+\ov{r^\sigma\psi(t)}$ is a basic solution with $\charac (w)=(\sigma
,\ind (\phi))$, the function $\widetilde{w}=r^\sigma (i\phi(t))+
\ov{r^\sigma i\psi(t)}$ is also a basic solution with $\charac
(w)=\charac (\widetilde{w})$ and $w$, $\widetilde{w}$ are
$\R$-independent.

When $\sigma =\tau \in \R$, and $w=r^\tau f(t)$ is a basic solution
with $\charac (w)=(\tau ,\ind (f))$, it is not always the case that
there is a second $\R$-independent basic solution with the same
exponent $\tau$. There is however a second $\R$-independent basic
solution $\widetilde{w}=r^{\tau'}g(t)$ with the same winding number
($\ind (f)=\ind (g)$) but with a different exponent $\tau'$ (see
Proposition 2.3).}

The following proposition follows from the constancy of the winding
number under continuous deformations.

{\prop Let $w_\ep (r,t)=r^{\sigma(\ep)}\phi (t,\ep)+
\ov{r^{\sigma(\ep)}\psi (t,\ep)}$ be a continuous family of basic
solutions of $\clep$ with $\ep\in I$, where $I\subset \R$ is an
interval. Then $\charac (w_\ep)$ depends continuously on $\ep$ and
$\ind (w_\ep )$ is constant.}

\subsection{The spectral  equation and $\mathrm{Spec}(\cl_0)$}
 We use the $2\times 2$ system of ordinary differential equations to obtain an
 equation for the spectral values in terms of the monodromy matrix.
 Results about the CR equation (1.9) are then used to list the
 properties of $\mathrm{Spec}(\cl_0)$.

In order for a function
$w(r,t)=r^\sigma\phi(t)+\ov{r^\sigma\psi(t)}$ to be a basic solution
of $\clep$, the $2\pi$-periodic and $\C^2$-valued function
$V(t)=\dis\left(\begar{l}\phi(t)\\
\psi(t)\stopar\right)$  must solve the periodic system of
differential equations
\[
\dot{V}=\MM (t,\sigma,\ep)V\tag{$\E$}
\]
where
\[
\MM (t,\sigma ,\ep)=\left(\begar{ll} \dis i\frac{\sigma
-\lam_\ep\nu}{\lam_\ep} & \dis\frac{c(t)}{\lam_\ep}\\
\dis\frac{\ov{c(t)}}{\ov{\lam_\ep}} & -i\dis\frac{\sigma
-\ov{\lam_\ep}\nu}{\ov{\lam_\ep}} \stopar\right)
\]
Note that since $\MM$ is linear in $\sigma\in\C$ and real analytic
in $\ep\in\R$, then any solution $V(t,\sigma ,\ep)$ is an entire
function in $\sigma$ and real analytic in $\ep$. The fundamental
matrix of $(\E)$ is the $2\times 2$ matrix $\VV (t,\sigma ,\ep)$
satisfying
\[
\dot\VV =\MM (t,\sigma,\ep)\VV,\qquad \VV (0,\sigma,\ep)=\II
\]
where $\II$ is the identity matrix. We know from Floquet theory that
\[
\VV (t,\sigma,\ep)=\PP (t,\sigma ,\ep)\exp(t\KK(\sigma ,\ep))
\]
where $\PP$ is a $2\pi$-periodic matrix (in $t$) and $\PP$ and $\KK$
are entire in $\sigma$ and real analytic in $\ep$. The monodromy
matrix of ($\E$) is
\[
\BB (\sigma,\ep)=\VV (2\pi ,\sigma,\ep)=\exp(2\pi\KK(\sigma,\ep))\,
.
\]
The Liouville-Jacobi formula gives
\[
\det(\VV(t,\sigma,\ep))=\exp\left(\int_0^t\tr(\MM(s,\sigma,\ep))ds\right)
=\exp\left(\frac{2b\ep}{|\lam_\ep|^2}\sigma t\right),
\]
where $\det(A)$ and $\tr(A)$ denote the determinant and the trace of
the matrix $A$. Hence,
 \begeq
\det(\BB(\sigma,\ep))=\exp\left(\frac{4\pi b\ep}{|\lam_\ep|^2}\sigma
\right)\, . \stopeq
 In order for system ($\E$) to have a periodic
solution, the corresponding monodromy matrix $\BB$ must have 1 as an
eigenvalue. Thus $\sigma$ must solve the spectral equation
\[
1-\tr(\BB(\sigma,\ep))+\det(\BB(\sigma,\ep))=0.
\]
or equivalently, $F(\sigma,\ep)=0$, where
 \begeq
F(\sigma,\ep)=\tr(\BB(\sigma,\ep))-1-\exp\left(\frac{4\pi
b\ep}{|\lam_\ep|^2}\sigma \right)\, .
 \stopeq
We first verify that $\speclep$ is a discrete set.

{\lem For every $\ep\in\R$,   $\ \speclep$ is a discrete subset of
$\C$.}

\vspace{.2cm}

 \Proof{ By contradiction, suppose that there exists
$\ep_0\in\R$ such that $\specb(\cl_{\ep_0})$ has an accumulation
point in $\C$. This means that the roots of the solutions of the
spectral equation $F(\sigma,\ep_0)=0$ have an accumulation point.
Since $F$ is an entire function, then $F(\sigma,\ep_0)\equiv 0$.
Thus, $\specb(\cl_{\ep_0})=\C$. Let
 $\dis\left(\begar{l}\phi(t,\sigma)\\ \psi(t,\sigma)\stopar\right)$
 be a continuous family of periodic solutions of
 $(\mathrm{E}_{\sigma,\ep_0})$. By Proposition 2.2, we can assume
 that $|\phi|>|\psi|$ for every $\sigma\in \R^++i\R$ and that
 $\ind(\phi)=j_0$ (is constant). Now the first equation of
 $(\mathrm{E}_{\sigma,\ep_0})$ gives
 \[
 \frac{\lam_{\ep_0}}{2\pi i}\int_0^{2\pi}
 \frac{\dot\phi(t,\sigma)}{\phi(t,\sigma)}dt =
 \sigma -\lam_{\ep_0}\nu +\frac{1}{2\pi i}\int_0^{2\pi}
 c(t)\frac{\psi(t,\sigma)}{\phi(t,\sigma)}dt\, .
 \]
 That is,
 \[
 \sigma =\lam_{\ep_0}(j_0+\nu) -\frac{1}{2\pi i}\int_0^{2\pi}
 c(t)\frac{\psi(t,\sigma)}{\phi(t,\sigma)}dt\qquad
 \forall \sigma\in \R^++i\R\, .
 \]
 This is a contradiction since
 $\dis\left|c\frac{\psi}{\phi}\right|<|c|$
 }

 The following proposition describes the spectrum of $\cl_0$.

 {\prop
For every $j\in\Z$, there exist $\tau_j^\pm\in\R$ with
$\tau_j^-\le\tau_j^+$ and $f_j^\pm\in C^{k+1}(\cir,\C)$, such that
$w_j^\pm(r,t)=r^{\tau_j^\pm}f_j^\pm(t)$ are $\R$-independent basic
solution of $\cl_0$ with
\[
\charac(w_j^\pm)=(\tau_j^\pm ,j)\, .
\]
Furthermore, $\mathrm{Spec}(\cl_0)=\{\tau_j^\pm ,\  j\in\Z\}$,
\[
\cdots <\tau_{-1}^-\le\tau_{-1}^+<\tau_0^-\le\tau_0^+<\tau_{1}^-\le
\tau_1^+<\cdots
\]
with $\dis\lim_{j\to -\infty}\tau_j^\pm =-\infty,\ \
 \lim_{j\to \infty}\tau_j^\pm =\infty$}

\vspace{.2cm}

 \Proof{ We have here $\lam_0=a >0$. The pushforward
of the equation $\cl_0 w=0$
 via the first integral $z=r^a\ei{it}$ of
 $L_0$ gives a CR equation with
 a singularity of the form studied in {\cite{Mez-CV3}}.
The spectral values $\tau_j^\pm$ of the CR equations are as in the
proposition. It remains only to verify that $\cl_0$ (or its
equivalent CR equation) has no complex spectral values. The Laurent
series representation for solutions of the CR equation (see
{\cite{Mez-CV3}}) imply that any solution of $\cl_0 w=0$ can be
written as
\[
w(r,t)=\sum_{j\in\Z}c_j^-r^{\tau_j^-}f_j^-(t)+c_j^+r^{\tau_j^+}f_j^+(t)
\]
with $c_j^\pm\in\R$. Now, if
$w=r^\sigma\phi(t)+\ov{r^\sigma\psi(t)}$ is a basic solution of
$\cl_0$, then it follows at once, from the series representation,
that $\sigma$ is one of the $\tau_j^\pm$'s}

\subsection{Existence of basic solutions}
 We use the spectral equation together with Proposition 2.3 to show
 the existence of basic solutions  for
 $\clep$ with any given winding number.
 More precisely, we have the following proposition.

{\prop{ For every $j\in\Z$, there exists
$\sigma_j^\pm(\ep)\in\speclep$ such that $\sigma_j^\pm (\ep)$
depends continuously on $\ep\in\R$, $\sigma_j^\pm(0)=\tau_j^\pm$,
and the corresponding basic solution
\[
w_j^\pm(r,t,\ep)=r^{\sigma_j^\pm(\ep)}\phi(t,\ep)+
\ov{r^{\sigma_j^\pm(\ep)}\psi(t,\ep)}
\]
is continuous in $\ep$ and
$\charac(w_j^\pm)=(\sigma_j^\pm(\ep),j)$.}}

\vspace{.2cm}

 \Proof{ For a given $j\in\Z$, it follows from Proposition
2.3 that the monodromy matrix $\BB (\tau_j^\pm ,0)$ admits 1 as an
eigenvalue.
 Since the spectral function $F(\sigma,\ep)$ given
by (2.6) is entire in $\C\times\R$ and since $F(\tau_j^\pm ,0)=0$,
then $F(\sigma,\ep)=0$ defines an analytic variety $\mathcal{V}$ in
$\C\times\R$ passing through the points $(\tau_j^\pm, 0)$. The
variable $\ep$ can be taken as a parameter for a branch of
$\mathcal{V}$ through the point $(\tau_j^\pm,0)$. This means that
the equation $F(\sigma,\ep)=0$ has a solution $\sigma =g(\ep)\in\C$,
with $g$ continuous and $g(0)=\tau_j^\pm$. In fact, $g$ is real
analytic except at isolated points.  The matrix $\BB(g(\ep),\ep)$ is
continuous and has 1 as an eigenvalue for every $\ep$. Let $E_0^\pm$
be an eigenvector of $\BB (\tau_j^\pm ,0)$ with eigenvalue 1. We can
select a continuous vector $E^\pm (\ep)\in\C^2$ such that
\[
\BB (g(\ep),\ep)E^\pm(\ep )=E^\pm(\ep)\quad \mathrm{and}\quad
E^\pm(0)=E_0^\pm
\]
Let $V(t,\ep)=\VV (t,g(\ep),\ep)E^\pm(\ep)$. Then $V(t,\ep)$ is a
periodic solution of the equation ($\mathrm{E}_{g(\ep),\ep}$). If we
set $V(t,\ep)=\dis
\left(\begar{l}\phi_j^\pm(t,\ep)\\
\psi_j^\pm(t,\ep)\stopar\right)$, then
\[
w(r,t,\ep)=r^{g(\ep)}\phi_j^\pm(t,\ep)+\ov{r^{g(\ep)}\psi_j^\pm(t,\ep)}
\]
is a basic solution of $\clep$ and it depends continuously on $\ep$.
Since for $\ep =0$, $\ w(r,t,0)$ has character $(\tau_j^\pm,j)$,
then by Proposition 2.2, the character of $w(r,t,\ep)$ is either
$(g(\ep),j)$ if $|\phi_j^\pm|>|\psi_j^\pm|$ or $(\ov{g(\ep)},j)$ if
$|\phi_j^\pm|<|\psi_j^\pm|$. In the first case,
$\sigma_j^\pm(\ep)=g(\ep)\in\speclep$ and, in the second,
$\sigma_j^\pm(\ep)=\ov{g(\ep)}\in\speclep$}

\subsection{Properties of the fundamental matrix of
($E_{\sigma,\ep}$)}
 We prove some symmetry properties of the fundamental matrix and of
 the monodromy matrix that will be used shortly.

{\prop There exist functions $f(t,\sigma,s),\ g(t,\sigma,s)$ of
class $C^{k+1}$ in $t\in\R$, analytic in $(\sigma, s)\in\C\times\R$,
such that the fundamental matrix $\VV(t,\sigma,\epsilon)$ of $(\E)$
has the form
 \begeq \VV(t,\sigma,\ep)=\left(\begar{ll}
f(t,\sigma,\ep^2)
&\ov{\lam_\ep g(t,\ov{\sigma},\ep^2)}\\
\lam_\ep g(t,\sigma,\ep^2) &
\ov{f(t,\ov{\sigma},\ep^2)}\stopar\right)
 \exp\left(\frac{\ep bt}{|\lam_\ep|^2}\sigma\right).
 \stopeq
Furthermore, $f$ and $g$ satisfy
 \begeq
f(t,\sigma,\ep^2)\ov{f(t,\ov{\sigma},\ep^2)}- |\lam_\ep|^2
g(t,\sigma,\ep^2)\ov{g(t,\ov{\sigma},\ep^2)}\equiv 1
 \stopeq}

\vspace{.2cm}

\Proof{ If we use the substitution $\dis V=Z\exp\left(\frac{\ep
b\sigma t}{|\lam_\ep|^2}\right)$ in equation ($\E$), then the system
for $Z$ is
 \begeq \dot{Z}=\Abf (t,\sigma,\ep^2)Z \stopeq
 with
 \[
\Abf(t,\sigma ,\ep^2)=\left(\begar{ll} i\mu & \ \ov{\lam_\ep
d(t,\ep^2)}\\ \lam_\ep d(t,\ep^2) & \ -i\mu\stopar\right)
 \]
 and where
 \[ \mu =\frac{a\sigma}{a^2+b^2\ep^2}-\nu\, ,
 \quad\mathrm{and}\quad d(t,\ep^2)=\frac{\ov{c(t)}}{a^2+b^2\ep^2}\]
The fundamental matrix $\ZZ (t,\sigma,\ep^2)$ of (2.9) with $\ZZ
(0,t,\ep^2)=\II$ is therefore of class $C^{k+1}$ in $t$ and analytic
in  $(\sigma ,s)$ with $s=\ep^2$.

For  functions $F(t,\mu,\ep^2)$ and $G(t,\mu,\ep^2)$, with
$\mu\in\R$, we use the notation
\[
\DD_F=\left(\begar{ll}F & 0\\ 0 &\ov{F}\stopar\right)
\quad\mathrm{and}\quad
 \JJ_{\lam_\ep G}=\left(\begar{ll} 0 & \ov{\lam_\ep G}\\ \lam_\ep G
& 0\stopar\right)
\]
With this notation, system (2.9) has the form
 \begeq
\dot{Z} =\left(\DD_{i\mu}+\JJ_{\lam_\ep d}\right)Z
 \stopeq
 Note that we have the following relations
 \[\begar{ll}
\DD_{F_1}\DD_{F_2}=\DD_{F_1F_2}, & \DD_{F}\JJ_{\lam_\ep
G}=\JJ_{\lam_\ep\ov{F}G},\\
 \JJ_{\lam_\ep G}\DD_{F}=\JJ_{\lam_\ep
FG}, & \JJ_{\lam_\ep G_1}\JJ_{\lam_\ep
G_2}=\DD_{|\lam_\ep|^2\ov{G_1}G_2} \stopar \]
 The fundamental matrix $\ZZ$ is obtained as the limit,
 $\ZZ =\lim_{k\to\infty}\ZZ_k$, where the matrices $\ZZ_k(t,\sigma
 ,\ep)$ are defined inductively by $\ZZ_0=\II$ and
 \[
\ZZ_{k+1}(t,\sigma,\ep)=\II +\int_0^t\left(\DD_{i\mu}+\JJ_{\lam_\ep
d}\right)\ZZ_k(s,\sigma ,\ep)ds
 \]
 Now, we prove by induction that $\ZZ_k=\DD_{F_k}+\JJ_{\lam_\ep
 G_k}$, where $F_k$ and $G_k$ are polynomials in the variable $\mu$,
 analytic and even in the variable $\ep$.
 The claim is obviously true for
 $\ZZ_0=\II$. Suppose that $\ZZ_k$ has the desired property for
 $k=0,\cdots ,n$, then
 \[\begar{ll}
\ZZ_{n+1}(t,\sigma ,\ep) & =\dis\II +\int_0^t
\left(\DD_{i\mu}+\JJ_{\lam_\ep
d}\right)\left(\DD_{F_n}+\JJ_{\lam_\ep G_n}\right)ds\\
& =\DD_{F_{n+1}}+\JJ_{\lam_\ep G_{n+1}} \stopar \]
 where
 \[\begar{ll}
 F_{n+1} (t,\mu ,\ep^2) & =\dis 1+\int_0^t \left( i\mu F_n(s,\mu
 ,\ep^2)+c(s)G_n(s,\mu ,\ep^2)\right) ds\\
 G_{n+1} (t,\mu ,\ep^2) & =\dis -\int_0^t \left( i\mu G_n(s,\mu
 ,\ep^2)+\ov{c(s)}F_n(s,\mu ,\ep^2)\right) ds
 \stopar \]
By taking the limit as $k\to\infty$, we get the fundamental matrix
\begeq \ZZ(t,\mu ,\ep)=\DD_{F_0(t,\mu,\ep^2)}+\JJ_{\lam_\ep
G_0(t,\mu,\ep^2)} \stopeq
 where $F_0$ and $G_0$ are entire functions
with respect to the real parameters $\mu$ and $\ep^2$. The
fundamental matrix $\ZZ (t,\sigma,\ep)$ of (2.9)  is therefore
$\DD_{f(t,\sigma,\ep^2)}+\JJ_{\lam_\ep g(t,\sigma ,\ep^2)}$, where
$f$ and $g$ are the holomorphic extensions of $F_0$ and $G_0$
(obtained by replacing $\mu$ by
$\dis\frac{a\sigma}{a^2+b^2\ep^2}-\nu$, with $\sigma\in\C$). The
proposition follows immediately, since $\VV
=\dis\ZZ\exp\left(\frac{\ep b\sigma t}{|\lam_\ep|^2}\right)$ }

A direct consequence of expression (2.7) of the proposition is the
following

{\prop Let $\VV(t,\sigma,\ep)$ be the fundamental matrix of $(\E)$,
then the fundamental matrix of equation
$(\mathrm{E}_{\ov{\sigma},\ep})$ is
 \begeq
 \VV(t,\ov{\sigma},\ep)=\JJ\ov{\VV(t,\sigma,\ep)}\JJ
 \stopeq
 where $\JJ =\dis\left(\begar{ll} 0 & 1\\ 1& 0\stopar \right)$.}

\vspace{.2cm}

The monodromy matrix of ($\E$) has the form
 \begeq
 \BB(\sigma,\ep)=\left(\begar{ll}p(\sigma,\ep^2) &\ov{\lam_\ep
 q(\ov{\sigma},\ep^2)}\\
  \lam_\ep q(\sigma,\ep^2) &\ov{p(\ov{\sigma},\ep^2)}\stopar\right)
  \exp\left(\frac{2\pi\ep b}{a^2+b^2\ep^2}\sigma\right)
 \stopeq
 with $p(\sigma ,\ep^2)=f(2\pi ,\sigma,\ep^2)$,
  $\ q(\sigma ,\ep^2)=g(2\pi ,\sigma,\ep^2)$ satisfying
  \begeq
p(\sigma ,\ep^2)\ov{p(\ov{\sigma} ,\ep^2)}- |\lam_\ep|^2
 q(\sigma,\ep^2)\ov{q(\ov{\sigma} ,\ep^2)}\equiv 1
  \stopeq
  Note that
  \begeq
\det(\BB(\sigma,\ep))=\exp\left(\frac{4\pi\ep
b}{a^2+b^2\ep^2}\sigma\right)
  \stopeq
We denote by $\specb(\BB (.,\ep))$ the set of spectral values of
  $\BB (.,\ep)$, i.e.
  \[
\specb(\BB (.,\ep))=\{\sigma\in\C ;\ \det(\BB(\sigma,\ep)-\II)=0\}.
  \]
It follows from (2.13) and Proposition 2.6 that
 \begeq
\specb(\BB (.,\ep))=\specb(\BB (.,-\ep))=\ov{\specb(\BB (.,\ep))}
\stopeq
 An element $\sigma\in\specb(\BB (.,\ep))$ is said to be a simple (or a double)
 spectral value if
 the corresponding eigenspace has dimension 1 (or 2).

 The spectral
function $F$ defined in (2.6) takes form
  \begeq
  F(\sigma ,\ep^2)=p(\sigma ,\ep^2)+\ov{p(\ov{\sigma},\ep^2)}-
  2\cosh\left(\frac{2\pi\ep b}{a^2+b^2\ep^2}\sigma\right)
  \stopeq
 Thus if $F(\sigma,\ep^2)=0$, then $F(\ov{\sigma},\ep)=0$. We have
 therefore
 {\cor{ If $\sigma\in\speclep$, then $\sigma\ \mathrm{or}\ \ov{\sigma}\,
  \in\mathrm{Spec}(\mathcal{L}_{-\ep})$.}}

\subsection{The system of equations for the adjoint operator $\clepstar$}
 The properties of the fundamental matrix of system $\E$ will be
 used to obtain those for the adjoint operator.
 The system of ordinary differential equations for
 the adjoint operator $\clepstar$ given in (1.7) is
\[
\dot{V}=\widetilde{\MM}(t,\mu,\ep)V\tag{$\WE$}
\]
where
 \begeq
 \widetilde{\MM}(t,\mu,\ep)=\left(\begar{ll}
 i\dis\frac{\mu +\lam_\ep\nu}{\lam_\ep} &\
 -\dis\frac{\ov{c(t)}}{\lam_\ep}\\
 -\dis\frac{c(t)}{\ov{\lam_\ep}} &\
 -i\dis\frac{\mu +\ov{\lam_\ep}\nu}{\ov{\lam_\ep}}
 \stopar\right)
 \stopeq
 Thus, if $\dis V(t)=\left(\begar{l} X(t)\\ Z(t)\stopar\right)$
 is a periodic solution of ($\WE$), then
 $w(r,t)=r^\mu X(t)+\ov{r^\mu Z(t)}$ is a basic solution
of $\clepstar$.

The relation between the fundamental matrices of this system and
those for $\E$ is given by the following proposition.

{\prop The fundamental matrix of $(\WE )$ is
 \begeq
 \widetilde{\VV}(t,\mu ,\ep)=\DD\ov{\VV(t,-\ov{\mu},-\ep)}\DD
 \stopeq
 where $\VV(t,-\ov{\mu},-\ep)$ is the fundamental matrix
 of $(\mathrm{E}_{-\ov{\mu},-\ep})$ and
 $\DD=\dis\left(\begar{ll}1 & 0\\ 0 & -1\stopar\right)$.}

\vspace{.2cm}

 \Proof{ If $\dis V=\left(\begar{l} X\\ Z\stopar\right)$ solves
($\WE$), then $\dis \ov{\DD V}=\left(\begar{l}\ \ov{X}\\ -\ov{
Z}\stopar\right)$ solves the equation
$(\mathrm{E}_{-\ov{\mu},-\ep})$. Therefore, if
 $\dis \left(\begar{ll} X_1 & X_2\\  Z_1 & Z_2\stopar\right)$
 is a fundamental matrix of $(\WE)$, then
 $\dis \left(\begar{ll} \ \ov{X_1} & -\ov{X_2}\\
  -\ov{ Z_1} &\ \ov{ Z_2}\stopar\right)$ is a fundamental matrix
  of $(\mathrm{E}_{-\ov{\mu},-\ep})$ }

 Immediate consequences are the following corollaries.

 {\cor The monodromy matrix of $(\WE)$ is
 \begeq
 \widetilde{\BB}(\mu,\ep)=\DD\ov{\BB(-\ov{\mu},-\ep)}\DD
 \stopeq
where $\BB(\sigma,\ep)$ is the monodromy matrix of $(\E)$.
Furthermore, if $\sigma\in\specb(\BB(.,\ep))$ and $\BB(\sigma
,\ep)E=E$, then $-\ov{\sigma}\in\specb(\widetilde{\BB}(.,-\ep))$ and
$\widetilde{\BB}(-\ov{\sigma},-\ep)\DD\ov{E}=\DD\ov{E}$.}

{\cor If $\sigma\in\speclep$, then either
$-\sigma\in\specb(\clepstar)$ or $-\ov{\sigma}\in\specb(\clepstar)$}

\subsection{Continuation of a simple spectral value}
 We start from a simple spectral value, when $\ep=0$, and use
 the properties of the fundamental matrix to obtain the behavior of
 $\sig(\ep)$ for $\ep $ near 0.

{\prop Suppose that $\tau\in\specb(\BB(.,0))$ and that $\tau$ is
simple. Then there exist $\delta >0$ and a unique function
$\sigma\in C^0([-\delta ,\ \delta],\R)$ such that $\sigma(0)=\tau$
and $\sigma(\ep)\in\specb(\BB(.,\ep))$ for every $\ep\in [-\delta ,\
\delta ]$}

\vspace{.2cm}

 \Proof{ The matrix $\BB(\tau ,0)$ has a
single eigenvector $U$ (up to a multiple) with eigenvalue 1. Since
$\det(\BB(\tau,0))=1$ (see (2.13) and (2.14)), then $\BB(\tau,0)$ is
similar to the matrix $\dis\left(\begar{ll} 1 & 1\\ 0 &
1\stopar\right)$. Let $\VV(t,\sigma ,\ep)$ be the fundamental matrix
of $(\E)$. The function
\[
\left(\begar{l}\phi(t)\\ \psi(t)\stopar\right) =\VV(t,\tau ,0)U,
\]
generates all periodic solutions of $(\mathrm{E}_{\tau,0})$. First,
we show that we can find a generator of the form
$\dis\left(\begar{l}f(t)\\ \ov{f}(t)\stopar\right)$ for some
function $f$. For this, note that since $\lam_0=a\in\R$, then it
follows from (2.2) that $\dis\left(\begar{l}\ov{\psi}(t)\\
\ov{\phi}(t)\stopar\right)$ is also a periodic solution of
$(\mathrm{E}_{\tau,0})$. Hence, there exists $c\in\C$, $\ |c|=1$
such that $\ov{\psi}(t)=c\phi (t)$ and $\ov{\phi}(t)=c\psi(t)$. If
$c=1$, then we can take $f=\phi$, if $c\ne -1$, we can take $f=\phi
+\ov{\psi}$, and if $c=-1$, we take $f=i\phi$. The vector
$U_0=\dis\left(\begar{l}f(0)\\ \ov{f}(0)\stopar\right)$ is the
eigenvector of $\BB(\tau ,0)$ that generates the solution
$\dis\left(\begar{l}f(t)\\ \ov{f}(t)\stopar\right)$.

We know from Proposition 2.4 that the spectral function $F(\sigma
,\ep^2)$ given in (2.6) has a root $\sigma(\ep)$ with
$\sigma(0)=\tau$. Furthermore,  $\sigma (\ep)$ is real analytic in a
neighborhood of $0\in\R$, except possibly at $\ep =0$. Now we show
that there is only one such function in a neighborhood of 0 and that
it is real-valued. Starting from $U_0$, we can find a continuous
vector $U(\ep)\in\C^2$ with $U(0)=U_0$ and such that
$\BB(\sigma(\ep),\ep)U(\ep)=U(\ep)$. The function
\[
\left(\begar{l}\phi(t,\ep)\\
\psi(t,\ep)\stopar\right)=\VV(t,\sigma(\ep),\ep)U(\ep)
\]
is a periodic solution of $(\mathrm{E}_{\sigma(\ep).\ep})$ such that
$\phi(t,0)=f(t)$ and $\psi(t,0)=\ov{f}(t)$. It follows from
Corollary 2.0.4 that $-\tau\in\specb(\widetilde{\BB}(., 0))$ and
$-\ov{\sigma(\ep)}\in\specb(\widetilde{\BB}(., \ep))$. Note that if
$V(t)$ is a periodic solution of $(\mathrm{E}_{\tau,0})$, then
$\DD\ov{V}(t)$ is a periodic solution of
$(\widetilde{\mathrm{E}}_{-\tau,0})$. Thus,
 $\dis\left(\begar{l}\ \ov{f}(t)\\ -f(t)\stopar\right)$ solves
$(\widetilde{\mathrm{E}}_{-\tau,0})$. Let $U_1(\ep)\in\C^2$ be a
continuous eigenvector of $\BB (\sigma(\ep),-\ep)$  such that
$U_1(0)=U_0$. Set $U^\star (\ep)=\DD\ov{U_1(\ep)}$. Then, it follows
from Corollary 2.0.3,  that
\[
\widetilde{\BB}(-\ov{\sigma(\ep)},\ep)U^\star
(\ep)=\DD\ov{\BB(\sigma(\ep),-\ep)}\DD
\DD\ov{U_1(\ep)}=\DD\ov{U_1(\ep)}=U^\star(\ep).
\]
Therefore,
\[
\left(\begar{l}X(t,\ep)\\
Z(t,\ep)\stopar\right)=\widetilde{\VV}(t,-\ov{\sigma(\ep)},\ep)
U^\star(\ep)
\]
is a periodic solution of
$(\widetilde{\mathrm{E}}_{-\ov{\sigma(\ep)}.\ep})$ with
$X(t,0)=\ov{f}(t)$ and $Z(t,0)=-f(t)$.

The corresponding basic solutions of $\clep$ and $\clepstar$ are
respectively,
\[\begar{l}
 w_\ep (r,t)=
 r^{\sigma(\ep)}\phi(t,\ep)+\ov{r^{\sigma(\ep)}\psi(t,\ep)},
 \quad\mathrm{and}\\
 w^\star_\ep (r,t)=
 r^{-\ov{\sigma(\ep)}}X(t,\ep)+r^{-\sigma(\ep)}\ov{Z(t,\ep)}
\stopar\]
 We apply Green's formula (1.8) to the pair $w_\ep$, $w^\star_\ep$ in
 the cylinder $A=[R_1,\ R_2]\times\cir$ (with $0<R_1<R_2$) to get
\[
\mathrm{Re}\left[ \int_{\pa A}w_\ep(r,t)w_\ep^\star (r,t)
\frac{dz_\ep}{z_\ep}\right] =0\, .
\]
That is,
 \begeq \mathrm{Re}\left[\int_0^{2\pi}\left((R_2^{\sigma
-\ov{\sigma}}- R_1^{\sigma -\ov{\sigma}})\phi X+ (R_2^{\ov{\sigma}
-\sigma}- R_1^{\ov{\sigma} -\sigma})\ov{\psi Z}\right)idt\right]
=0\, . \stopeq
 Suppose that $\sigma(\ep)$ is not $\R$-valued in a
neighborhood of 0. Then $\sigma(\ep)=\alpha(\ep)+i\beta(\ep)$ with
$\beta(\ep)>0$ (or $<0$) in a an interval $(0,\ \ep_0)$. If we set
$p=\log R_2$ and $q=\log R_1$, we get
\[
R_2^{\sigma -\ov{\sigma}}- R_1^{\sigma -\ov{\sigma}} =\ei{2i\beta
p}-\ei{2i\beta q}=2i\sin(\beta(p-q))\ei{i\beta(p+q)}
\]
and (2.21) becomes (with $x=p+q$ arbitrary)
 \begeq
\mathrm{Re}\left[\int_0^{2\pi}\left(\ei{i\beta x}\phi(t,\ep)
X(t,\ep)-i\ei{-i\beta x}\ov{\psi(t,\ep) Z(t,\ep)}\right)dt\right]
=0\, ,
 \stopeq
Let $\dis
P(\ep)+iQ(\ep)=\!\!\int_0^{2\pi}\!\!\!\phi(t,\ep)X(t,\ep)dt\,$ and
$\,\dis
R(\ep)+iS(\ep)=\!\!\int_0^{2\pi}\!\!\!\ov{\psi(t,\ep)Z(t,\ep)}dt$
 From (2.22), we have
 \[
\cos(\beta(\ep)x)(P(\ep)-R(\ep))-\sin(\beta(\ep)x)(Q(\ep)+S(\ep))=0,
\qquad \forall x\in\R\, .
 \]
Therefore,
\[
P(\ep)-R(\ep)=0,\quad Q(\ep)+S(\ep)=0,\qquad\forall \ep\in (0,\
\ep_0).
\]
By continuity, we get $P(0)=R(0)$ and $Q(0)=-S(0)$. But,
 \[\begar{ll}
P(0)+iQ(0) &
=\dis\int_0^{2\pi}\phi(t,0)X(t,0)dt=\int_0^{2\pi}|f(t)|^2dt\\
R(0)+iS(0) &
=\dis\int_0^{2\pi}\ov{\psi(t,0)Z(t,0)}dt=-\int_0^{2\pi}|f(t)|^2dt
 \stopar\]
 and it follows from $P(0)=R(0)$ that $\dis\int_0^{2\pi}|f|^2dt=0$.
 This is a contradiction since $f\ne 0$. This means that
 $\sigma(\ep)$ is an $\R$-valued function in a neighborhood of $\ep
 = 0$.

 Now we show that $\sigma(\ep)$ is unique near $\ep =0$. By
 contradiction, suppose that there is another real valued
 solution $\sigma_1(\ep)$, with
 $\sigma(\ep)<\sigma_1(\ep)$ in an interval $(0,\ \ep_0)$,
  and $\sigma(0)=\sigma_1(0)=\tau$.
  Let $\phi(t,\ep)$ and $\psi(t,\ep)$ be as above. Let $U'(\ep)$ be
  an eigenvector (with eigenvalue 1) of $\BB(\sigma_1(\ep),-\ep)$
  such that $U'(0)=-iU_0$, where $U_0$ is the eigenvector used
  above. Let $U_1^\star(\ep)=\DD\ov{U'(\ep)}$. Then
  \[
 \widetilde{\BB}(-\sigma_1(\ep),\ep)U^\star(\ep)=U_1^\star(\ep)
 \quad\mathrm{and}\quad U_1^\star(0)=i\DD U_0\, .
  \]
  To $U^\star(\ep)$ corresponds the $2\pi$-periodic solution
\[
\left(\begar{l}X_1(t,\ep)\\
Z_1(t,\ep)\stopar\right)=\widetilde{\VV}(t,-\ov{\sigma_1(\ep)},\ep)
U_1^\star(\ep)
\]
of $(\widetilde{\mathrm{E}}_{-\ov{\sigma_1(\ep)}.\ep})$ with
$X_1(t,0)=i\,\ov{f}(t)$ and $Z_1(t,0)=-if(t)$. The corresponding
basic solution of $\clepstar$ is
\[
w_{1,\ep}^\star(r,t)=r^{-\sigma_1(\ep)}(X_1(t,\ep)+\ov{Z_1(t,\ep)})
\]
The Green's formula, applied to the pair $w_\ep$, $w_{1,\ep}^\star$
in the cylinder $(R_1,\ R_2)\times\cir$,  gives
\[
\mathrm{Re}\left[\int_0^{2\pi}(R_2^{\sigma -\sigma_1}- R_1^{\sigma
-\sigma_1})(\phi +\ov{\psi})(X_1+\ov{Z_1})idt \right] =0\, .
\]
Thus,
\[
\mathrm{Re}\left[\int_0^{2\pi}(\phi(t,\ep)
+\ov{\psi(t,\ep)})(X_1(t,\ep)+\ov{Z_1(t,\ep)})idt \right] =0 \qquad
\forall \ep\in (0,\ \ep_0)\, .\]
 By letting $\ep\to 0$, we get again
 $\dis\int_0^{2\pi}|f(t)|^2dt=0$,
 which is a contradiction. This shows that $\sigma(\ep)$ is unique
 for $\ep$ near 0 }

 \subsection{Continuation of a double spectral value}
This time we study the behavior of $\sig(\ep)$ when $\sig(0)$ has
multiplicity 2. Hence assume that $\tau\in\specb(\BB(.,0))$ has
multiplicity 2. Therefore, $\BB (\tau,0)=\II$. We start with the
following proposition.

{ \prop{ If $\BB(\tau,0)=\II$, then $\dd{\tr\BB}{\sigma}(\tau,0)=0$
and $\dis\frac{\pa^2\tr\BB}{\pa\sigma^2}(\tau,0)\ne 0$.}}

\vspace{.2cm}

 The proof of this proposition makes use of the the
following  lemma.

{\lem Given $M>0$, there is a positive constant $C$ such that
 \[
(1+2x)(1+2y)-4\sqrt{xy(1+x)(1+y)}\ge C,\qquad \forall x,y\in [0,\ M]
 \]}

\vspace{.2cm}

 \noindent\textit{Proof.} Consider the function
 $g(x,y)=(1+2x)^2(1+2y)^2-16xy(1+x)(1+y)$. It can be easily verified
 that $g(x,y)\ge 1$ in the square $[0,\ M]^2$.
 This implies in turn that
 \[
(1+2x)(1+2y)-4\sqrt{xy(1+x)(1+y)}\ge\frac{1}{1+8M+8M^2}\quad\Box
 \]
\vspace{.2cm}

 \noindent\textit{Proof of Proposition $2.9$.} Let $\VV(t,\sigma,\ep)$ be
 the fundamental matrix of equation $(\E)$ given by (2.7). Its
 derivative $\VV_\sigma$, with respect to $\sigma$, satisfies the system
 \begeq
\dot\VV_\sigma =\MM\VV_\sigma +\MM_\sigma\VV,\quad \VV_\sigma
(0,\sigma,\ep)=0
 \stopeq
 (the last condition follows from $\VV(0,\sigma,\ep)=\II$). Note
 that
 \[
 \dis\MM_\sigma =\DD_{i/\lam_\ep}=\left(\begar{ll} \dis\frac{i}{\lam_\ep} & 0
 \\\  0 & \dis\frac{-i}{\ov{\lam_\ep}}\stopar\right)\, .\]
 We consider (2.23) as a nonhomogeneous system in $\VV_\sigma$
 and we get
 \begeq
\VV_\sigma(t,\sigma,\ep)=\VV(t,\sigma,\ep)\int_0^t\VV(s,\sigma,\ep)^{-1}
\DD_{i/\lam_\ep}\VV(s,\sigma,\ep)ds\, .
 \stopeq
By using formula (2.7),  we have
\[
\VV^{-1}\DD_{i/\lam_\ep}\VV =i \left(\begar{ll} N_{11} & N_{12}
 \\ N_{21} & N_{22}\stopar\right)
\]
where
\[\begar{l}
N_{11}=\dis\frac{f(t,\sigma,\ep^2)\ov{f(t,\ov{\sigma},\ep^2)}}{\lam_\ep}
+\lam_\ep g(t,\sigma,\ep^2)\ov{g(t,\ov{\sigma},\ep^2)} \\
N_{12} =\dis
\frac{2a}{\lam_\ep}\ov{f(t,\ov{\sigma},\ep^2)g(t,\ov{\sigma},\ep^2)}\\
N_{21} =\dis
-\frac{2a}{\ov{\lam_\ep}}f(t,\sigma,\ep^2)g(t,\sigma,\ep^2)\\
N_{22}=-\dis\frac{f(t,\sigma,\ep^2)\ov{f(t,\ov{\sigma},\ep^2)}}{\ov{\lam_\ep}}
-\ov{\lam_\ep} g(t,\sigma,\ep^2)\ov{g(t,\ov{\sigma},\ep^2)}
 \stopar\]
 In particular
 \begeq
\tr(\VV^{-1}\DD_{i/\lam_\ep}\VV )
 (t,\sigma,0)=iN_{11}(t,\sigma,0)+iN_{22}(t,\sigma, 0)\equiv 0.
 \stopeq
 If we set $t=2\pi$ in (2.24), we get
 \begeq
\dd{\BB}{\sigma}(\sigma,\ep)=\BB(\sigma,\ep)\int_0^{2\pi}
\VV(s,\sigma,\ep)^{-1}
\DD_{i/\lam_\ep}\VV(s,\sigma,\ep)ds
 \stopeq
 Since, $\BB(\tau,0)=\II$, then it follows at once from (2.25) and
 (2.26) that $\dis\dd{\tr\BB}{\sigma}(\tau, 0)=0.$
 Now we compute $\VV_{\sigma\sigma}$. We have
 \begeq
\dot\VV_{\sigma\sigma}=\MM\VV_{\sigma\sigma}+2\DD_{i/\lam_\ep}
\VV_\sigma\, ,\quad \VV_{\sigma\sigma}(0,\sigma ,\ep)=0
 \stopeq
 and  after integrating this nonhomogenous system and using
 (2.24), we obtain
 \begeq\begar{ll}
\VV_{\sigma\sigma}(t,\sigma,\ep) & = \dis 2\VV(t,\sigma,\ep)
\int_0^t \VV^{-1}\DD_{i/\lam_\ep}\VV_\sigma(s,\sigma,\ep)ds\\
&
=2\dis\VV(t,\sigma,\ep)\int_0^t\int_0^s\LL(s,\sigma,\ep)\LL(u,\sigma,\ep)duds
 \stopar \stopeq
where $\LL(t,\sigma,\ep)=\VV^{-1}\DD_{i/\lam_\ep}\VV(t,\sigma,\ep)$.
We have in particular
 \begeq
 \LL(t,\sigma,0)=\frac{i}{a}\left(\begar{ll}
 P(t,\sigma) & 2\ov{Q(t,\ov{\sigma})}\\
 -2Q(t,\sigma) & -P(t,\sigma)
 \stopar\right)
 \stopeq
 where
 \begeq\begar{l}
 P(t,\sigma)=f(t,\sigma,0)\ov{f(t,\ov{\sigma},0)}+a^2
 g(t,\sigma,0)\ov{g(t,\ov{\sigma},0)}\\
 Q(t,\sigma)=af(t,\sigma,0)g(t,\sigma,0)
 \stopar\stopeq
 If we set $t=2\pi$ in (2.28), we get
 \[
 \BB_{\sigma\sigma}(\sigma ,\ep)=2\BB(\sigma,\ep)
 \int_0^{2\pi}\int_0^s\LL(s,\sigma,\ep)\LL(u,\sigma,\ep)duds
 \]
 and since $\BB(\tau,0)=\II$, we have
 \begeq
 \frac{\pa^2\tr(\BB)}{\pa\sigma^2}(\tau ,0)=2\int_0^{2\pi}\int_0^s
 \tr(\LL(s,\tau,0)\LL(u,\tau,0))duds
 \stopeq
It follows from (2.29) that
 \begeq
\frac{-a^2}{2}\tr(\LL(s,\tau,0)\LL(u,\tau,0))= P(s,\tau)P(u,\tau)-
4\mathrm{Re}\left[Q(s,\tau)\ov{Q(u,\tau)}\right]
 \stopeq
  Let $g(t,\tau,0)=\dis\frac{\rho(t)}{a}\ei{i\beta (t)}$
 (thus, $\rho =a|g|$ and $\beta$ is the argument of $g$) and then
since $f$ and $g$ satisfy (2.8) we have
$f(t,\tau,0)=(1+\rho(t)^2)^{1/2}\ei{i\alpha(t)}$. With this
notation, the functions $P$ and $Q$ become
 \[
 P(t,\tau)=1+2\rho(t)^2\quad\mathrm{and}\quad
 Q(t,\tau)=\rho(t)\sqrt{1+\rho(t)^2}\ei{i(\alpha(t)+\beta(t))}
 \]
If we set $x=\rho(s)^2$, $y=\rho(u)^2$, and
$\ta=\alpha(s)+\beta(s)-\alpha(u)-\beta(u)$, then formula (2.32)
becomes
 \begeq
\frac{-a^2}{2}\tr(\LL(s,\tau,0)\LL(u,\tau,0))= (1+2x)(1+2y)- 4
\sqrt{xy(1+x)(1+y)}\cos\ta
 \stopeq
 Since, $x,y$ are positive and bounded ($g$ is bounded over the
 interval $[0,\ \pi]$), then Lemma 2.2 implies that there is a
 positive constant $C$ such  that
 \[
 \tr(\LL(s,\tau,0)\LL(u,\tau,0)) \le -C\qquad\forall u,s\in [0,\
 2\pi].
 \]
 Therefore, by (2.31), we have $\tr(\BB_{\sigma\sigma}(\tau,0)\ne
 0$.
 This completes the proof of the Proposition $\quad \Box$

\vspace{.2cm}

 The behavior of the spectral values of $\clep$ is given by the following
 proposition.

 {\prop Suppose that $\BB (\tau, 0)=\II$. Then there exists
 $\ep_0>0$ such that the spectral values of
 $\clep$ through  $\tau$
 satisfy one of the followings:
 \begin{enumerate}
 \item there is a unique continuous function $\sigma(\ep)$
 defined in $[-\ep_0,\ \ep_0]$ such that $\sigma(\ep)\in\speclep$,
  $\sigma(\ep)\in\C\backslash\R$ for $\ep\ne 0$ and $\sigma(0)=\tau$;
 \item there are two continuous $\R$-valued functions $\sigma_1(\ep)$,
 $\sigma_2(\ep)$ defined in $[-\ep_0,\ \ep_0]$,
 such that $\sigma_1(\ep),\sigma_2(\ep)\in\speclep$,
 $\sigma_1(\ep)<\sigma_2(\ep)$ for $\ep\ne 0$, and
 $\sigma_1(0)=\sigma_2(0)=\tau$;
 \item there is a unique continuous $\R$-valued function
 $\sigma(\ep)$ defined in $[-\ep_0,\ \ep_0]$ such that
 $\sigma(\ep)\in\speclep$ and $\sigma(0)=\tau$
 \end{enumerate}}

\vspace{.2cm}

 \Proof{ It follows from Proposition 2.9 that the
spectral function $F(\sigma,\ep)$ defined in (2.6) satisfies
\[
F(\tau,0)=0,\quad \dd{F}{\sigma}(\tau,0)=0,\quad \mathrm{and}\quad
\frac{\pa^2F}{\pa\sigma^2}(\tau,0)\ne 0\, .
\]
Since $F$ is analytic in both variables, then by the Weierstrass
Preparation Theorem (see {\cite{Gol-Gui}}) we can find analytic
functions $G(\sigma,\ep)$, $A_0(\ep)$ and $A_1(\ep)$ with
$(\sigma,\ep)$ near $(\tau,0)\in\C\times \R$, such that
$G(\tau,0)\ne 0$, $A_1(0)=A_0(0)=0$ and
 \[
F(\sigma,\ep)=G(\sigma,\ep)\left[(\sigma-\tau)^2-2A_1(\ep)(\sigma-\tau)
+A_0(\ep)\right]\, .
 \]
Thus, there exists $\ep_0>0$ such that the roots of the spectral
equation $F(\sigma,\ep)=0$ in a neighborhood of $(\tau,0)$ are
 given by the quadratic formula
 \[
\sigma_{1,2}(\ep)=\tau +A_1(\ep)\pm\sqrt{A_1^2(\ep)-A_0(\ep)}\, .
 \]
The conclusion of the proposition follows depending on the sign of
the discriminant $A_1^2-A_0$ }

\subsection{Purely imaginary spectral value}
 We study here the behavior of the monodromy matrix at possible
 spectral value on the imaginary axis.

{\prop Suppose that for some $\ep_0\in\R^\ast$, the operator
$\cl_{\ep_0}$ has a spectral value $\sigma_0$ of the form
 \begeq
\sigma_0=i\frac{|\lam_{\ep_0}|^2}{2b\ep_0}k,\quad \mathrm{with}\quad
k\in \Z^\ast\, .
 \stopeq
Then the monodromy matrix $\BB(\sigma_0,\ep_0)$ is similar to
 $\dis\left(\begar{ll}1 & 1\\ 0 & 1\stopar\right)$.}

\vspace{.2cm}

 \Proof{ If $\sigma_0$ given by (2.34) is a spectral
value, then $\det\BB(\sigma_0,\ep_0)=1$, by (2.5). Hence 1 is an
eigenvalue of $\BB(\sigma_0,\ep_0)$ and, a priori, it could have
multiplicity 2. In which case $\BB(\sigma_0,\ep_0)=\II$. We are
going to show that this case does not happen.  By contradiction,
suppose that $\BB(\sigma_0,\ep_0)=\II$. First we prove that
$\tr(\BB_\sigma(\sigma_0,\ep_0))\ne 0$. From formulas (2.26),
(2.25), and (2.8) we have
\[
\dd{\tr(\BB)}{\sigma}(\sigma_0,\ep_0)=\int_0^{2\pi}m(s)ds
\]
where
 \[\begar{ll}
 m(s) & =\dis \frac{\ov{\lam_{\ep_0}}-\lam_{\ep_0}}{|\lam_{\ep_0}|^2}
 \left[f(s,\sigma_0,\ep_0)\ov{f(s,\ov{\sigma_0},\ep_0)}-
 |\lam_{\ep_0}|^2g(s,\sigma_0,\ep_0)\ov{g(s,\ov{\sigma_0},\ep_0)}\right]
 \\
 & =\dis\frac{-2ib\ep_0^2}{|\lam_{\ep_0}|^2}
 \stopar\]
 This shows that $\tr(\BB_\sigma(\sigma_0,\ep_0))\ne 0$. Hence,
 it follows  that the spectral function
satisfies
\[
\dd{F}{\sigma}(\sigma_0,\ep_0)=\dd{\tr(\BB)}{\sigma}(\sigma_0,\ep_0)\ne
0\, .
\]
By the implicit function theorem, the germ of the analytic variety
$F(\sigma,\ep)=0$ through $(\sigma_0,\ep_0)$ is smooth and there is
a unique analytic function $\sigma(\ep)$ defined near $\ep_0$ with
$\sigma(\ep_0)=\sigma_0$ such that $F(\sigma(\ep),\ep)\equiv 0$. It
follows that $\BB(.,\ep)$ has a unique spectral value through
$(\sigma_0,\ep_0)$.

Let $U_1(\ep)$ be a continuous eigenvector (with eigenvalue 1) of
$\BB (\sigma(\ep),\ep)$ defined in an interval $(\ep_0-\delta ,\
\ep_0+\delta)$ for some
$\delta>0$. We can assume that $U_1(\ep_0)=\dis\left(\begar{l}\alpha\\
\beta\stopar\right)$ with $\alpha\ne 0$. Now, consider the equation
 \[
 G(\sigma,\ep,z)=\left(\BB(\sigma,\ep)-\II\right)\,\left(\begar{l} z\\
1\stopar\right)\, =0
 \]
 in a neighborhood of $(\sigma_0,\ep_0,0)\in\C\times\R\times\C$.
 This equation defines a germ of an analytic variety of real
 dimension 1 in $\C\times\R\times\C$ that passes through the point
 $(\sigma_0,\ep_0,0)$. Therefore, there exists $\delta_1>0$ and
 continuous functions $\sigma'(\ep)$ and $z(\ep)$ defined in
 $[\ep_0-\delta_1,\ \ep_0+\delta_1]$ such that
 $\sigma'(\ep_0)=\sigma_0$, $z(\ep_0)=0$ and
 $G(\sigma'(\ep),\ep,z(\ep)\equiv 0$. The continuous vector
 $\dis U_2(\ep)=\left(\begar{l} z(\ep)\\
1\stopar\right)$ is therefore an eigenvector with eigenvalue 1 of
$\BB (\sigma'(\ep),\ep)$. Moreover, $U_1(\ep)$ and $U_2(\ep)$ are
independent for $\ep$ close to $\ep_0$. By the uniqueness of the
spectral value established above, we get $\sigma'(\ep)=\sigma(\ep)$.
This means that $\BB(\sigma(\ep),\ep)\equiv \II$ for  $\ep$ close to
$\ep_0$. From this and (2.5) we get that
\[
\det(\BB(\sigma(\ep),\ep))=\exp\left(\frac{4\pi
b\ep}{|\lam_\ep|^2}\sigma(\ep)\right)\equiv 1
\]
and therefore, $\sigma(\ep)=\dis i\frac{|\lam_{\ep}|^2}{2b\ep}k$ for
every $\ep$. Since $\BB$ is analytic, we get that
$\BB(\sigma(\ep),\ep)=\II$ for every $\ep\in\R$. Now if we go back
to the system $(\mathrm{E}_{\sigma(\ep),\ep})$, we can assume, by
continuity and Corollary 2.0.1,  that for the solution
$(\phi,\psi)$, the function $\phi$ is dominating for every $\ep$,
i.e. $|\phi(t,\ep)|>|\psi(t,\ep)|$. The winding number
$j_0=\ind(\phi)$ is then constant and we get from the first equation
of $(\mathrm{E}_{\sigma(\ep),\ep})$ that
\[
\frac{1}{2\pi}\lam_\ep\int_0^{2\pi}\frac{\phi'(t,\ep)}{\phi(t,\ep)}dt=
\lam_\ep j_0=(\sigma(\ep)-\lam_\ep\nu)+\frac{1}{2\pi i}\int_0^{2\pi}
c(t)\frac{\psi(t,\ep)}{\phi(t,\ep)}dt
\]
By taking the limit, we obtain
 \[
 \lim_{\ep\to 0}\sigma(\ep)=\lim_{\ep\to 0}i\frac{|\lam_{\ep}|^2}{2b\ep}k
=\lim_{\ep\to 0}\left(\lam_\ep(j_0+\nu)+\frac{1}{2\pi
i}\int_0^{2\pi} c(t)\frac{\psi(t,\ep)}{\phi(t,\ep)}dt\right)
 \]
Since the right hand side is bounded and $\lam_\ep =a+ib\ep$  with
$a>0$ and $b\ne 0$, then necessarily $k=0$ and this is a
contradiction. This proves that $\BB (\sigma_0,\ep_0)\ne\II$ }

The following corollary is a direct consequence of Proposition 2.11
and formula (2.5)

 {\cor If 1 is an eigenvalue of the monodromy matrix
 $\BB(\sigma,\ep)$ with $\ep\ne 0$, then it has
 multiplicity one}

\subsection{Main result about basic solutions}
 The following theorem summarizes the main properties of the basic
 solutions of $\clep$.

{\theo For every $j\in\Z$ there are exactly two $\R$-independent
basic solutions $w_j^+(r,t,\ep)$ and $w_j^-(r,t,\ep)$ of $\clep$
with $\charac (w_j^\pm)=(\sigpmj ,j)$ such that the spectral values
$\sigpmj\in\speclep$ satisfy
\begin{itemize}
\item $\sigpmj(\ep)$ depends continuously on $\ep$ and
\item if for some $\ep_0\in\R$, $\sigpj(\ep_0)\in\C\backslash\R$,
then $\sigmj(\ep_0)=\sigpj(\ep_0)$
\end{itemize}}

\vspace{.2cm}

 \Proof{ Consider the analytic variety $\Gamma
=\{(\sigma,\ep)\in \C\times\R;\ F(\sig,\ep)=0\}$ where $F$ is the
spectral function given in (2.6). Thus the real spectral values
$\tau_k^\pm$ of $\cl_0$ (see Proposition 2.3) satisfy $(\tau_k^-,0)$
and $(\tau_k^+,0)\in\Gamma$. Let $\Gamma_j^\pm$ be connected
components of $\Gamma$ containing $(\tau_j^\pm,0)$ and
$\Gamma_j=\Gamma_j^-\cup\Gamma_j^+$.

For $j\ne k$ we have $\Gamma_k\cap\Gamma_j=\emptyset$. Indeed, if
there is $(\sigma_0,\ep_0)\in\Gamma_k\cap\Gamma_j$ (with $\ep_0\ne
0$), then equation $(\mathrm{E}_{\sigma_0,\ep_0})$ would have
periodic solutions,
\[
\left(\begar{l}\phi_j(t,\sigma_0,\ep_0)\\
\psi_j(t,\sigma_0,\ep_0)\stopar\right)\quad\mathrm{and}\quad
\left(\begar{l}\phi_k(t,\sigma_0,\ep_0)\\
\psi_k(t,\sigma_0,\ep_0)\stopar\right)
\]
giving rise to basic solutions
\[
w_j=r^{\sigma_0}\phi_j+\ov{r^{\sigma_0}\psi_j}\quad\mathrm{and}\quad
w_k=r^{\sigma_0}\phi_k+\ov{r^{\sigma_0}\psi_k}
\]
with winding numbers $j$ and $k$, respectively. But the monodromy
matrix $\BB(\sigma_0,\ep_0)$ has only one eigenvector with
eigenvalue 1 (Corollary 2.0.5). This means $\phi_k=c\phi_j$ and
$\psi_k=c\psi_j$ for some constant $c$. Hence,
 $w_k=r^{\sigma_0}c\phi_j+\ov{r^{\sigma_0}c\psi_j}$ has also winding
 number $j$. This is a contradiction.

 To complete the proof, we need to show that for every $j\in\Z$
 and for every $\ep_0\in\R$
 \[
 \left|\Gamma_j\cap\{(\sigma,\ep_0);\ \sigma\in\C\}\right| \, \le\,
 2 \, , \]
 where $|S|$ denotes the cardinality of the set $S$.
 By contradiction, suppose that there exists $j_0\in\Z$ and
 $\ep_0\in\R^\ast$ such that
\[
 \left|\Gamma_{j_0}\cap\{(\sigma,\ep_0);\ \sigma\in\C\}\right| \, \ge\,
 3 \, .\]
 Let $(\sigma_1,\ep_0),\ (\sigma_2,\ep_0)$, and $(\sigma_3,\ep_0)$
 be three distinct points in $\Gamma_{j_0}$. Hence, $\Gamma_{j_0}$
 has three distinct components $C_1,\ C_2$, and $C_3$ over a
 neighborhood of $\ep_0$. They are defined by functions $f_1(\ep)$,
 $f_2(\ep)$ and $f_3(\ep)$, that are analytic, except possibly at
 $\ep_0$. By analytic continuation, $\Gamma_{j_0}$
  has three
 distinct branches $C_1$, $C_2$, $C_3$, parametrized by $\ep$.
 That is $C_l=\{ (f_l(\ep),\ep);\ \ep\in\R\}$ with $f_l\in C^0(\R)$
 and analytic everywhere, except on a set of isolated points. In
 particular for $\ep =0$, we get
 \[
\{ f_1(0),\, f_2(0),\, f_3(0)\}=\{\tau_{j_0}^-,\, \tau_{j_0}^+\}
 \]
 In the case $\tau_{j_0}^- < \tau_{j_0}^+$, we can assume
 that $f_1(0)=f_2(0)$ and this contradicts Proposition 2.8.
 In the case $\tau_{j_0}^- = \tau_{j_0}^+$, we would have
 $f_1(0)=f_2(0)=f_3(0)=\tau_{j_0}^+$ and this would contradict
 Proposition 2.10}

For the adjoint operator we have the following theorem.

 {\theo Let
 \[w_\ep(r,t)=r^{\sigma(\ep)}\phi(t,\ep)+
 \ov{r^{\sigma(\ep)}\psi(t,\ep)}\]
  be a
 basic solution of $\clep$ with
 $\charac(w_\ep)=(\sigma(\ep) ,j).$
 Then $\clepstar$ has a basic solution
 \[w^\star_\ep(r,t)=r^{-\sigma(\ep)}X(t,\ep)+\ov{r^{-\sigma(\ep)}Z(t,\ep)}\]
 with
 $\charac(w^\star_\ep)=(-\sigma(\ep) ,-j).$}

 \vspace{.2cm}

 \Proof{ For $\sigma(\ep)\in\specb(\BB(.,\ep))$, it follows from
 (2.17) and (2.16) that $\sigma(-\ep)=\sigma(\ep)$ and
 $\ov{\sigma(\ep)}=\ov{\sigma(-\ep)}\in\specb(\BB(.,-\ep))$.
 Let $U(\ep)$ be a continuous eigenvector with eigenvalue 1 of
 $\BB(\ov{\sigma(-\ep)},-\ep)$ such that
 \[
U(0)=\left(\begar{l}\phi(0,0)\\ \psi(0,0)\stopar\right)\ .
 \]
 Then $\VV(t,\ov{\sigma(-\ep)},-\ep)U(\ep)$ is a periodic solution
 of $(\mathrm{E}_{\ov{\sigma(-\ep)},-\ep})$. The function
 \[
\left(\begar{l}X(t,\ep)\\ Z(t,\ep)\stopar\right) = \DD
 \ov{\VV(t,\ov{\sigma(-\ep)},-\ep)U(\ep)}=\widetilde{\VV}(t,-\sigma(\ep),\ep)
 \DD\ov{U(\ep)}
 \]
 is a periodic solution of the adjoint system
 $(\widetilde{\mathrm{E}}_{-\sigma(\ep),\ep})$. Furthermore,
 $|X|>|Z|$ and $\ind(X)=-j$
 so that the character of the associated basic solution
 is $(-\sigma(\ep),-j)$ }

 \section{Example}
 We give here an example in which the basic solutions can be explicitly
 determined. This is
 the case when $c(t)=ic_0\ei{ikt}$ with $c_0\in\C^\ast$. For
  simplicity, we assume that $\nu =0$. The system of equations
 $(\E)$ is
 \[\begar{ll}
 \lam_\ep\dot\phi(t) & = i\sigma\phi(t)+ic_0\ei{ikt}\psi(t)\\
 \ov{\lam_\ep}\dot\psi(t) & =-i\sigma\psi(t)-i\ov{c_0}
 \ei{-ikt}\phi(t)
 \stopar \]
 In this case we can use Fourier series to determine
 the spectral values and the periodic solutions. For a
 given $j\in\Z$, the system has a solution of
 the form $\phi (t)=\ei{ijt}$, $\psi (t)=D\ei{i(j-k)t}$
 with  $\sigma$ and $D$ satisfying
 \[\begar{l}
 \lam_\ep j  =\sigma +c_0D\\
 \ov{\lam_\ep}(j-k)D=-\sigma D-\ov{c_0}
 \stopar \]
 The elimination of $D$, gives the following quadratic
 equation for the spectral value $\sigma$
 \[
 \sigma^2-[(\lam_\ep-\ov{\lam_\ep})j+k\ov{\lam_\ep}]
 \sigma -[j(j-k)|\lam_\ep|^2+|c_0|^2]=0
 \]
 After replacing $\lam_\ep$ by $a+ib\ep$ we get
 $\sigma$ and $D$:
 \[\begar{l}\dis
\sigma_j=ib\ep j+\frac{(a-ib\ep)k}{2}+
 \sqrt{\left(aj-\frac{(a-ib\ep)k}{2}\right)^2+|c_0|^2}\\
 \dis D_j=\frac{(a+ib\ep)j-\sigma_j}{c_0}
 \stopar\]
 The corresponding basic solution of $\clep$ is
 \[
 w_j(r,t,\ep)=r^{\sigma_j}\ei{ijt}+
 \ov{r^{\sigma_j}D_j\ei{i(j-k)t}}\,.
 \]
 The character of $w_j$ is $(\sigma_j,j)$ if $|D_j|<1$
 and it is $(\ov{\sigma_j},k-j)$ if $|D_j|>1$.

 Note that in order for $D_j$ to have norm 1, for
 some $j_0$, say
 $D_{j_0}=\ei{i\alpha}$, the exponent $\sigma_{j_0}$ needs to satisfy
 \[
\sigma_{j_0}=\lam_\ep j_0-c_0\ei{i\alpha},\quad\mathrm{and} \quad
 \sigma_{j_0}= -\ov{\lam_\ep} (j_0-k)-\ov{c_0}\ei{-i\alpha}\, .
 \]
 Consequently, $\lam_\ep (2j_0-k)=\sigma_{j_0}-\ov{\sigma_{j_0}}$.
 Since $\mathrm{Re}(\lam_\ep)=a>0$, then necessarily $k=2j_0$ is an
 even integer.

 Thus for $k$ odd, $|D_j|\ne 1$ for every $j\in\Z$ and
 for an even $k$, $k=2j_0$, we have $|D_j|\ne 1$ for every
 $j\ne j_0$. At the level $j_0$, we get
 \[
 \sigma_{j_0}=aj_0+\sqrt{-b^2\ep^2j_0^2+|c_0|^2}
 \quad\mathrm{and}\quad
 D_{j_0}=\frac{ib\ep j_0-\sqrt{-b^2\ep^2j_0^2+|c_0|^2}}{c_0}
 \]
and the character of the corresponding basic solution is
$(\sigma_{j_0},j_0)$.

\section{Asymptotic behavior of the basic solutions of $\cl$}
 In this section, we determine the asymptotic behavior of the basic
 solutions. This behavior will be needed in the next section
 to estimate the kernels. From now on, there is no need anymore for
 the parameter $\ep$. So we will denote $\cl_1$ by $\cl$ and the
 associated system of differential equations $(\mathrm{E}_{\sig,1})$
 by $(\mathrm{E}_\sig)$. We will assume here that $\lam=a+ib$ with
 $a>0$ and $b\ne 0$, since the asymptotic behavior in case $b=0$ is
 known
 from {\cite{Mez-CV3}}. Hence,
\[
 \cl u =Lu+i\lam\nu u-c(t)\ov{u}\, ,
\]
where $L$ is the vector field given in (1.1).
 Let
 \begeq
 \gamma =\frac{1}{4a\pi}\int_0^{2\pi}|c(t)|^2dt\quad
 \mathrm{and}\quad
 k(t)=\frac{1}{\lam}
 \left[\gamma t-\frac{1}{2a}\int_0^t|c(s)|^2ds\right]
 \stopeq
 Note that $k(t)$ is $2\pi$-periodic. We have the following
 theorem

 {\theo For $j\in\Z$, the operator $\cl$ has basic solution
 \[
w_j(r,t)=r^{\sigma_j}\phi_j(t)+\ov{r^{\sigma_j}\psi_j(t)}
 \]
 with character $(\sigma_j,j)$ such that, as $|j|\to\infty$,
 we have
 \begin{eqnarray}
 \sigma_j & = & \dis\lam (j+\nu)+\frac{\gamma}{j}+
 O(j^{-2})\\
 \phi_j(t) & = & \dis \ei{ijt}
 \left( 1+i\frac{k(t)}{j}\right) +O(j^{-2})\\
 \psi_j(t) & = & -\dis \ei{ijt}\frac{\ov{c(t)}}{2aj}
  +O(j^{-2})
 \end{eqnarray}
  where $\gamma$ and $k(t)$ are given in $(4.1)$.
Furthermore, any basic solution  (with $|j|$ large) has the form
\[
 w(r,t)=r^{\sigma_j}a\phi_j(t)+\ov{r^{\sigma_j}a\psi_j(t)},
\]
 with $a\in\C$}

\vspace{.2cm}

\Rem{\bf 4.1}{ It follows from Theorems 2.2 and 4.1 that for
$|j|\in\Z$ large, the basic solutions of $\clstar$ have the form,
\[
w^\ast(r,t)=r^{-\sigma_j}X_{-j}(t)+\ov{r^{-\sigma_j}Z_{-j}(t)}
\]
with $\charac(w^\ast) =(-\sigma_j,-j)$ where
 \[\begar{ll}
 X_{-j}(t) & =  \dis \ei{-ijt}
 \left( 1-i\frac{k(t)}{j}\right) +O(j^{-2})\\
 Z_{-j}(t) & =  -\dis \ei{-ijt}\frac{c(t)}{2aj}
  +O(j^{-2})
 \stopar\]
 }

 The remainder of this section deals with the proof of Theorem 4.1.
 The proof will be divided into 3 steps.
 To simplify the expressions, we will use the following variables
 \begeq
 \mu =\frac{\sigma -\lam\nu}{\lam}, \
 \ei{ix}={\lam}/{\ov{\lam}},\ \delta =(\ei{-ix}+1)\nu,\
 \mathrm{and}\ c_1(t)=\frac{c(t)}{\lam}
 \stopeq
The system of equations  $(\mathrm{E}_\sigma)$ becomes then
 \begeq\begar{ll}
 \dot\phi & =i\mu \phi +c_1(t)\psi\\
 \dot\psi & =-i\ei{ix}(\mu +\delta )\psi +\ov{c_1(t)}\phi
 \stopar\stopeq
 Now, we proceed with the proof of the theorem.

 \subsection{Estimate of $\sigma$}
 For a periodic solution $(\phi,\psi)$ of (4.6) with
  $|\psi|<|\phi|$, we can assume that
  $\max|\phi|=1$ and $\ind(\phi)=j$. Let
  \begeq
 T(t)=\frac{\psi(t)}{\phi(t)}
  \stopeq
 It follows from the first equation of (4.6) that
\[
\mu +\frac{1}{2\pi i}\int_0^{2\pi}c_1(t)T(t)dt=\frac{1}{2\pi i}
 \int_0^{2\pi}\frac{\dot\phi (t)}{\phi(t)}dt=j
\]
Hence,
 \begeq
 \mu =j-M_j\quad\mathrm{with}\quad
 M_j=\frac{1}{2\pi i}\int_0^{2\pi}c_1(t)T(t)dt\, .
 \stopeq
 Note that $|M_j|\le\dis\frac{1}{2\pi}\int_0^{2\pi}|c_1(t)|dt$.
 To obtain a better estimate of $M_j$ we use the second equation of
 (4.6) to get (after multiplying by $c_1$, dividing
 by $\phi$ and integrating over $[0,\ 2\pi ]$)
 \begeq
 -i\ei{ix}(\mu+\delta)M_j =\frac{-1}{2\pi i}\int_0^{2\pi}
 |c_1(t)|^2dt+\frac{1}{2\pi i}\int_0^{2\pi}
 c_1(t)\frac{\psi'(t)}{\phi(t)}dt
 \stopeq
 We use integration by parts in the last integral together with
 (4.6) to obtain
 \[\begar{ll}
 \dis\int_0^{2\pi}\!\!\!
 c_1(t)\frac{\psi'(t)}{\phi(t)}dt &=
 \dis -\int_0^{2\pi}\!\!\! c_1'(t)T(t)dt+\int_0^{2\pi}\!\!\!
 c_1(t)
 \frac{(i\mu\phi(t)+c_1(t)\psi(t))\psi(t)}{\phi(t)^2}dt\\
 & =\dis \int_0^{2\pi}(-c_1'(t)T(t)+i\mu
 c_1(t)T(t)+c_1(t)^2T(t)^2)dt
 \stopar\]
 From this and (4.9), we deduce that
 \[
 -i\left[(1+\ei{ix})\mu+\delta\ei{ix}\right]M_j =
 \frac{-1}{2\pi
 i}\int_0^{2\pi}(|c_1(t)|^2+c_1'(t)T(t)-c^2_1(t)T^2(t))dt.
 \]
Consequently,
 \begeq
 |(1+\ei{ix})\mu +\delta\ei{ix}|\, |M_j|\le A_1\, ,
 \stopeq
where
 $\dis A_1=\frac{1}{2\pi}\int_0^{2\pi}(2|c_1(t)|^2+|c_1'(t)|)dt$.
 Note that since $1+\ei{ix}=2a/\ov{\lam}$ satisfies
 $\mathrm{Re}(1+\ei{ix})>0$ and $|\delta |<2$, then it follows from
 (4.8) and the boundedness of $M_j$ that
 there exists $J_0\in\Z^+$ such that
 \begeq
 |(1+\ei{ix})\mu +\delta\ei{ix}|\ge\frac{|j|}{2},\qquad
 \forall j\in\Z,\ \ |j|\ge J_0.\stopeq

 {\lem Let $N_j=jM_j$, then
 \begeq
 \lim_{|j|\to\infty}N_j=\frac{-1}{2\pi (1+\ei{ix})}\int_0^{2\pi}
 |c_1(t)|^2dt=-\frac{\gamma}{\lam}
 \stopeq
 where $\gamma$ is given in $(4.1)$}

 \vspace{.2cm}

\Proof{ It follows from (4.10) and (4.11) that $|N_j|\le  2A_1$. Let
 \begeq
 \phi(t)=\ei{ijt}\phi_1(t)\quad\mathrm{and}\quad
 \psi(t)=\ei{ijt}\psi_1(t)
 \stopeq
 Hence, $\ind(\phi_1)=0$ and $T=\psi_1/\phi_1$. From
 (4.6), we get the system for $\phi_1,\ \psi_1$
\begeq\begar{ll}
 \dot\phi_1 & =-\dis i\frac{N_j}{j}\phi_1 +c_1(t)\psi_1\\
 \dot\psi_1 & =-iE_j\psi_1 +\ov{c_1(t)}\phi_1
 \stopar\stopeq
where
 \begeq
 E_j=(1+\ei{ix})j-\ei{ix}\frac{N_j}{j}+\ei{ix}\delta\, .
 \stopeq
Note that since $|N_j|\le 2 A_1$ and $\ei{ix}\ne -1$, then for $|j|$
large ($|j|\ge J_0$), we have
 \begeq
|E_j|\ge |1+\ei{ix}|\frac{|j|}{2}\, .
 \stopeq
 Now, it follows from the first equation of (4.14) and from
 $\ind(\phi_1)=0$ that
 \begeq
\frac{N_j}{j}=\frac{1}{2\pi i}\int_0^{2\pi}c_1(t)T(t)dt\, .
 \stopeq
 We use the second equation of (4.14) to estimate the
 integral appearing in (4.17)
 \begeq
 \dis\int_0^{2\pi}\!\!\!
 c_1Tdt =\dis\frac{-1}{iE_j}\int_0^{2\pi}\!\!\!
 c_1\frac{\psi'_1-\ov{c_1}\phi_1}{\phi_1}dt
  =\dis\frac{1}{iE_j}\int_0^{2\pi}\!\!\left(
 |c_1|^2dt-
c_1 \frac{\psi'_1}{\phi_1}\right)dt
 \stopeq
 We use integration by parts and system (4.14) to evaluate
 the last integral appearing in (4.18).
 \begeq\begar{ll}
 \dis\int_0^{2\pi}c_1\frac{\psi_1'}{\phi_1}dt & =\dis
 -\int_0^{2\pi}\left[c'_1T+c_1
\psi_1 \frac{(-iN_j/j)\phi_1+c_1\psi_1}{\phi^2_1}\right]dt\\
 & =\dis\int_0^{2\pi}(-c'_1T+c_1^2T^2)dt-\frac{iN_j}{j}
 \int_0^{2\pi}c_1Tdt
 \stopar\stopeq
 Thus,
 \begeq
 \int_0^{2\pi}c_1Tdt=\frac{1}{iE_j}\int_0^{2\pi}(|c_1|^2+c'_1T-
 c_1^2T^2)dt+\frac{N_j}{jE_j}\int_0^{2\pi}c_1Tdt
 \stopeq
 Therefore, from (4.20) and (4.16), we have
 \begeq
 \int_0^{2\pi}c_1Tdt=\frac{1}{iE_j}\left[
 \int_0^{2\pi}|c_1|^2dt +I_1-I_2\right] +O(\frac{1}{j^2})
 \stopeq
where $I_1$ and $I_2$ are given by
\[
 I_1=\int_0^{2\pi}c_1'(t)T(t)dt \quad\mathrm{and}\quad
   I_2=\int_0^{2\pi}c_1^2(t)T^2(t)dt .
\]
Now we show that $I_1=O(1/j)$ and $I_2=O(1/j)$. For $I_1$, we have,
after using system (4.14), that
 \[\begar{ll}
 I_1 & =\dis\frac{-1}{iE_j}\int_0^{2\pi}c_1'
 \frac{\psi_1'-\ov{c_1}\phi_1}{\phi_1}dt\\
 & =\dis\frac{-1}{iE_j}\left[\int_0^{2\pi} |c_1|^2dt+
 \int_0^{2\pi}c_1''Tdt-\int_0^{2\pi}c_1'\psi_1
 \frac{\phi_1'}{\phi_1^2}dt \right]\\
 & =\dis\frac{-1}{iE_j}\left[\int_0^{2\pi} (|c_1|^2+c_1''T)dt-
 \int_0^{2\pi}
 c_1'\psi_1
 \frac{(-iN_j/j)\phi_1+c_1\psi_1}{\phi_1^2}dt \right]
 \stopar\]
So
 \begeq
 iE_jI_1=\int_0^{2\pi}(|c_1|^2+c_1''T-c_1'c_1T^2)dt+
 \frac{iN_j}{j}\int_0^{2\pi}c_1'Tdt
 \stopeq
 Since $|T|<1$ and $|N_j|<2A_1$, we get $I_1=O(1/j)$.
For $I_2$, we use system (4.14) and integration by parts, to obtain
 \begeq\begar{ll}
 -iE_jI_2 & =
 \dis\int_0^{2\pi}\!\!\!\!\! c_1^2\frac{\psi_1}{\phi_1^2}
 (\psi_1'-\ov{c_1}\phi_1)dt=-\int_0^{2\pi}\!\!\!\!\!
 c_1^2\ov{c_1}T dt
 +\frac{1}{2}\int_0^{2\pi}\!\!\!\!\!
 c_1^2\frac{(\psi_1^2)'}{\phi_1^2}dt\\
 & =\dis -\int_0^{2\pi}\!\!\! c_1^2\ov{c_1}T dt-\int_0^{2\pi}
 \!\!\! c_1c'_1T^2dt+\int_0^{2\pi}c_1^2\psi_1^2
 \frac{\phi_1'}{\phi_1^3}dt\\
 & =\dis -\int_0^{2\pi}\!\!\!
 (c_1^2\ov{c_1}T+c_1c_1'T^2)dt+
 \int_0^{2\pi}\!\!\! c_1^2\psi_1^2
 \frac{(-iN_j/j)\phi_1+c_1\psi_1}{\phi_1^3}dt\\
 & =\dis -\int_0^{2\pi}\!\!\!
 (c_1^2\ov{c_1}T+c_1c_1'T^2-c_1^3T^3)dt-\frac{iN_j}{j}
 \int_0^{2\pi}\!\!\! c_1^2T^2dt
 \stopar\stopeq
and again $|E_jI_2|$ is bounded and therefore $I_2=O(1/j)$. With
these estimates for $I_1$ and $I_2$, expressions (4.17), (4.20), and
(4.21) give
 \[
N_j=-\frac{j}{2\pi E_j}\left[\int_0^{2\pi}\!\!\! |c_1|^2dt
+O(|j|^{-1})\right]\, .
 \]
 Since $\dis\lim_{|j|\to\infty}\frac{j}{E_j}=\frac{1}{1+\ei{ix}}$,
 the lemma follows}

 Using this lemma we get
 \[
 M_j=\frac{N_j}{j}=\frac{-\ov{\lam}}{ 2\pi j(\lam +\ov{\lam})}
 \int_0^{2\pi}\!\!\!\frac{|c(t)|^2}{|\lam|^2}dt +O(j^{-2})=
 -\frac{\gamma}{\lam j}+O(j^{-2}).
 \]
 Consequently, $\dis\mu =j-M_j=j+\frac{\gamma}{j\lam}+O(j^{-2})$
 and $\sigma =\lam (\mu +\nu)$ gives estimate (4.2) of the
 theorem.

 \subsection{First estimate of $\phi$ and $\psi$}
 We begin with the estimates of the components $\phi$ and $\psi$ of the
 basic solutions.
Let $\phi_1$, and $\psi_1$ be the functions defined in (4.13).
 Note that $\max|\phi_1|\le 1$ and $|\psi_1|<|\phi_1|$.
 We have the following lemma that gives an estimate of $\psi_1$.

{\lem There exist $J_0\in\Z^+$ and $K>0$ such that  the function
$\psi_2 =j\psi_1$ satisfies
 \begeq
 |\psi_2(t)|=|j\psi_1(t)|\le K ,\qquad\forall t\in\R,\
 \forall j\in\Z,\ \ |j|\ge J_0\, .
 \stopeq }

\vspace{.2cm}

\Proof{ The system (4.14) implies that
 $T=\psi /\phi =\psi_1/\phi_1$ satisfies the equation
 \begeq
 T'(t)=-L_jT(t)-c_1(t)T^2(t)+\ov{c_1}(t)
 \stopeq
where
 \begeq
 L_j=i\left(E_j-\frac{N_j}{j}\right)=
 i\left[(1+\ei{ix})j+\ei{ix}\delta +O(1/j)\right]\, .
 \stopeq
Note that the real part, $p$, of $L_j$ is given by
\[
 p=\mathrm{Re}(L_j)=-j\sin x +\mathrm{Re}(\delta\ei{ix})+
 O(1/j)\, .
\]
and, since $\sin x\ne 0$ (because $b\ne 0$),  there exists
$J_0\in\Z^+$ such that
 \begeq
 |p|\ge \frac{|\sin x|}{2}|j|\qquad\forall j\in\Z,\ |j|\ge J_0
 \stopeq
 Let $\rho (t)=|T(t)|$, $A (t)=\arg T(t)$,
 $\vartheta (t)=\arg c_1(t)$ and $M=\dis\max_{0\le t\le
 2\pi}|c_1(t)|$.
 Let us rewrite equation (4.25) as
 \begeq
 \rho'+iA'\rho =-L_j\rho -|c_1|\rho^2\ei{i(A+\vartheta)}+
 |c_1|\ei{-i(A+\vartheta)}.
 \stopeq
 By taking the real part, we obtain
 \begeq
 \rho'=-p\rho-|c_1|\rho^2\cos(A+\vartheta)+|c_1|\cos(A+\vartheta).
 \stopeq
 Since $0\le \rho <1$, we get $-2M\le \rho'+p\rho\le 2M$. Equivalently,
 \begeq
 -2M\ei{pt}\le (\rho(t)\ei{pt})'\le 2M\ei{pt}.
 \stopeq
 In the case where $p>0$ (and so $p>|j\sin x|/2$), we obtain, after
 integrating (4.30) from 0 to $t$, with $t>0$, that
 \begeq
 \rho(t)\le \left(\rho(0)-\frac{2M}{p}\right)\ei{-pt}+\frac{2M}{p}.
 \stopeq
 Let $t_0>0$ be such that $\ei{-pt_0}\le 2M/p$. Then, it follows from
 (4.31), that
 \begeq
 0\le \rho(t)\le \frac{2M}{p}\left[
 1+\rho(0)-\frac{2M}{p}\right]\le \frac{K}{j},\qquad
 \forall t\ge t_0
 \stopeq
 where $K$ is a constant independent on $j$. Since the function
 $\rho$ is periodic, then inequality (4.32) holds for
 every $t\in\R$. When $p<0$, an analogous argument
 (integrating (4.30) from $t$ to 0 with $t<0$) yields the estimate
 (4.32). Hence, $|T|<K/j$ and this completes the proof of the lemma
}

 Throughout the remainder of this section, we will use Fourier
 series. For a function $F\in\ L^2(\cir ,\C)$, we denote by $F^{(l)}$
 its
 $l$-th Fourier coefficient:
 \[
 F^{(l)}=\frac{1}{2\pi i}\int_0^{2\pi}F(t)\ei{-ilt}dt
 \]
An estimate of $\phi_1$ is given in the following lemma.
 {\lem There exist $J_0\in\Z^+$ and  positive constants $K_1,\ K_2$
 such that for $j\in\Z$, $|j|>J_0$, the function
 $\phi_1$ has the form
 \begeq
 \phi_1(t)=\phi_1^{(0)}+\frac{\phi_2(t)}{j}
 \stopeq
 with $|\phi_2(t)|\le K_1$ for every $t\in\R$, $\phi_2^{(0)}=0$, and
 $\dis\frac{1}{2}\le |\phi_1^{(0)}|\le K_2$.}

 \vspace{.2cm}

 \Proof{ We  estimate the Fourier coefficients of $\phi_1$
 in terms of those of the differentiable function
 $c_1\psi_2$, where $\psi_2$ is the function in Lemma 4.2.
 For this, we replace $\psi_1$ by $\psi_2/j$ in
  (4.14) and rewrite the first equation as
 \begeq
 \phi_1'(t)=-i\frac{N_j}{j}\phi_1(t)+c_1(t)\frac{\psi_2(t)}{j}\, .
 \stopeq
 In terms of the Fourier coefficients, we get then
 \begeq
 i\left(l+\frac{N_j}{j}\right)\phi_1^{(l)}=\frac{1}{j}
 (c_1\psi_2)^{(l)}\qquad l\in\Z
 \stopeq
 Since $N_j\sim -\dis\frac{\gamma}{\lam}$
 (by Lemma 4.1), then there exists $J_0\in\Z^+$ such that
 \[
\left|\mathrm{Re}\left(l+\frac{N_j}{j}\right)\right|\ge
\frac{|l|}{2}\qquad\forall l\in\Z^\ast
 \]
 It follows  from (4.35) that
 \begeq
 |\phi_1^{(l)}|\le \frac{2}{|j|}\left(\frac{|(c_1\psi_2)^{(l)}|}{|l|}
 \right)\qquad\forall l\in\Z^\ast
 \stopeq
 The function $\phi_2(t)=j(\phi_1(t)-\phi_1^{(0)})$ satisfies
 the conditions of the lemma.
 For $l=0$, we get
 \[
\phi_1^{(0)}=\frac{-i}{N_j}(c_1\psi_2)^{(0)}
 \]
 and since $\psi_2$ is bounded, then $|\phi_1^{(0)}|\le K_2$. Finally,
 $|\phi_1^{(0)}|\ge 1/2$ for $|j|$ large since $|\phi_2|\le K_1$
 and
 $\max|\phi_1|=1$ }

 \subsection{End of the proof of Theorem 4.1}
  Now we estimate the functions $\psi_2$ and $\phi_2$ that are
  defined in the previous lemmas

{\lem There exist $J_0\in\Z^+$ and $K>0$ such that for $j\in\Z$,
$|j|\ge J_0$, the function $\psi_2$, given in Lemma $4.2$, has the
form
 \begeq
 \psi_2(t)=\phi_1^{(0)}\frac{\ov{c_1}(t)}{i(1+\ei{ix})}+\frac{\psi_4(t)}{j}
 \stopeq
with $\psi_4$ satisfying $|\psi_4(t)|\le K$ for every $t\in\R$ }

\vspace{.2cm}

 \Proof{ With $\phi_2$ and $\psi_2$ as in Lemmas 4.2
and 4.3, we rewrite the system (4.14) and obtain
 \begeq\begar{ll}
 \phi_2'(t) & =\dis -iN_j\left(\phi_1^{(0)}+\frac{\phi_2(t)}{j}\right)
 +c_1(t)
 \psi_2(t)\\
 \psi_2'(t) & =\dis -iE_j\psi_2(t)+\ov{c_1}(t)(j\phi_1^{(0)}+\phi_2(t))
 \stopar\stopeq
 From (4.38), we see that the Fourier coefficients satisfy
 \begeq
 \psi_2^{(l)}=\phi_1^{(0)}\frac{j}{i(l+E_j)}\ov{c_1}^{(l)}+
 \frac{(\ov{c_1}\phi_2)^{(l)}}{i(l+E_j)},\qquad l\in\Z
 \stopeq
 Let
 \begeq
 \psi_3(t)=\psi_2(t)-\phi_1^{(0)}\frac{\ov{c_1}(t)}{i(1+\ei{ix})}.
 \stopeq
 It follows at once from (4.39) that the Fourier coefficients
 of $\psi_3$ satisfy
 \begeq
 \psi_3^{(l)}=\phi_1^{(0)}\frac{(1+\ei{ix})j-(l+E_j)}{i(l+E_j)(1+\ei{ix})}
 \ov{c_1}^{(l)}+\frac{(\ov{c_1}\phi_2)^{(l)}}{i(l+E_j)},\qquad
 l\in\Z
 \stopeq
Since for $|j|$ large, we have
\[
 |(1+\ei{ix})j-E_j|=|\frac{N_j}{j}-\delta|<3\quad\mathrm{and}\quad
 |l+E_j|\ge |\mathrm{Im}(E_j)|\ge \frac{j}{C_1}
\]
for some positive constant $C_1$, then from (4.41) we get
 \begeq
 |\psi_3^{(l)}|\le \frac{K}{j}
 \left(|\phi_1^{(0)}l\ov{c_1}^{(l)}|+|(\ov{c_1}\phi_2)^{(l)}|\right)
 \stopeq
The term $il\ov{c_1}^{(l)}$ is the $l$-th Fourier coefficient of
$\ov{c_1}'\in C^{k-1}$, and therefore it follows from (4.42) that
$j\psi_3$ is uniformly bounded (in $j$). Thus the function
$\psi_4=\dis\frac{\psi_3}{j}\in C^{k-1}$ is bounded
 }

{\lem There exist $J_0\in\Z^+$ and $K>0$ such that for $j\in\Z$,
$|j|\ge J_0$, the function $\phi_2$, given in Lemma $4.3$, has the
form
 \begeq
 \phi_2(t)=i\phi_1^{(0)}k(t)+\frac{\phi_4(t)}{j}
 \stopeq
where $k(t)$ is the function given in $(4.1)$, and $\phi_4\in C^k$
satisfies $|\phi_4(t)|\le K$ for every $t\in\R$ }

\vspace{.2cm}

 \Proof{ We use Lemma 4.1 to write
\[
 N_j=-\frac{\gamma}{\lam} +\frac{C_j}{j}
\]
with  $C_j\in \C$ bounded. By using Lemmas 4.3 and 4.4 in (4.34), we
obtain an equation for $\phi_2$:
 \begeq
 \phi_2'(t)=(i\frac{\gamma}{\lam}-\frac{iC_j}{j})(\phi_1^{(0)}+
 \frac{\phi_2(t)}{j})+\frac{\phi_1^{(0)}|c_1(t)|^2}{i(1+\ei{ix})}+
 \frac{c_1(t)\psi_4(t)}{j}
 \stopeq
 Equivalently,
\begeq
 \phi_2'(t)=i\phi_1^{(0)}\left[\frac{\gamma}{\lam}-
 \frac{|c_1(t)|^2}{1+\ei{ix}} +\frac{C_j}{j\phi_1^{(0)}}\right]
+\frac{1}{j}\left[i(\frac{\gamma}{\lam}-\frac{C_j}{j})\phi_2(t)+
c_1(t)\psi_4(t)\right]
 \stopeq
 Let $\phi_3(t)=\phi_2(t)-i\phi_1^{(0)}k(t)$, where $k(t)$ is the
 function defined in (4.1). Note that since
  $k'(t)=\dis\frac{\gamma}{\lam}-\frac{|c_1(t)|^2}{1+\ei{ix}}$, then
 it follows from (4.45) that
 \begeq
 \phi_3'(t)=\frac{1}{j}\left[ -iC_j+c_1(t)\psi_4(t)+i
 \left(\frac{\gamma}{\lam}-\frac{C_j}{j}\right)(\phi_3(t)+
 i\phi_1^{(0)}k(t))\right]
 \stopeq
 The Fourier coefficients of $\phi_3$ satisfy
 \begeq
 i\left[l-\frac{\gamma}{\lam j}+\frac{C_j}{j^2}\right]\phi_3^{(l)}=
 \frac{1}{j}\left[-iC_j+(c_1\psi_4)^{(l)}
  -
 \left(\frac{\gamma}{\lam}-\frac{C_j}{j}\right)\phi_1^{(0)}k^{(l)}
 \right]
 \stopeq
 This shows that $j\phi_3 $ is bounded}

 From these lemmas, we get for the functions $\phi$ and $\psi$
 \[\begar{ll}
 \phi(t) & =\dis\ei{ijt}\left(\phi_1^0+i\phi_1^0\frac{k(t)}{j}+
 \frac{\phi_4(t)}{j^2}\right)\\
 \psi(t) & =\dis\ei{ijt}\left(\phi_1^0\frac{\ov{c_1}(t)}{ij(1+\ei{ix})}+
 \frac{\psi_4(t)}{j^2}\right)\\
 \stopar\]
Since $\dis\frac{1}{2}\le\phi_1^{(0)}\le K_2$, then we can divide
the above functions by $\phi_1^{(0)}$ and obtain (4.3) and (4.4) of
the Theorem $\ \Box$

\section{The kernels}
 We use the basic solutions to construct kernels for the operator
 $\cl$.
 For $j,k\in\Z$, let
 \begeq\begar{ll}
 w_j^\pm(r,t) & =\dis
 r^{\sigpmj}\phipmj(t)+\ov{r^{\sigpmj}\psipmj(t)}\\
 w^{\ast\pm}_k(r,t) & =
 \dis r^{\mu_k^\pm}X_k^\pm(t)+\ov{r^{\mu_k^\pm}Z_k^\pm(t)}
 \stopar
 \stopeq
 be the basic solutions of $\cl$ and $\clstar$, respectively,
 with $\charac(w_j^\pm)=(\sigpmj,j)$,
 $\charac(w_k^{\ast\pm})=(\mu_k^\pm,k)$
 and such that $\phipj(0)=1$, $\phimj(0)=i$, $X_k^+(0)=i$,
 and $X_k^-(0)=1$. Note that, by Theorem 2.2, we have
  $\mu^\pm_k=-\sigma^\pm_{-j}$.

 In the remainder of this paper we will use the
 following notation, $A^\pm B^\pm =A^+B^+ +A^-B^-$. Functions $f(r,t)$
 and $g(\rho ,\ta)$ will be denoted by $f(z)$ and $g(\zeta)$, where
 $z=r^\lam\ei{it}$ and $\zeta =\rho^\lam\ei{i\ta}$. Define
 functions
 $\Om_1(z,\zeta)$ and $\Om_2(z,\zeta)$ as follows:
 \begeq
 \Om_1(z,\zeta)=\left\{\begar{ll}
 \dis\frac{1}{2}\sum_{\mathrm{Re}(\sigpmj)\ge 0}
 w_j^\pm(z)w^{\ast\pm}_{-j}(\zeta) & \qquad\mathrm{if}\
 r<\rho\\
 \dis -\frac{1}{2}\sum_{\mathrm{Re}(\sigpmj)< 0}
 w_j^\pm(z)w^{\ast\pm}_{-j}(\zeta) & \qquad\mathrm{if}\
 r>\rho
 \stopar\right.
 \stopeq
 and
 \begeq
 \Om_2(z,\zeta)=\left\{\begar{ll}
 \dis\frac{1}{2}\sum_{\mathrm{Re}(\sigpmj)\ge 0}
 \ov{w_j^\pm(z)}w^{\ast\pm}_{-j}(\zeta)
 & \qquad\mathrm{if}\
 r<\rho\\
 \dis -\frac{1}{2}\sum_{\mathrm{Re}(\sigpmj)< 0}
 \ov{w_j^\pm(z)}w^{\ast\pm}_{-j}(\zeta)
 & \qquad\mathrm{if}\
 r>\rho
 \stopar\right.
 \stopeq
Let $K(t,\ta)$ and $L(z,\zeta)$ be defined by
 \begeq
 K(t,\ta)=k(t)-k(\ta)\quad\mathrm{and}\quad
 L(z,\zeta)=\left\{\begar{ll}
 \dis\log\frac{\zeta}{\zeta -z} & \mathrm{if}\ r<\rho\\
 \dis\log\frac{z}{z -\zeta} & \mathrm{if}\ r>\rho
 \stopar\right.
 \stopeq
where $k(t)$ is the function defined in (4.1) and where $\log$
denotes the principal branch of the logarithm in $\C\backslash\R^-$.
In the next theorem, we will use the notation
\[
 \Delta_1=\{(r,t,\rho,\ta);\ 0<r\le\rho\},
 \quad
 \Delta_2=\{(r,t,\rho,\ta);\ 0<\rho\le r\}
\]
and $\mathrm{Int}(\Delta_1)$, $\mathrm{Int}(\Delta_2)$ will denote
their interiors.

{\theo The functions $C_1(z,\zeta)$ and $C_2(z,\zeta)$ defined in
$\mathrm{Int}(\Delta_1)\cup\mathrm{Int}(\Delta_2)$ by
 \begeq\begar{ll}
 C_1(z,\zeta)& =\dis
 \Om_1(z,\zeta)-i\left(\frac{r}{\rho}\right)^{\lam\nu}
 \left[\frac{\zeta}{\zeta -z}+iK(t,\ta)L(z,\zeta)\right]\\
 C_2(z,\zeta)& =\dis
 \Om_2(z,\zeta)-
 \frac{\ov{c(t)}}{2a}\left(\frac{r}{\rho}\right)^{\lam\nu}
 L(z,\zeta)-
 \frac{\ov{c(\ta)}}{2a}\ov{\left(\frac{r}{\rho}\right)^{\lam\nu}
 L(z,\zeta)}
 \stopar\stopeq
 are in $C^1(\Delta_1)\cup C^1(\Delta_2)$, meaning that
 the restrictions of $C_{1,2}$ to $\mathrm{Int}(\Delta_1)$ (or
 to $\mathrm{Int}(\Delta_2)$) extend as $C^1$ functions to
 $\Delta_1$ (or $\Delta_2$). Furthermore, for any $R>0$,
 the functions $C_{1,2}$ are bounded for $r\le R$ and $\rho\le R$}

\vspace{.2cm}

 To prove this theorem, we need two lemmas.

 \subsection{Two lemmas}

 {\lem For $|j|$ large, we have
 \begeq\begar{ll}
 w_j^\pm(z)w^{\ast\pm}_{-j}(\zeta) & =\dis
 2i\left(\frac{r}{\rho}\right)^{\sigma_j}\ei{ij(t-\ta)}
 \left( 1 +i\frac{K(t,\ta)}{j}+O(j^{-2})\right)\\
  \ov{w_j^\pm(z)}w^{\ast\pm}_{-j}(\zeta) & =\dis
  \left(\frac{r}{\rho}\right)^{\sigma_j}\frac{\ov{c(t)}}{aj}
  \ei{ij(t-\ta)}+
  \left(\frac{r}{\rho}\right)^{\ov{\sigma_j}}
  \frac{\ov{c(\ta)}}{aj}
  \ei{-ij(t-\ta)}+O(j^{-2})
 \stopar\stopeq}

\vspace{.2cm}

 \Proof{ It follows from Theorem 4.1 that for large
$|j|$, the corresponding spectral values are in $\C\backslash\R$
(for $b\ne 0$). We have,
 $\sigma_j^-=\sigma_j^+=\sig_j$ and $\sig_j$  is given by (4.2).
 Furthermore, it follows from section 2 that
 \[
 w_j^+(r,t)=r^{\sigma_j}\phi_j(t)+\ov{r^{\sigma_j}\psi_j(t)}\quad
 \mathrm{and}\quad
 w_j^-(r,t)=i(r^{\sigma_j}\phi_j(t)-\ov{r^{\sigma_j}\psi_j(t)})
 \]
 with $\phi_j(0)=1$. For the basic solutions of the adjoint operator,
 we have (from Theorem 2.2) that
 \[\begar{l}
 w_{-j}^{\ast +}(\rho,\ta)=i(\rho^{-\sigma_j}X_{-j}(\ta)-
 \ov{\rho^{-\sigma_j}Z_{-j}(\ta)})\quad
 \mathrm{and}\\
 w_{-j}^{\ast -}(\rho,\ta)=\rho^{-\sigma_j}X_{-j}(\ta)+
 \ov{\rho^{-\sigma_j}Z_{-j}(\ta)}
 \stopar\]
 with $X_{-j}(0)=1$. Hence
 \begeq\begar{ll}
 w_j^\pm(z)w^{\ast\pm}_{-j}(\zeta) & =\dis
 2i\left(\left(\frac{r}{\rho}\right)^{\sigma_j}\phi_j(t)X_{-j}(\ta)
 -\ov{\left(\frac{r}{\rho}\right)^{\sigma_j}\psi_j(t)Z_{-j}(\ta)}
  \right)\\
 \ov{w_j^\pm(z)}w^{\ast\pm}_{-j}(\zeta) & =\dis
 2i\left(\left(\frac{r}{\rho}\right)^{\sigma_j}\psi_j(t)X_{-j}(\ta)
 -\ov{\left(\frac{r}{\rho}\right)^{\sigma_j}\phi_j(t)Z_{-j}(\ta)}
  \right)
 \stopar\stopeq
 Now, the asymptotic expansions (4.3) and (4.4) give the
 following products
 \begeq\begar{ll}
 \phi_j(t)X_{-j}(\ta) & =\dis
 \ei{ij(t-\ta)}\left(1+i\frac{K(t,\ta)}{j}+O(j^{-2})\right)\\
 \psi_j(t)X_{-j}(\ta) & =\dis -i\frac{\ov{c(t)}}{2aj}
 \ei{ij(t-\ta)}+O(j^{-2})\\
 \phi_j(t)Z_{-j}(\ta) & =\dis -i\frac{c(\ta)}{2aj}
 \ei{ij(t-\ta)}+O(j^{-2})\\
 \psi_j(t)Z_{-j}(\ta) & =\dis O(j^{-2})
 \stopar\stopeq
 Estimates (5.6) of the lemma follow from (5.7) and (5.8) }

 {\lem For $j\in\Z^+$ large and $\sig_j$ as in $(4.2)$, consider the
 function
 \[
 f_j(t)=t^{\sig_j}-t^{\lam(j+\nu)},\qquad 0<t<1\, .
 \]
 Then there are $J_0>0$ and $C>0$ such that
 \[
 |f_j(t)|\le \frac{C}{j^2},\qquad \forall t\in  (0,\  1),\ j\ge J_0
 \]}

\vspace{.2cm}

 \Proof{ By using the asymptotic expansion for $\sig_j$ given
 in (4.2) and $\lam =a+ib$ ($a>0$), we write
 \[
 \sig_j=\left[a(j+\nu)+\frac{\alpha}{j}\right]+i
 \left[b(j+\nu)+\frac{2\beta}{j^2}\right]
 \]
 where $\alpha >0$ and $\beta\in\R$, depend on $j$, but are bounded.
 Hence,
 \[
f_j(t)=t^{a(j+\nu)+(\alpha /j)} t^{i[b(j+\nu)+2(\beta /j^2)]}-
t^{a(j+\nu)}t^{ib(j+\nu)}
 \]
 We decompose $f_j$ as $f_j=g_j+h_j$ with
 \[\begar{ll}
 g_j(t) & =\dis t^{a(j+\nu)}
  t^{i(b(j+\nu)+2(\beta /j^2))}\left( t^{\alpha /j}-1\right)\\
 h_j(t) & =t^{a(j+\nu)}t^{ib(j+\nu)}\left(
  t^{2i\beta /j^2}-1\right)
 \stopar\]
It is enough to verify that both $|g|$ and $|h|$ are $O(1/j^2)$.
Since $a>0$ and $j$ large, then we can assume that $g$ and $h$ are
defined at 0 and that $g(0)=h(0)=g(1)=h(1)=0$.

For the function $g$, we have
 $|g(t)|=t^{a(j+\nu)}-t^{a(j+\nu)+(\alpha /j)}$.
The maximum of $|g|$ occurs at the point
\[
t_\ast =\left(\frac{a(j+\nu)}{a(j+\nu)+(\alpha
/j)}\right)^{j/\alpha}
\]
and
 \[
|g(t_\ast)|=t_\ast^{a(j+\nu)}\left(
 1-\frac{a(j+\nu)}{a(j+\nu)+(\alpha /j)}\right) \le
 \frac{\alpha}{j(a(j+\nu)+(\alpha /j))}
 \le\frac{A_1}{j^2} \]
 for some positive constant $A_1$.

 For the function $h$, note that if $\beta =0$, then $h=0$.
 So assume that $\beta\ne 0$.  We have
 \[
 |h(t)|^2=t^{2a(j+\nu)}\left|t^{i\beta /j^2}-
 t^{-i\beta /j^2}\right|^2 =4t^{2a(j+\nu)}
 \sin^2\left(\frac{\beta\ln t}{j^2}\right)
 \]
 For $0<t\le 1/2$, we get
 \[
 |h(t)|\le 2t^{a(j+\nu)}\le \frac{2}{2^{a(j+\nu)}}\le
 \frac{A_2}{j^2}
 \]
for some $A_2>0$. For $1\ge t\ge (1/2)$, we have
\[
\frac{d}{dt}(|h(t)|^2)= 8t^{a(j+\nu)-1}
 \sin(\frac{\beta\ln t}{j^2})
 \left[ a(j+\nu)\sin(\frac{\beta\ln t}{j^2})
 +\frac{\beta}{j^2}\cos(\frac{\beta\ln t}{j^2})
 \right]\, .
\]
  For $j$ sufficiently large, the critical points of
  $|h|^2$ in the interval $(1/2\, ,\ 1)$ are  the solutions of the
  equation
  \[
 \tan\left(\frac{\beta\ln t}{j^2}\right)=-\frac{\beta}{aj^2(j+\nu)}.
  \]
Hence,
\[
 \frac{\beta\ln t}{j^2}=-\arctan( \frac{\beta}{aj^2(j+\nu)})
 +k\pi,\qquad k\in\Z.
\]
 However, since $1\ge t\ge 1/2$ and $j$ is large, the only
 possible value of the integer $k$ is $k=0$. Hence, $|h|$
 has a single critical point in $(1/2,\ 1)$:
 \[
 t_\ast =\exp\left(-\frac{j^2}{\beta}
 \arctan( \frac{\beta}{aj^2(j+\nu)})\right).
 \]
 The function $|h|$ has a maximum at $t_\ast$ and
 \[
|h(t_\ast)|^2=4t_\ast^{2a(j+\nu)}\sin^2(\arctan(
\frac{\beta}{aj^2(j+\nu)}))\le 4\arctan^2(
\frac{\beta}{aj^2(j+\nu)})\le \frac{A_3}{j^6}
 \]
 for some $A_3>0$}

\subsection{Proof of Theorem 5.1}
We use the series expansions
\[\begar{ll}
 \dis\frac{\zeta}{\zeta -z} & =\left\{\begar{ll}\dis
 \sum_{j\ge 0}(r/\rho)^{\lam j}\ei{ij(t-\ta)} &
 \mathrm{if}\ r<\rho \\
 \dis -\sum_{j\ge 1}(r/\rho)^{-\lam j}\ei{-ij(t-\ta)} &
 \mathrm{if}\ r<\rho
  \stopar\right. \\
 \dis L(z,\zeta ) & =\left\{\begar{ll}\dis
 \sum_{j\ge 1}(r/\rho)^{\lam j}\frac{\ei{ij(t-\ta)}}{j} &
 \mathrm{if}\ r<\rho \\
 \dis \sum_{j\ge 1}(r/\rho)^{-\lam j}
 \frac{\ei{-ij(t-\ta)}}{j} &
 \mathrm{if}\ r<\rho
  \stopar\right.
 \stopar\]
 together with (5.2) and (5.6) to decompose $C_1(z,\zeta)$ as
 follows. For a large integer $J_0$ and $r<\rho$,
 \begeq
 C_1=P_1+i\sum_{j\ge J_0}
 \left[\left(\frac{r}{\rho}\right)^{\sigma_j}-
 \left(\frac{r}{\rho}\right)^{\lam(j+\nu)}\right]
 \ei{ij(t-\ta)}\left(1+i\frac{K}{j}+O(j^{-2})\right)
 \stopeq
 where $P_1(z,\zeta)$ consists of the finite collection of terms in the
 series with
 index $j<J_0$. Thus $P\in C^1(\Delta_1)$. The second term
 ($\sum_{j\ge J_0}\cdots$) on the right of (5.9) is also in
 $C^1(\Delta_1)$ since
 \[
 \left|\left(\frac{r}{\rho}\right)^{\sigma_j}-
 \left(\frac{r}{\rho}\right)^{\lam(j+\nu)}\right| =O(j^{-2})
 \]
 by Lemma 5.2. When $r>\rho$, the decomposition of $C_1$ takes the
 form
 \begeq
 C_1=\widetilde{P}_1+i\sum_{j\ge J_0}
 \left[\left(\frac{r}{\rho}\right)^{-\lam(j-\nu)}-
 \left(\frac{r}{\rho}\right)^{\sigma_{-j}}\right]
 \ei{-ij(t-\ta)}\left(1-i\frac{K}{j}+O(1/j^2)\right)
 \stopeq
 As before, the finite sum $\widetilde{P}_1(z,\zeta)\in
 C^1(\Delta_2)$. Since $\sig_{-j}=\lam(-j+\nu)-(\gamma /j)+O(j^{-2})$
 and $r>\rho$, we have
 \[
 \left|\left(\frac{r}{\rho}\right)^{-\lam(j-\nu)}-
 \left(\frac{r}{\rho}\right)^{\sigma_{-j}}\right|=
 \left|\left(\frac{\rho}{r}\right)^{\lam(j-\nu)}-
 \left(\frac{\rho}{r}\right)^{-\sigma_{-j}}\right|
 =O(j^{-2}).
 \]
 Again, the infinite sum on the right of (5.10) is in
 $C^1(\Delta_2)$. This proves the theorem for the function $C_1$.
 Similar arguments can be used for the function $C_2$ $\quad\Box$

\subsection{Modified kernels}
 The following modifications to the kernels $\Om_1$ and $\Om_2$ will
 be used to establish a similarity principle in section 8.
 For $j_0\in\Z$, we define $\Om_{j_0,1}^\pm (z,\zeta)$ and
$\Om_{j_0,2}^\pm (z,\zeta)$  by
 \begeq
 \Om_{j_0,1}^\pm(z,\zeta)=\left\{\begar{ll}
 \dis\frac{1}{2}\sum_{\mathrm{Re}(\sigpmj)\ge \mathrm{Re}(\sig_{j_0}^\pm)}
 w_j^\pm(z)w^{\ast\pm}_{-j}(\zeta) & \qquad\mathrm{if}\
 r<\rho\\
 \dis -\frac{1}{2}\sum_{\mathrm{Re}(\sigpmj)< \mathrm{Re}(\sig_{j_0}^\pm)}
 w_j^\pm(z)w^{\ast\pm}_{-j}(\zeta) & \qquad\mathrm{if}\
 r>\rho
 \stopar\right.
 \stopeq
 and
 \begeq
 \Om_{j_0,2}^\pm (z,\zeta)=\left\{\begar{ll}
 \dis\frac{1}{2}\sum_{\mathrm{Re}(\sigpmj)\ge \mathrm{Re}(\sig_{j_0}^\pm)}
 \ov{w_j^\pm(z)}w^{\ast\pm}_{-j}(\zeta)
 & \qquad\mathrm{if}\
 r<\rho\\
 \dis -\frac{1}{2}\sum_{\mathrm{Re}(\sigpmj)< \mathrm{Re}(\sig_{j_0}^\pm)}
 \ov{w_j^\pm(z)}w^{\ast\pm}_{-j}(\zeta)
 & \qquad\mathrm{if}\
 r>\rho
 \stopar\right.
 \stopeq

 {\theo The functions $C_1(z,\zeta)$ and $C_2(z,\zeta)$ given by
 \begeq
 \Om_{j_0,1}^\pm(z,\zeta) =i\left(\frac{r}{\rho}\right)^{\sig_{j_0}^\pm}
 \ei{ij_0(t-\ta)}\left[
 \frac{\zeta}{\zeta -z}+iK(t,\ta)L(z,\zeta)+C_1(z,\zeta)\right]
 \stopeq
 and
 \begeq\begar{ll}
\Om_{j_0,2}^\pm(z,\zeta)
  & =\dis\frac{\ov{c(t)}}{2a}\left(\frac{r}{\rho}\right)^{\sig_{j_0}^\pm}
 L(z,\zeta)+
 \ov{ \frac{c(\ta)}{2a}\left(\frac{r}{\rho}\right)^{\sig_{j_0}^\pm}
 L(z,\zeta)}\\
 &\dis \qquad\qquad\qquad\qquad
 +\left(\frac{r}{\rho}\right)^{\sig_{j_0}^\pm}C_2(z,\zeta)
 \stopar\stopeq
 have the following properties:
 \begin{enumerate}
 \item $C_{1,2}\in
 C^1(\mathrm{Int}(\Delta_1)\cup\mathrm{Int}(\Delta_2))$,
 \item
 for a given $z=r^\lam\ei{it}$ with $0<r<R$, the functions $C_{1,2}(z,\,
 \cdot) $ are in $L^p(\{(\rho ,\ta);\ \rho <R\})$, for every $p>0$, and
 \item the functions $|z-\zeta|^\ep C_{1,2}$ are bounded
 in the region $r<R$ and $\rho <R$, for any $\ep >0$.
 \end{enumerate}}

\vspace{.2cm}

 \Proof{ The proof follows similar arguments as those
used in the proof of Theorem 5.1. We describe briefly how the
properties of $C_1$ can be established in the region $r<\rho$. We
write
 \[
 \Om_{j_0,1}^\pm -i\left(\frac{r}{\rho}\right)^{\sig_{j_0}^\pm}
 \ei{ij_0(t-\ta)}\left[
 \frac{\zeta}{\zeta -z}+iK(t,\ta)L(z,\zeta)\right] =P_1(z,\zeta)+
 \sum_{j\ge J_0}
 \]
where $P_1$ is the finite sum consisting of terms in the series of
$\Om_{j_0,1}^{\pm}$ with indices $j\le J_0$ and
$\mathrm{Re}(\sig_j^\pm)\ge \mathrm{Re}(\sig_{j_0}^\pm)$, and the
terms with indices $j\le J_0$ in the series expansions of $\zeta
/(\zeta -z)$ and $L(z,\zeta)$.  The infinite sum
 $\dis\sum_{j\ge
J_0}$ can be written as
 \[
 \sum_{j\ge J_0}=
 i\left(\frac{r}{\rho}\right)^{\sig_{j_0}^\pm}
 \sum_{j\ge J_0} A_j(z,\zeta)
 \]
 where
 \[
 A_j(z,\zeta)= \left[\left(\frac{r}{\rho}\right)^{\sig_j-\sig_{j_0}^\pm}
 -\left(\frac{r}{\rho}\right)^{\lam(j-j_0)} \right]
 \ei{ij(t-\ta)}\left(1+i\frac{K(t,\ta)}{j}+O(j^{-2})\right)
 \]
 Since $\sig_j$ satisfies the asymptotic expansion $(4.2)$ and
 since $r<\rho$,  arguments similar to those used in the proof
 of Lemma 5.2 show that
 \[
\left|\left(\frac{r}{\rho}\right)^{\sig_j-\sig_{j_0}^\pm}
 -\left(\frac{r}{\rho}\right)^{\lam(j-j_0)} \right| =
 O\left(\frac{1}{j}\right).
 \]
 Thus  $\dis\sum_{j\ge J_0}$ has the desired properties
 of the theorem in $\Delta_1$. Analogous arguments can be used
 in the region $r>\rho$ and also for the function $C_2$}

 \section{The homogeneous equation $\cl u=0$}
 In this section, we use the kernels defined in section 5 to obtain
 series and integral representations of the solutions of the
 equation $\cl u=0$. Versions of the Laurent series expansion, in
 terms of the basic solutions, and the Cauchy integral formula are
 derived. Some consequences of these representations are given.

 \subsection{Representation of solutions in a cylinder}
For $R,\delta \in \R^+$ with $\delta <R$, consider the cylinder
$A(\delta,R)= (\delta ,\ R)\times \cir$. Again, let
$z=r^\lam\ei{it}$,  $\ \zeta =\rho^\lam\ei{i\ta}$ and $\Om_1$,
$\Om_2$  denote the functions defined in section 5. We have the
following theorem.

{\theo Let $u\in C^0(\ov{A(\delta ,R)})$ be a solution of $\cl u=0$.
Then
 \begeq
 u(z)=\frac{-1}{2\pi}\int_{\pa A(\delta ,R)}
 \frac{\Om_1(z,\zeta)}{\zeta}u(\zeta)d\zeta +
 \frac{\ov{\Om_2(z,\zeta)}}{\ov{\zeta}}\,\ov{u(\zeta)}\,\ov{d\zeta}
 \stopeq}

\vspace{.2cm}

 \Proof{ Let
 \[
 L=\lam\dd{}{t}-ir\dd{}{r} \quad\mathrm{and}\quad
 L^\ast =\lam\dd{}{\ta}-i\rho\dd{}{\rho},
 \]
so that $\cl u= Lu +i\lam\nu u -c(t)\ov{u}$ and $-\clstar v =L^\ast
v -i\lam\nu v+\ov{c(\ta)}\ov{v}$. It follows from the definitions of
the kernels $\Om_{1,2}$ given in (5.2) and (5.3) and from the fact
that $\cl w_j^\pm=0$ and $\clstar w_k^{\ast \pm}=0$ that the kernels
satisfy the following relations
 \begeq\begar{l}\\
 L^\ast\Om_1(z,\zeta)
 =i\lam\nu\Om_1(z,\zeta)-\ov{c(\ta)}\,\ov{\Om_2(z,\zeta)}\\
 L^\ast\Om_2(z,\zeta)
 =i\lam\nu\Om_2(z,\zeta)-\ov{c(\ta)}\,\ov{\Om_1(z,\zeta)}
 \stopar\stopeq
 Consider the functions $P_1=\Om_1+\Om_2$ and $P_2=-i\Om_1+i\Om_2$.
 Then (6.2) implies that
 \begeq
 \clstar P_1=\clstar P_2 =0\, .
 \stopeq
 Let $(r_0,t_0)\in A(\delta ,R)$ and $z_0=r_0^\lam\ei{it_0}$.
 For $\ep >0$, let
 \[
 D_\ep =\{ (\rho ,\ta)\in \R^+\times\cir ,\ \ |\zeta -z_0|<\ep\}\, .
 \]
 Hence, for $\ep$ small, $D_\ep$ is diffeomorphic to the disc and is
 contained in $A(\delta ,R)$. We apply Green's identity (1.8)
 in the region $A(\delta ,R)\backslash D_\ep$  to each
 pair $u(\zeta),\ P_k(z_0,\zeta)$,   with $k=1,2$. Since, $\cl u=0$
 and $\clstar P_k=0$, then
 \[\begar{l}
\dis\mathrm{Re}\left[\int_{\pa A}\!\!\!\!\!
P_k(z_0,\zeta)u(\zeta)\frac{d\zeta}{\zeta}+
\ov{P_k(z_0,\zeta)}\,\ov{u(\zeta)}\,\frac{d\ov{\zeta}}{\ov{\zeta}}\right]
= \\
\dis \qquad\qquad\qquad\qquad\qquad \mathrm{Re}\left[\int_{\pa
D_\ep}\!\!\!\!\! P_k(z_0,\zeta)u(\zeta)\frac{d\zeta}{\zeta}+
\ov{P_k(z_0,\zeta)}\,\ov{u(\zeta)}\,\frac{d\ov{\zeta}}{\ov{\zeta}}\right]
 \stopar\]
 Then, after multiplying by $i$ the above identity with $k=2$ and
 adding it to the identity with $k=1$,
   we obtain
 \[
 \int_{\pa A}\!\!\!\!\! (P_1+iP_2)u\frac{d\zeta}{\zeta}+
 (\ov{P_1}+i\ov{P_2})\ov{u}\frac{d\ov{\zeta}}{\ov{\zeta}}
=
  \int_{\pa D_\ep}\!\!\!\!\! (P_1+iP_2)u\frac{d\zeta}{\zeta}+
 (\ov{P_1}+i\ov{P_2})\ov{u}\frac{d\ov{\zeta}}{\ov{\zeta}}
 \]
 Since $2\Om_1=P_1+iP_2$ and $2\Om_2=P_1-iP_2$, we get
 \begeq\begar{l}\dis
\int_{\pa A}\!\!\!\!\! \Om_1(z_0,\zeta)u(\zeta)\frac{d\zeta}{\zeta}+
 \ov{\Om_2(z_0,\zeta)}\,\ov{u(\zeta)}\,\frac{d\ov{\zeta}}{\ov{\zeta}}
=\\ \qquad\qquad\qquad\qquad\qquad\dis
 \int_{\pa D_\ep}\!\!\!\!\! \Om_1(z_0,\zeta)u(\zeta)\frac{d\zeta}{\zeta}+
 \ov{\Om_2(z_0,\zeta)}\,\ov{u(\zeta)}\,\frac{d\ov{\zeta}}{\ov{\zeta}}
 \stopar \stopeq
 Now, we let $\ep\to 0$ in the right side of (6.4). From the
 estimates (5.5) of the kernels, it follows that the only term that provides
 a nonzero contribution (as $\ep\to 0$) is the term containing
 $\zeta /(\zeta -z)$ since $C_1$, $C_2$ are bounded and $L(z_0,\zeta)$
 has a logarithmic growth. That is,
 \[
 \lim_{\ep\to 0}\int_{\pa D_\ep}\!\!\!\!\!
 \Om_1u\frac{d\zeta}{\zeta}+
 \ov{\Om_2}\,\ov{u}\,\frac{d\ov{\zeta}}{\ov{\zeta}}=
 \lim_{\ep\to 0}\int_{\pa D_\ep}\!\!\!\!\! i(r_0/\rho)^{\lam\nu}
 \frac{u(\zeta)}{\zeta -z_0}d\zeta =-2\pi u(z_0)
 \]
 This proves the Theorem }

 {\theo Suppose that $u$ is a solution of $\cl u=0$ in the cylinder $A(\delta
 ,R)$ with $0\le \delta \le R\le \infty$. Then, $u$ has the Laurent
 series expansion
 \begeq
 u(r,t)=\sum_{j\in\Z}a_j^\pm w_j^\pm(r,t)
 \stopeq
where $a_j^\pm\in\R$ are given by
 \begeq
 a_j^\pm =\frac{-1}{2\pi}\mathrm{Re}\int_0^{2\pi}w^{\ast\pm}_{-j}
  (R_0,\ta)
 u(R_0,\ta)id\ta
 \stopeq
where $R_0$ is any point in $(\delta ,\ R)$. Furthermore, there
exists $C>0$ such that
 \begeq
 |a_j^\pm|\le \frac{C}{R_0^{\mathrm{Re}(\sig_j^\pm)}}
 \max_\ta |u(R_0,\ta)|\qquad\forall j\in\Z
 \stopeq}

\vspace{.2cm}

\Proof{ Let $R_0\in (\delta ,\ R)$, and $\delta_1, R_1$ be such that
 $\delta <\delta_1<R_0<R_1<R$. For $r\in (\delta_1,\ R_1)$, we apply
 the integral representation (6.1) in the cylinder $A(\delta_1,R_1)$
 to get
 \begeq -2\pi u(r,t)=I_1-I_2\stopeq
 where
 \[\begar{ll}
 I_1 & =\dis\int_{\rho =R_1}\!\!\!\Om_1(z,\zeta)u(\zeta)\frac{d\zeta}{\zeta}
 +\ov{\Om_2(z,\zeta)}\,\ov{u(\zeta)}\,\frac{d\ov{\zeta}}{\ov{\zeta}}\\
 I_2 & =\dis\int_{\rho =\delta_1}\!\!\!\Om_1(z,\zeta)u(\zeta)\frac{d\zeta}{\zeta}
 +\ov{\Om_2(z,\zeta)}\,\ov{u(\zeta)}\,\frac{d\ov{\zeta}}{\ov{\zeta}}
 \stopar\]
 The series (5.2) and (5.3) for $\Om_1$ and $\Om_2$  give
 \[
 I_1=\sum_{\mathrm{Re}(\sig_j^\pm)\ge 0}
 w_j^\pm(r,t)\mathrm{Re}\int_{\rho =R_1}w_{-j}^{\ast \pm}(\zeta)u(\zeta)
 \frac{d\zeta}{\zeta}\, .
 \]
 Since $\cl u=0$ and $\clstar w_{j}^{\ast \pm}=0$, then Green's
 identity gives
 \[
\mathrm{Re}\int_{\rho =R_1}w_{-j}^{\ast \pm}(\zeta)u(\zeta)
 \frac{d\zeta}{\zeta} =\mathrm{Re}\int_0^{2\pi}w_{-j}^{\ast \pm}
 (R_0,\ta)u(R_0,\ta)id\ta\, .
 \]
 Hence,
 \[
 I_1=-2\pi\sum_{\mathrm{Re}(\sig_j^\pm)\ge 0}a_j^\pm w_j^\pm (r,t)
 \]
 where $a_j^\pm$ is given by (6.6). A similar calculation
 shows that
 \[
 I_2 =2\pi\sum_{\mathrm{Re}(\sig_j^\pm)< 0}a_j^\pm w_j^\pm (r,t)\, .
 \]
 To estimate the coefficients $a_j^\pm$, recall that
 \[
 w_{-j}^{\ast \pm}(R_0,\ta) =
 R_0^{-\sig_j^\pm}X_{-j}^\pm (\ta)+
 \ov{R_0^{-\sig_j^\pm}Z_{-j}^\pm
 (\ta)}\, .
 \]
 Thus
 \[
 |w_{-j}^{\ast \pm}(R_0,\ta)| \le \frac{1}{
 R_0^{\mathrm{Re}(\sig_j^\pm)}}\left(
 |X_{-j}^\pm (\ta)|+|Z_{-j}^\pm
 (\ta)|\right)\, \le \frac{C}{ R_0^{\mathrm{Re}(\sig_j^\pm)}}
 \]
 where
 $C=\dis\sup_{k,\ta}\left(|X_k^\pm(\ta)|+|Z_k^\pm(\ta)|\right)$.
 This gives estimate (6.7)}

The following theorem is a direct consequence of Theorem 6.2.

 {\theo Let $u$ be a bounded solution of $\cl u=0$ in the
cylinder $A(0 ,R)$. Then $u$ has the
 series expansion
 \begeq
 u(r,t)=\sum_{\mathrm{Re}(\sig_j^\pm)\ge 0}a_j^\pm w_j^\pm(r,t)
 \stopeq
 where $a_j^\pm$ are given by (6.6). If, in addition, $u$ is continuous
 on $\ov{A(0,R)}$, then the above summation is taken over the
 spectral values $\sig_j^\pm$ satisfying $\mathrm{Re}(\sig_j^\pm)>0$
 or $\sig_j^\pm =0$. }

 \subsection{Cauchy integral formula}

For a subset $U\subset \R\times\cir$, we set $\pa_0U=\pa U\backslash
S_0$, where  $S_0=\{0\}\times\cir$. We have the following
 integral representation that generalizes the classical Cauchy
 integral formula.

 {\theo Let $U$ be an open and bounded subset of $\R^+\times\cir$  such that
 $\pa U$ consists of finitely many simple closed and piecewise smooth
 curves. Let $u\in C^0(\ov{U}\backslash S_0)$ be such that $\cl u=0$.
 Then, for $(r,t)\in U$, we have
 \begeq
 u(r,t)=\frac{-1}{2\pi}\int_{\pa_0U}\!\!\!\Om_1(z,\zeta)u(\zeta)
 \frac{d\zeta}{\zeta}+\ov{\Om_2(z,\zeta)}\, \ov{u(\zeta)}\,
 \frac{d\ov{\zeta}}{\ov{\zeta}}
 \stopeq
  }

 \vspace{.2cm}

  \Proof{ For $\delta >0$, define
  $U_\delta =U\backslash A(0,\delta)$. Let $(r_0,t_0)\in U$.
  Choose $\ep >0$ and $\delta >0$ small enough so that
  $(r_0,t_0)\in U_\delta$ and $D_\ep\subset U_\delta$, where
  $D_\ep =\{(\rho,\ta);\ |\zeta -z_0|<\ep\}$.
  Arguments similar to those used in the proof of Theorem 6.1 show
  that
  \begeq
 u(r_0,t_0)=\frac{-1}{2\pi}\int_{\pa U_\delta}\!\!\!
 \Om_1(z,\zeta)u(\zeta)
 \frac{d\zeta}{\zeta}+\ov{\Om_2(z,\zeta)}\, \ov{u(\zeta)}\,
 \frac{d\ov{\zeta}}{\ov{\zeta}}\, .
  \stopeq
 If $S_0\cap\pa U=\emptyset$, then for $\delta$ small enough
 $U_\delta =U$ and the theorem is proved in this case.
 If $S_0\cap\pa U\ne\emptyset$, let $\Gamma_\delta =\pa
 U_\delta\cap\left(\{\delta\}\times \cir\right)$. That is, $\Gamma_\delta$
 is the part of $\pa U_\delta$ contained in the circle $r=\delta$.
 Going back to the definition of the kernels, we obtain
 \begeq
 \int_{\Gamma_\delta}\!\!\!
 \Om_1(z,\zeta)u(\zeta)
 \frac{d\zeta}{\zeta}+\ov{\Om_2(z,\zeta)}\, \ov{u(\zeta)}\,
 \frac{d\ov{\zeta}}{\ov{\zeta}} =-\!\!\!\!\!
 \sum_{\mathrm{Re}(\sig_j^\pm)<0}\!\!\!\!\!
 w_j^\pm(r_0,t_0)\mathrm{Re}\int_{\Gamma_\delta}\!\!\!\!\!
 w_{-j}^{\ast \pm}(\zeta)u(\zeta)\frac{d\zeta}{\zeta}
 \stopeq
 Since there exists $C>0$ such that
 \[
 | w_{-j}^{\ast \pm}(\rho,\ta)|\le
 C\rho^{-\mathrm{Re}(\sig_j^\pm)}\qquad\forall j\in \Z
 \]
 then, it follows from (6.12), that
 \[
 \left|\int_{\Gamma_\delta}\!\!\!
 \Om_1(z,\zeta)u(\zeta)
 \frac{d\zeta}{\zeta}+\ov{\Om_2(z,\zeta)}\, \ov{u(\zeta)}\,
 \frac{d\ov{\zeta}}{\ov{\zeta}} \right|\le 2\pi C ||u||_0\!\!\!\!\!
 \sum_{\mathrm{Re}(\sig_j^\pm)<0}\!\!\!\!\!
 \rho^{-\mathrm{Re}(\sig_j^\pm)}
 |w_j^\pm(r_0,t_0)|
 \]
Therefore,
\[
\lim_{\delta\to 0}\int_{\Gamma_\delta}\!\!\!
 \Om_1(z,\zeta)u(\zeta)
 \frac{d\zeta}{\zeta}+\ov{\Om_2(z,\zeta)}\, \ov{u(\zeta)}\,
 \frac{d\ov{\zeta}}{\ov{\zeta}} =0
\]
 and (6.10) follows from (6.11) when we let $\delta\to 0$}

 The following theorem extends the Cauchy integral formula to
 include the points on the characteristic circle $S_0$, when the
 solution is continuous up to the boundary.

{\theo Suppose that $\cl$ has no spectral values in $i\R^\ast$. Let
$U$ be an open, bounded subset of $\R^+\times\cir$, such that
 $\pa U$ consists of finitely many simple closed and piecewise smooth
 curves and with $S_0\subset\pa U$.
 Let $u\in C^0(\ov{U})$ be such that $\cl u=0$. Then,
 for $(r,t)\in U\cup S_0$, we have
 \begeq
 u(r,t)=\frac{-1}{2\pi}\int_{\pa_0U}\!\!\!\Om_1(z,\zeta)u(\zeta)
 \frac{d\zeta}{\zeta}+\ov{\Om_2(z,\zeta)}\, \ov{u(\zeta)}\,
 \frac{d\ov{\zeta}}{\ov{\zeta}}
 \stopeq}

 \vspace{.2cm}

\noindent{\textit{Proof.}} We know from Theorem 6.4 that (6.13)
holds for $r>0$. We need to verify that it holds at the points
$(0,t)\in S_0$. Since $\cl$ has no purely imaginary spectral values,
then $\Om_1$ and $\Om_2$ are well defined on $S_0$. We have
$\Om_1(0,t)=\Om_2(0,t)=0$ if 0 is not a spectral value of $\cl$ and
if 0 is a spectral value, with say multiplicity 2, then
 \[\begar{ll}
 \Om_1(0,t,\rho,\ta) & =f_{j_0}^+(t)g_{-j_0}^+(\ta)+
 f_{j_0}^-(t)g_{-j_0}^-(\ta)\\
 \Om_2(0,t,\rho,\ta) & =\ov{f_{j_0}^+(t)}g_{-j_0}^+(\ta)+
 \ov{f_{j_0}^-(t)}g_{-j_0}^-(\ta)
 \stopar\]
where $f_{j_0}^\pm(t)$ are the basic solutions of $\cl$ with
exponent 0 and $g_{-j_0}^\pm(t)$ are the basic solutions of
$\clstar$ with exponent 0.

When 0 is not a spectral value, (6.13) holds for $r=0$ by letting
$r\to 0$ in (6.10). In this case $u\equiv 0$ on $S_0$. When 0 is
spectral value, then since $S_0\subset\pa U$,  $\, u(r,t)$ has a
Laurent series expansion in a cylinder $A(0,\delta)\subset U$. In
particular,
\[\begar{ll}
u(0,t) & =a^+f_{j_0}^+(t)+a^-f_{j_0}^-(t)\\
 & =\dis\frac{-1}{2\pi}\int_{\pa_0U}\!\!\!\Om_1(0,t,\zeta)u(\zeta)
 \frac{d\zeta}{\zeta}+\ov{\Om_2(0,t,\zeta)}\, \ov{u(\zeta)}\,
 \frac{d\ov{\zeta}}{\ov{\zeta}}\quad\Box
\stopar\]

\vspace{.2cm}

\Rem{6.1} {It follows from this Theorem that if $\cl$ has no
spectral values on $i\R^\ast$, then all bounded solutions of $\cl
u=0$ in a cylinder $A(0,R)$ are continuous up to $S_0$. If $\cl$ has
spectral values on $i\R^\ast$, then there are bounded solutions on
$A(0,R)$ that are not continuous up to $S_0$. In fact, the basic
solution $r^{i\tau}\phi(t)+r^{-i\tau}\psi(t)$ is such a solution
when $i\tau\in\specl$. Note also that the number of spectral values
on $i\R^\ast$ is at most finite (this follows from the asymptotic
expansion of $\sig_j$).}

\subsection{Consequences}
 We give here some consequences of the above representations
 theorems. First we define the order of a solution along $S_0$.
We say that a solution, $u$, of $\cl u=0$ in a cylinder $A(0,R)$ has
a zero or a pole of order $s=\mathrm{Re}(\sig_{j_0})$ (with
$\sig_{j_0}=\sig_{j_0}^+$ or $\sig_{j_0}^-$) along the circle $S_0$
if, in the Laurent series expansion of $u$, all the coefficients
$a_j^\pm$ corresponding to $\mathrm{Re}(\sig_j^\pm)<s$ are zero.
That is
\[
u(r,t)=\sum_{\mathrm{Re}(\sigpmj)\ge s}a_j^\pm w_j^\pm(r,t)\, .
\]
We have the following uniqueness result.

{\theo Suppose that the spectral values of $\cl$ satisfy the
following condition
 \begeq
\mathrm{Re}(\sig_j)=\mathrm{Re}(\sig_k)\,\Longrightarrow\,
\sig_j=\sig_k, \qquad\forall\sig_j,\sig_k\in\specl\, .
 \stopeq
Let $u$ be a solution of $\cl u=0$ in a cylinder $A(0,R)$. Suppose
that $u$ is of finite order along $S_0$ and that there is a sequence
of points $(r_k,t_k)\in A(0,R)$ such that $r_k\to 0$ and
$u(r_k,t_k)=0$ for every $k\in Z^+$. Then $u\equiv 0$}

\vspace{.2cm}

\Proof{ By contradiction, suppose that $u\not\equiv 0$. Let $s$ be
the order of $u$ on $S_0$. First, consider the case where $s\in\R$
is a spectral value (say of multiplicity 2). Let $r^sf^+(t)$ and
$r^sf^-(t)$ be the corresponding basic solutions. The function $u$
has the form
\[
 u(r,t)=r^s(a^-f^-(t)+a^+f^+(t))+o(r^s).
\]
 The functions $f^+$ and $f^-$ are independent solutions of the
 first order
differential equation (2.3) and $a^\pm\in\R$ (not both zero). We can
assume $t_k\to t_0$ as $k\to\infty$. It follows from the above
representation of $u$ and from the hypothesis $u(r_k,t_k)=0$ that
$\lim_{k\to\infty}(u(r_k,t_k)/r_k^s)=0$. Consequently,
\[
 a^-f(t_0)+a^+f^+(t_0)=0.
\]
Thus, the solution $a^-f^-+a^+f^+$ of (2.3) is identically zero (by
uniqueness). Hence, $a^-=a^+=0$ which is a contradiction.

If $s$ is not a spectral value, then (by condition (6.14)), it must
be the real part of a unique spectral value
$\sigma\in\C\backslash\R$. The corresponding $\R$-independent basic
solutions are
\[
r^{s+i\beta}\phi(t)+\ov{r^{s+i\beta}\psi(t)}\quad\mathrm{and}\quad
 i(r^{s+i\beta}\phi(t)-\ov{r^{s+i\beta}\psi(t)})
\]
with $\beta\in\R^\ast$ and $|\phi(t)|>|\psi(t)|$ for every $t\in\R$.
The Laurent series of $u$ starts as a linear combination of these
two basic solution and $u$ can then be written as
\[
 u(r,t)=r^s\left( (a^++ia^-)r^{i\beta}\phi(t)+
 (a^+-ia^-)\ov{r^{i\beta}\psi(t)}\right)
 +r^\tau\Phi(r,t)
\]
with $a^\pm\in\R$ (not both zero), $\tau >s$ and $\Phi$  a bounded
function. It follows from the assumption $u(r_k,t_k)=0$ that for
every $k\in\Z$, we have
\[
(a^++ia^-)\phi(t_k)+
 (a^+-ia^-)\ov{r_k^{2i\beta}\psi(t_k)}
 +r_k^{\tau-s-i\beta}\Phi(r_k,t_k)=0,
\]
But this is only possible when $a^++ia^-=0$ (since $|\phi|>|\psi|$,
$r_k\to 0$ and $\Phi$ bounded) }

The next theorem deals with sequences of solutions that converge on
the distinguished boundary $\pa_0U =\pa U\backslash S_0$.

 {\theo Let $U$ be an open and bounded subset
of $\R^+\times\cir$ whose boundary consists of finitely many simple,
closed and piecewise smooth curves. Let $u_n(r,t)$ be a sequence of
bounded functions with $u_n\in C^0(\ov{U}\backslash S_0)$ such that
$\cl u_n=0$ for every $n$. If $u_n$ converges uniformly on
$\pa_0U=\pa U\backslash S_0$, then $u_n$ converges uniformly on
$\ov{U}\backslash S_0$ to a solution $u$ of $\cl u=0$.}

 \vspace{.2cm}

\Proof{ For $(r,t)\in\pa_0U$, let
$\Phi(r,t)=\lim_{n\to\infty}u_n(r,t)$. The function $u(r,t)$ defined
in $U$ by
\[
 u(r,t)=\frac{-1}{2\pi}\int_{\pa_0U}\!\!\!
 \Om_1(r,t,\zeta)\Phi(\zeta)\frac{d\zeta}{\zeta}+
 \ov{\Om_2(r,t,\zeta)}\,\ov{\Phi(\zeta)}\,
 \frac{d\ov{\zeta}}{\ov{\zeta}}
\]
solves $\cl u=0$ (since for each fixed $\zeta$, $\cl\Om_1(z,\zeta)=
\cl\,\ov{\Om_2(z,\zeta)}=0$).
 Now, the Cauchy integral formula applied to $u_n$, shows that
$u$ is the uniform limit of $u_n$ inside $U$}

The following Liouville property is a direct consequence of the
Laurent series expansion and estimate (6.7) of the coefficients.

{\theo Let $u$ be a bounded solution of $\cl u=0$ in $\R^+\times
\cir$. Then
\[
u(r,t)=\sum_{\mathrm{Re}(\sigpmj)=0}a_j^\pm w_j^\pm(r,t)\, .
\]
In particular, if $\cl$ has no spectral values on $i\R$, then
$u\equiv 0$. }

 \vspace{.2cm}

Another consequence of the Laurent series representation is to patch
together solutions from both sides of the characteristic circle
$S_0$. More precisely we have the following theorem.

{\theo Suppose that $\cl $ has no spectral values in $i\R^\ast$.
Then we have the following.
\begin{enumerate}
\item If 0 is not a spectral value, then any bounded solution
of $\cl u=0$ in the cylinder $(-R,\ R)\times\cir$ is continuous
 on the circle $S_0$.
 \item If 0 is a spectral value (say with multiplicity $2$),
 let $g^\pm(t)$ be the basic solutions of $\clstar$ with
 exponent 0. Then a bounded solution $u$ of $\cl u=0$ in
 $\left((-R,\ 0)\cup (0,\ R)\right)\times\cir$ is continuous on
  $(-R,\ R)\times\cir$ if and only if
  \[
 \mathrm{Re}\int_0^{2\pi}g^\pm(\ta)u(\delta ,\ta)d\ta =
 \mathrm{Re}\int_0^{2\pi}g^\pm(\ta)u(-\delta ,\ta)d\ta
  \]
  for some $\delta\in (0,\ R)$.
\end{enumerate}}

\section{The nonhomogeneous equation $\cl u=F$}
 After we extend the Cauchy integral formula to include the
 nonhomogeneous case, we define an integral operator for the
 nonhomogeneous equation $\cl u=F$.
 Throughout the remainder of this paper $U$ will denote an open and
 bounded set in $\R^+\times\cir$ whose boundary consists of finitely
 many simple,
closed and piecewise smooth curves.

\subsection{Generalized Cauchy Integral Formula}
 The following generalization of the Cauchy integral formula will be
 used later.

{\theo Suppose that $F(r,t)$ is a function in $U$ such that
$\dis\frac{F}{r}\in L^p(U)$ with $p\ge 1$. If equation $\cl u=F$ has
a solution $u\in C^0(\ov{U})$, then
 \begeq\begar{ll}
 u(r,t) & =\dis\frac{-1}{2\pi}\int_{\pa_0U}
 \!\!\!\Om_1(r,t,\zeta)u(\zeta)\frac{d\zeta}{\zeta}+
 \ov{\Om_2(r,t,\zeta)}\,
 \ov{u(\zeta)}\,\frac{d\ov{\zeta}}{\ov{\zeta}}\, -\\
 & \dis\qquad\qquad
 \frac{1}{2\pi}\iint_{U}\!\!\!\left[
 \Om_1(r,t,\zeta)F(\zeta)+ \ov{\Om_2(r,t,\zeta)}\,\ov{F(\zeta)}
 \right]\frac{d\rho d\ta}{\rho}
 \stopar\stopeq}

\vspace{.2cm}

 \Proof{ For $\delta >0$, let $U_\delta =U\backslash A(0,\delta)$.
 Let $z_0\in U$ and choose $\delta >0$ so that $z_0\in U_\delta$.
 Green's identity (1.8)
 and  arguments similar to those
 used in the proof of the Cauchy integral formula show that
 \[\begar{ll}
 -2\pi u(z_0)=&\dis
 \int_{\pa U_\delta}
 \!\!\!\Om_1(r,t,\zeta)u(\zeta)\frac{d\zeta}{\zeta}+
 \ov{\Om_2(r,t,\zeta)}\,
 \ov{u(\zeta)}\,\frac{d\ov{\zeta}}{\ov{\zeta}}\, +\\
 & \dis\qquad\qquad
 \iint_{U_\delta}\!\!\!\left[
 \Om_1(r,t,\zeta)F(\zeta)+ \ov{\Om_2(r,t,\zeta)}\,\ov{F(\zeta)}
 \right]\frac{d\rho d\ta}{\rho}
 \stopar\]
 Since $(F/r)\in L^p$ with $p\ge 1$, then  the limits
 of the above integrals as $\delta\to 0$ give
 (7.1) }

 For the adjoint operator, we have the following.

{\theo Let $v(\rho,\ta)\in C^0(\ov{U})$ be such that
 $\dis\frac{\clstar v}{\rho}\in
L^p(U)$ with $p\ge 1$. Then
 \begeq\begar{ll}
 v(\rho ,\ta) & =\dis\frac{-1}{2\pi}\int_{\pa_0U}
 \!\!\!\Om_1(z,\rho,\ta)v(z)\frac{dz}{z}+
 \ov{\Om_2(z,\rho,\ta)}\,
 \ov{v(z)}\,\frac{d\ov{z}}{\ov{z}}\, +\\
 & \dis\qquad
 \frac{1}{2\pi}\iint_{U}\!\!\!\left[
 \Om_1(z,\rho,\ta)\clstar v(z)+ \ov{\Om_2(z,\rho,\ta)}\,
 \ov{\clstar v(z)}
 \right]\frac{dr dt}{r}
 \stopar\stopeq}

 \vspace{.2cm}

\Proof{ Notice that the kernels $\Om_1(z,\zeta)$ and
$\Om_2(z,\zeta)$ satisfy
\[
 L\Om_1=-i\lam\nu\Om_1+c(t)\Om_2\quad\mathrm{and}\quad
 L\ov{\Om_2}=-i\lam\nu\ov{\Om_2}+c(t)\ov{\Om_1}
\]
where $L=\dis\dd{}{ t}-ir\dd{}{r}$. Arguments similar to those used
in the proofs of Theorems 6.4 and 7.1 lead to (7.2). The functions
$P_1$ and $P_2$ used in the proof of Theorem 6.4 need now to be
replaced by the functions $Q_1=\Om_1+\ov{\Om_2}$ and
$Q_2=-i\Om_1+i\,\ov{\Om_2}$}

\subsection{The integral operator $T$}
 We define the operator $T$ and the appropriate $L^p$-spaces in
 which it acts to produce H\"{o}lder continuous solutions.
 For an open set $U\subset\R^+\times\cir$ as before and such that
$S_0\subset\pa U$, we denote by $L^p_a(U)$ the Banach space of
functions $F(r,t)$ such that $\dis\frac{F(r,t)}{r^a}$ is integrable
in $U$ with the norm
\[
||F||_{p,a}=\left(\iint_U
\left|\frac{F(r,t)}{r^a}\right|^pr^{2a-1}drdt\right)^{\frac{1}{p}}\,
.
\]
Note that if $\Phi :\, \R^+\times\cir\,\longrightarrow\,\C^\ast$ is
the diffeomorphism induced by the first integral $z$. That is,
$\Phi(r,t)=r^\lam\ei{it}$, then $F\in L^p_a(U)$ means that the push
forward $\widetilde{F}=F\circ\Phi^{-1}$ satisfies
$\dis\frac{\widetilde{F}(z)}{z}\in L^p(\Phi(U))$.

 We define the integral operator $T$ by
 \begeq
 TF(r,t)=\frac{-1}{2\pi}\iint_U\left[\Om_1(r,t,\zeta)F(\zeta)+
 \ov{\Om_2(r,t,\zeta)}\,\ov{F(\zeta)} \right]
 \frac{d\rho d\ta}{\rho}
 \stopeq
 When $\cl$ has no purely imaginary
spectral values, i.e. $\specl \cap i\R^\ast =\emptyset$, we have the
following theorem.

{\theo Assume $\specl \cap i\R^\ast =\emptyset$. Let $U\subset
A(0,R)$ be an open set as above. The function $TF$ defined by
$(7.3)$ satisfies the followings.
\begin{enumerate}
 \item[$1.$] There exist positive constants $C$ and $\delta$,
  independent on $U$ and $R$,
  such that for every $(r,t)\in U$
  \begeq
  |TF(r,t)|\le CR^\delta ||F||_{p,a}
  \stopeq
  for every $F\in
L^p_a(U)$ with $p>2/(1-\nu)$;
 \item[$2.$] the function $TF$ satisfies the equation $\cl TF=F$; and
 \item[$3.$] the function $TF$ is H\"{o}lder continuous on $\ov{U}$;
\end{enumerate}
Furthermore, if $0$ is not a spectral value, then $TF(0,t)\equiv 0$}

 \vspace{.2cm}

\Proof{ We use the estimates on $\Om_1$ and $\Om_2$ of Theorem 5.1
to write
 \[
-2\pi u(r,t) =I_1+I_2+I_3+I_4\, ,
 \]
 where
 \[\begar{ll}
 I_1 & =i\dis\iint_U\left(\frac{r}{\rho}\right)^{\lam\nu}
 \frac{\zeta}{\zeta -z}F(\zeta)\frac{d\rho
 d\ta}{\rho}\\
 I_2 & =-\dis\iint_U\left(\frac{r}{\rho}\right)^{\lam\nu}
  K(t,\ta)L(z,\zeta)F(\zeta)\frac{d\rho
 d\ta}{\rho}\\
 I_3 & =\dis\frac{1}{2a}\iint_U
 \left[
 \left(\frac{r}{\rho}\right)^{\ov{\lam}\nu}
  c(t)\ov{L(z,\zeta)}+c(\ta)\left(\frac{r}{\rho}\right)^{\lam\nu}
 L(z,\zeta)\right]\ov{F(\zeta)}
  \frac{d\rho
 d\ta}{\rho}\\
 I_4 & =\dis\iint_U\left[ C_1(z,\zeta)F(\zeta)
+\ov{C_2(z,\zeta)}\,\ov{F(\zeta)}\right]
 \frac{d\rho
 d\ta}{\rho}
 \stopar\]
 We use the substitution $\zeta
 =\Phi(\rho,\ta)=\rho^\alpha\ei{i\ta}$ to estimate $I_1$. We find
 \[
 |I_1|\,\le\, r^{a\nu}\iint_{\Phi(U)}
 \frac{|\widetilde{F}(\zeta)|}{|\zeta -z||\zeta|^{1+\nu}}d\xi d\eta\,
 ,\]
  where we have set $\widetilde{F}=F\circ\Phi^{-1}$ and
  $\zeta =\xi +i\eta$. Since
  $\dis\frac{\widetilde{F}}{\zeta}\in L^p(\Phi(U))$
  with $p>2/(1-\nu)$, and since $\Phi(U)$
   is contained in the disc $D(0,R^a)\subset\C$,
    then H\"{o}lder inequality can be used to show that
    there are constants $C$ and $\delta$ so
  that $|I_1|\le r^{a\nu}CR^\delta ||\widetilde{F}||_p$. Furthermore
  these constants are independent on $\Phi(U)$ and $R$.
  Because of the
  logarithmic type growth of $L(z,\zeta)$ and the boundedness of the
  functions $C_1$ and $C_2$, analogous arguments can be used to show
  that $|I_k|\le CR^\delta ||F||_{p,a}$ for $k=2,3,4$.

 Now we verify that $u=TF$ solves $\cl u=F$ in the sense of
 distributions. Let $\psi\in C_0^1(U)$ be a test function.
 The generalized
 Cauchy integral formula (7.2) applied to $\psi$ gives
 \begeq
  \psi(\rho ,\ta) =
 \frac{1}{2\pi}\iint_{U}\!\left[
 \Om_1(z,\rho,\ta)\clstar \psi(z)+ \ov{\Om_2(z,\rho,\ta)}\,
 \ov{\clstar \psi(z)}
 \right]\frac{dr dt}{r}
 \stopeq
 The definition (7.3) of the operator $T$ and estimate (7.4) give
 \[\begar{ll}
 2<TF,\clstar \psi > & =\dis
 \iint_{U}\!\left[TF(z)\clstar\psi(z)+\ov{TF(z)}\,\ov{\clstar\psi(z)}
 \right]\frac{drdt}{r}\\
 & =\dis
 \iint_{U}\!\left[F(\zeta)\psi(\zeta)+\ov{F(\zeta)}\,
 \ov{\psi(\zeta)} \right]\frac{d\rho d\ta}{\rho}\\
 & =2<F,\psi>
 \stopar\]
 This shows that $\cl TF=F$.

 Next, we prove that $TF$ is H\"{o}lder continuous. Since the equation
 is elliptic away from the circle $S_0$,  it is enough to
 prove the regularity of $TF$ on $S_0$. For this, we consider
 the case when 0 is not a spectral value. Then $\, i\R\cap\specl
 =\emptyset\,$ and $\Om_1(0,t,\zeta)=\Om_2(0,t,\zeta)=0$. Hence
 $TF(0,t)=0$. Since $TF$ satisfies $\cl TF=F$, then its pushforward
 $V(z)=TF\circ\Phi^{-1} (z)$ via the first integral satisfies the
 generalized CR equation
 \[
 V_{\ov{z}}=\frac{i\lam\nu}{2ia\ov{z}}V-
 \frac{\widetilde{c}(z)}{2ia\ov{z}}\ov{V}-
 \frac{\widetilde{F}(z)}{2ia\ov{z}}
 \]
where $\widetilde{c}$ and $\widetilde{F}$ are the pushforwards of
$c$ and $F$. We will use the classical results on the CR equation
(see {\cite{Vek}}) to show that $V$ is H\"{o}lder continuous. We
rewrite the above equation as
 \begeq
 V_{\ov{z}}=\frac{G(z)}{z}
 \stopeq
 where
 \[
 G(z)=\frac{i\lam\nu z}{2ia\ov{z}}\,V(z)-
 \frac{z\widetilde{c}(z)}{2ia\ov{z}}\,\ov{V(z)}-
 \frac{z\widetilde{F}(z)}{2ia\ov{z}}\, .
 \]
 Note that since $\widetilde{c}$ and $V$ are bounded functions
 and since $(\widetilde{F}/z)\in L^p$, then
 $G\in L^p(\Phi(U))$ with $p>2$. The solution of (7.6) can then be
 written as $\dis V(z)=\frac{W(z)}{z}$ where $W$ is the solution of
 the equation $W_{\ov{z}}=G$. We know that $W$ is H\"{o}lder
 continuous and has the form
 \[
 W(z)=H(z)-\frac{1}{\pi}\iint_{\Phi(U)}\frac{G(\zeta)}{\zeta -z}
 d\xi d\eta
 \]
 where $H$ is a holomorphic function in $\Phi(U)$. Since
 $V(z)=W(z)/z$ satisfies $V(0)=0$, then necessarily $W(0)=0$ and it
 vanishes to an order $>1$ at 0. Thus $|V(z)|\le K|z|^\tau$
 for some positive constants $K$ and $\tau$. This means
 that $TF$ is H\"{o}lder continuous on $S_0$.

 Finally we consider the case when 0 is a spectral value of $\cl$
 (say, with multiplicity 2).
 Let $f_{j_0}^\pm(t)$ and $g_{-j_0}^\pm(\ta)$ be the basic solutions
 of $\cl$ and $\clstar$ with exponents 0. We have then
 \[
 \Om_1(0,t,\zeta)=\frac{1}{2}f_{j_0}^\pm(t)g_{-j_0}^\pm(\ta)
 \quad\mathrm{and}\quad
 \Om_2(0,t,\zeta)=\frac{1}{2}\ov{f_{j_0}^\pm(t)}\,g_{-j_0}^\pm(\ta).
 \]
 The value of $TF$ on $S_0$ is found to be
 \[
 TF(0,t)=A^+f_{j_0}^+(t)+A^-f_{j_0}^-(t)\, ,
 \]
where
 \[
 A^\pm =\frac{-1}{2\pi}\mathrm{Re}\iint_{U}g_{-j_0}^\pm(\ta)
 F(\zeta)\frac{d\rho d\ta}{\rho}\, .
 \]
 Hence $TF(0,t)$ solves the homogeneous equation $\cl u=0$.
 Let $v(r,t)=TF(r,t)-TF(0,t)$. The function $v$ satisfies
 $\cl v=F$ and $v(0,t)=0$. The push forward arguments, used in the case
 when 0 is not a spectral value, can be used again for the function
 $v$ to establish that $|v(r,t)|\le Cr^\tau$ with $\tau$ and $C$
 positive}

 In general, when $\cl$ has spectral values on $i\R$, we can define
 $\widehat{\Om}_1$ and $\widehat{\Om}_2$ as in (5.2) and (5.3) except that
 the terms corresponding to $\sigpmj$ that are in $i\R$ are missing from
 the sums. That is, if $w^\pm_1,\,\cdots , \, w_p^\pm$ denotes the
 collection of basic solutions of $\cl$ with exponents in $i\R$,
 and $w^{\ast\pm}_1,\,\cdots\,, w^{\ast\pm}_p$ the corresponding
 collection of basic solutions of the adjoint $\cl^\ast$, then
 \begeq
 \widehat{\Om}_k(z,\zeta)=\Om_k(z,\zeta)-\sum_{k=1}^p
 w_k^\pm(z)w^{\ast \pm}_k(\zeta)\, .
 \stopeq
 We define the modified operator $\widehat{T}$ by
 \begeq
 \widehat{T}F(r,t)=\frac{-1}{2\pi}\iint_U\left[\widehat{\Om}_1(r,t,\zeta)F(\zeta)+
 \ov{\widehat{\Om}_2(r,t,\zeta)}\,\ov{F(\zeta)} \right]
 \frac{d\rho d\ta}{\rho}\, .
 \stopeq
 Arguments similar to those used in the proof of Theorem 7.3
 establish the following result.

 {\theo Let $U\subset A(0,R)$ be as above. Then the function $\widehat{T}F$
 defined by $(7.8)$ satisfies properties $1$, $2$, and $3$ of Theorem $7.3$
 and $\widehat{T}F(0,t)=0$.}

  \subsection{Compactness of the operator $T$}

 {\theo Suppose that $\cl$ has no spectral values in
 $i\R^\ast$. Then, for $p>2/(1-\nu)$, the operator
 $T\, :\, L_a^p(U)\,\longrightarrow\,
 C^0(\ov{U})$ is compact.}

 \vspace{.2cm}

 \Proof{ Let $R>0$ be such that $U\subset A(0,R)$.
  A function in $L_a^p(U)$
 can be considered in $L^p_a(A(0,R))$ by extending as 0 on
 $A(0,R)\backslash U$. Denote by $T_R$ the operator $T$ on the cylinder
 $A(0,R)$ and set
 \[
 \widehat{T}_RF(r,t)=T_RF(r,t)-T_RF(0,t)
 \]
 Thus
 \[
 \widehat{T}_RF(r,t)=\frac{-1}{2\pi}\iint_{A(0,R)}\!
 \left[\widehat{\Om}_1(z,\zeta)F(\zeta)+
 \ov{\widehat{\Om}_2(z,\zeta)}\,\ov{F(\zeta)}\right]
 \frac{d\rho d\ta}{\rho}
 \]
 where $\widehat{\Om}_1$ and $\widehat{\Om}_2$ are defined by
 (5.2) and (5.3), respectively, except that the terms corresponding to the
 spectral value $\sig_j=0$ are missing. In particular, if 0 is not
 a spectral value, then $\widehat{\Om}_1=\Om_1$ and
 $\widehat{\Om}_2=\Om_2$. Note that $\widehat{\Om}_1(0,t)=0$ and
 $\widehat{\Om}_2(0,t)=0$. The operator $\widehat{T}_RF$ satisfies
 the properties of Theorem 7.3 and
 $\widehat{T}_RF(0,t)=0$. To show that $T$ is compact, it is  enough
 to show that $\widehat{T}_R$ is compact.

  Let $B\subset L^p_a(U)$ be a bounded set. We need to show that
  $\widehat{T}_R(B)$ is relatively compact in $C^0(\ov{U})$.
  Let $M>0$ be such that $||F||_{p,a}\le M$ for every $F\in B$. It follows
  from Theorem 7.3 that $\widehat{T}_R(B)$ is bounded (by $CR^\delta
  M$).
  Now we show the equicontinuity of $\widehat{T}_R(B)$. First along $S_0$.
  For $\ep>0$, let $r_0>0$ be
  such that $Cr_0^\delta M<(\ep/2)$. We have then
  \[
 \widehat{T}_RF(z)=\widehat{T}_{r_0}F(z)+\widehat{T}_{A(r_0,R)}F(z)
  \]
  where $\widehat{T}_{A(r_0,R)}$ denotes the integral operator over
  the cylinder $A(r_0,R)$. Let $r_0$ be small enough so that
  \[
 E=\max_{r<(r_0/2),\ r_0<\rho <R}
 (|\widehat{\Om}_1(r,t,\rho,\ta)|+
 |\widehat{\Om}_2(r,t,\rho,\ta)|)\, <\,
 \frac{\ep}{2M(\pi R^{2a})^{1/q}}
  \]
  where $q$ is such that $\dis\frac{1}{p}+\frac{1}{q}=1$.
  For $r<r_0/2$, we have then
  \[
 |\widehat{T}_{r_0}F(r,t)|\, \le\, C r_0^\delta M\, \le\, \frac{\ep}{2}
  \]
  and
  \[
 |\widehat{T}_{A(r_0,R)}F(r,t)|\le
 \iint_{A(r_0,R)}\!\!\!(|\widehat{\Om}_1|+
 |\widehat{\Om}_2|) \frac{|F(\zeta)|}{\rho}\,d\rho d\ta
 \le E||F||_{p,a}\le \frac{\ep}{2}
  \]
  This estimate is obtained from H\"{o}lder's inequality
  and the above estimate on $E$. Thus, $ |\widehat{T}_RF(r,t)|\le
  \ep$ and so $\widehat{T}_RB$ is equicontinuous on $S_0$.

  Next, let $(r_1,t_1)\in U$ with $r_1<r_0/4$. Set
  $z_1=r_1^\lam \ei{it_1}$ and $z=r^\lam \ei{it}$. If
  $|z-z_1|<r_0^a/4$, then $r <r_0/2$ and the above argument gives
  \[
 |\widehat{T}_RF(z)-\widehat{T}_F(z_1)|\le
 |\widehat{T}_R(z)|+|\widehat{T}_R(z_1)|\le 2\ep,\qquad
 \forall F\in B
  \]
  Finally, suppose that $r_1>r_0/4$. Let $b$ be such that
  $0<b< r_0/4$. We write
  \[
 \widehat{T}_RF(z)-\widehat{T}_RF(z_1)=
\widehat{T}_bF(z)-\widehat{T}_bF(z_1)+
\widehat{T}_{A(b,R)}F(z)-\widehat{T}_{A(b,R)}F(z_1)
  \]
  After using H\"{o}lder's inequality we obtain
  \[
 |\widehat{T}_bF(z)-\widehat{T}_bF(z_1)|\,\le\, CS(b)||F||_{p,a}
  \]
  where
  \[
 S(b)=\max_{P}
 (|\widehat{\Om}_1(z,\zeta)-\widehat{\Om}_1(z_1,\zeta)|+
 |\widehat{\Om}_2(z,\zeta)-\widehat{\Om}_2(z_1,\zeta)|)\, ,
  \]
 and where the maximum is taken over the set $P$ of points satisfying
 $\rho <b$, $|r-r_1|<b$, $r>r_0/4$ and $r_1>r_0/4$.
 The continuity of the kernels in the region $\rho <b$ and $r>r_0/4$
 implies that if $b$ is small enough, then $S(b)<\ep/(2MC)$ and
 consequently
 \[
 |\widehat{T}_bF(z)-\widehat{T}_bF(z_1)|\,\le\,\frac{\ep}{2}
 \]
 Finally, for $\widehat{T}_{A(b,R)}F$, it suffices to notice that it
 solves the equation $\cl u=F$ in the cylinder $A(r_0/4,R)$.
 In this cylinder, the equation is elliptic and the classical theory
 of generalized analytic function ({\cite{Vek}} Chapter 7) implies that
 the family $\widehat{T}_{A(b,R)}B$ is equicontinuous.}

 \section{The semilinear equation}
 In this section, we make use of the operator $T$ and of its
 modified version, through the kernels $\Om_{j,1}$ and $\Om_{j,2}$
 (defined in section 5), to establish a correspondence between the
 solutions of the homogeneous equation $\cl u=0$ and the solutions
 of a semilinear equation.

 {\theo Assume that $\cl$ has no spectral values in $i\R^\ast$.
 Let $G(u,r,t)$ be a bounded function defined in
 $\C\times A(0,R_0)$, for some $R_0>0$,
 and let $\tau >a\nu$. Then, there are $R>0$ and a one to one map
 between the space of continuous solutions of the equation
 $\cl u=0$  in $A(0,R)$ and the space of continuous solutions of the
 equation
 \begeq
 \cl u =r^\tau |u|G(u,r,t)\, .
 \stopeq
 Furthermore, if $v$ is a bounded solution of $(8.1)$ in a cylinder
 $A(0,R)$, then $v$ is continuous up to the circle $S_0$.}

 \vspace{.2cm}

 \Proof{ First note that since $\tau >a\nu$, then $r^\tau\in
 L^p_a(A(0,R))$, with a $p$ satisfying $p>2/(1-\nu)$, and
 $||r^\tau||_{p,a}=C_1R^{\delta_1}$ with $C_1$ and $\delta_1$ positive.
 Consider the operator
 \[
 \mathcal{P}:\,
 C^0(\ov{A(0,R)})\,\longrightarrow\,C^0(\ov{A(0,R)})\, ;\quad
 \mathcal{P}(f)=T_R(r^\tau |f|G(f,r,t))
 \]
 where, as before, $T_R$ denotes the integral operator on the cylinder
 $A(0,R)$. Note that since $G$ is a bounded function, then $\mathcal{P}$ is
 well defined. It follows from the properties of $T$ given
 in Theorem 7.3, from the boundedness of $G$,  and from $r^\tau \in L^p_a$
 that the operator $\mathcal{P}$ satisfies
 \[
 \cl(\mathcal{P}(f))=r^\tau |f|G(f,r,t)
 \]
 and
 \[
 |\mathcal{P}(f)(r,t)|\le CR^\delta ||r^\tau |f|G(f,r,t)||_{p,a}
 \le \,
 C'R^{\delta'} ||f||_0\qquad \forall f\in C^0(\ov{A(0,R)})
 \]
 with $C'$ and $\delta'$ positive. Hence, if $R>0$ is small enough,
 $||\mathcal{P}||\le C'R^{\delta'} <1$,  and $\mathcal{P}$ is thus a
 contraction. Let $\mathcal{F}=(I-\mathcal{P})^{-1}$. The operator
 $\mathcal{F}$ realizes the one to one mapping between the space
 of continuous solutions of $\cl u=0$ and those of equation (8.1).

 Now, we show that if $v$ a bounded solution of (8.1) in a cylinder
 $A(0,R)$, then it is continuous.
 For a bounded solution $v$, the function $r^\tau |v|G(v,r,t)$ is bounded and
 is in $L^p_a(A(0,R))$. Consequently, $\mathcal{P}(v)$ is continuous up
 to the boundary $S_0$.
 The function $u=v-\mathcal{P}(v)$
 is a bounded solution of $\cl u=0$ and so it is
 continuous up to $S_0$ (Remark 6.1). It
  follows that
 $v=u+\mathcal{P}(v)$ is also continuous up to $S_0$ }

 Let $\sig_j =\sigpj$ (or $\sig_j=\sigmj$) be a spectral value
 of $\cl$ such that $\mathrm{Re}(\sig_j)>0$.
 Consider the
 Banach spaces $r^{\sig_j}L^p_a(A(0,R))$ and
 $r^{\sig_j}C_b^0(A(0,R))$ defined as follows:
 $f\in r^{\sig_j}L^p_a(A(0,R))$ if
 $\dis\frac{f}{r^{\sig_j}}\in L^p_a$ and $g\in r^{\sig_j}C_b^0(A(0,R))$
 if $\dis\frac{g}{r^{\sig_j}}\in C^0(A(0,R))$ and is
 bounded. The norms in these spaces are defined by
 \[
 ||f||_{p,a,\sig_j}=||\frac{f}{r^{\sig_j}}||_{p,a}\quad
 \mathrm{and}\quad
 ||g||_{0,\sig_j}=||\frac{g}{r^{\sig_j}}||_0
 \]
 Consider the operator $T_R^j$ defined by
 \begeq
 T_R^jF(r,t)=\frac{-1}{2\pi}\iint_{A(0,R)}\!
 \left[\Om_{j,1}(z,\zeta)F(\zeta)+
 \ov{\Om_{j,2}(z,\zeta)}\, \ov{F(\zeta)}\right]\frac{d\rho d\ta}{\rho}
 \stopeq
where $\Om_{j,1}$ and $\Om_{j,2}$ denote the modified kernels
defined in (5.11) and (5.12). Note that the estimates of Theorem 5.2
on the modified kernels imply that $T_R^jF$ is in
$r^{\sig_j}C_b^0(A(0,R))$ when $F$ is in $r^{\sig_j}L^p_a(A(0,R))$.
 Arguments similar to those used in the proof of Theorem 7.3 can be
used to establish the following theorem.

{\theo For $p>2$, the operator
 \[
 T_R^j:\, r^{\sig_j}L^p_a(A(0,R))\, \longrightarrow\,
 r^{\sig_j}C_b^0(A(0,R))
 \]
 satisfies
 \begeq
 \cl T_R^jF(r,t)=F(r,t) \quad\mathrm{and}\quad
  ||T_R^jF||_{0,\sig_j}\le CR^\delta ||F||_{p,a,\sig_j}
 \stopeq
 where $C$ and $\delta$ are positive constants.}

\vspace{.2cm}

 Two functions $u$ and $v$ defined in the cylinder $A(0,R)$ are
 said to be similar if $u/v$ is continuous in $A(0,R)$
 and there exist positive constants $C_1$ and $C_2$ such that
 \[
 C_1\le \left|\frac{u(r,t)}{v(r,t)}\right|\le C_2,\qquad
 \forall (r,t)\in A(0,R)\, .
 \]

 {\theo Let $\cl$, $\tau$,  and $G$ be as in Theorem $8.1$.
 Then there exists
 $R>0$ such that each continuous solution of $\cl u=0$ in
 $\ov{A(0,R)}$ is
 similar to a continuous solution
 of equation $(8.1)$.}

 \vspace{.2cm}

\Proof{ Let $u$ be a continuous solution of $\cl u=0$ on
$\ov{A(0,R)}$. Let $\mu\ge 0$ be the order of $u$ along $S_0$. If
$\mu
>0$, then $\mu =\mathrm{Re}(\sig_j^\pm)$ for some spectral value
$\sig_j^\pm$. Assume that $\sig_j^-=\sig_j^+=\sig_j$. Then it
follows from the Laurent series expansion that $u$ is similar to a
linear combination
 \[
 u_0(r,t)=a^-w_j^-(r,t)+a^+w_j^+(r,t)
 \]
 of the basic solutions
$w_j^+$ and $w_j^-$ in a cylinder $A(0,R_1)$ with $0<R_1<R$.
Consider the operator
  \[
 \mathcal{P}_j:\,
 r^{\sig_j}C_b^0(A(0,R))\,\longrightarrow\, r^{\sig_j}C_b^0(A(0,R))
 \]
 defined by
 $\mathcal{P}_j(f)=T^j_R(r^\tau |f|G(f,r,t))$
 It follows from the hypotheses on $\tau$, on $ G$, and from Theorem
 8.2 that
 \[
 \cl\mathcal{P}_j(f)=r^\tau |f|G(f,r,t) \quad\mathrm{and}\quad
 ||\mathcal{P}_j(f)||_{0,\sig_j}\le KR^\delta ||f||_{0,\sig_j}
 \]
 for some positive constants $K$ and $\delta$. In particular,
 $\mathcal{P}_j$ is a contraction, if $R$ is small enough. Hence,
 the function $v=(I-\mathcal{P}_j)^{-1}(u)$ is a solution of
 equation (8.1) and it is also similar to $u_0$.

If $\mu =0$ (then we are necessarily in the case where 0 is a
spectral value), let $\mathcal{F}$ be the resolvent of $\mathcal{P}$
used in the proof of Theorem 8.1.  Then $v=\mathcal{F}(u)$ is
similar to $u$ and solves (8.1)}

A direct consequence of Theorem 8.1 and Theorem 6.6 is the following
uniqueness result for the solutions of (8.1)

 {\theo Suppose that
the spectral values of $\cl$ satisfy the following condition
 \[
\mathrm{Re}(\sig_j)=\mathrm{Re}(\sig_k)\,\Longrightarrow\,
\sig_j=\sig_k, \qquad\forall\sig_j,\sig_k\in\specl\, .
 \]
Let $u$ be a bounded solution of $(8.1)$ in a cylinder $A(0,R)$.
Suppose that  there is a sequence of points $(r_k,t_k)\in A(0,R)$
such that $r_k\to 0$ and $u(r_k,t_k)=0$ for every $k\in Z^+$. Then
$u\equiv 0$}

\section{The second order equation: Reduction}
 This section deals with the second order operator $P=L\ov{L}
 +\mathrm{Re}(aL)$. We show that the equation $Pu=F$, with $u$ and
 $F$ real-valued, can be reduced to an equation of the form
 $\cl u=G$.

 As before, let $\lam=a+ib\in\R^++i\R$, $\, L$  be the vector
field
 given by (1.1) and let $\beta(t)\in C^m(S^1,\C)$, with
 $m\ge 2$, satisfies
 \begeq
 \frac{1}{2\pi i}\int_0^{2\pi}\beta (t)dt =k\, \in\, \Z\, .
 \stopeq
  Consider the second order
 operator $P$ defined as
 \begeq
 P=L\ov{L}+\ov{\lam}\beta(t)L+\lam\ov{\beta(t)}\,\ov{L}\, .
 \stopeq
 Then,
 \[
 Pu=|\lam|^2u_{tt}-2bru_{rt}+r^2u_{rr}+
 |\lam|^2(\beta+\ov{\beta})u_t+
 [1-i(\ov{\lam}\,\beta -\lam\ov{\beta})]ru_r.
 \]
 Note that $P$ is elliptic except along the circle
 $S_0=\{0\}\times\cir$, and that $Pu$ is $\R$-valued when $u$ is
 $\R$-valued.

With the operator $P$ we associate a first order operator $\cl$ and
show that the equation $Pu=F$, with $F$ real-valued, is equivalent
to an equation of the form $\cl w =G$. Let
 \begeq
 B(t)=\exp\int_0^t\ov{\beta(s)}ds\, .
 \stopeq
 It follows from (9.1) that $B$ is periodic with $\ind(B)=-k$ and
 satisfies $LB=\lam\ov{\beta}\, B$. Define the function $c(t)$ by
 \begeq
 c(t)=-\ov{\lam}\, \beta(t)\frac{B(t)}{\ov{B(t)}}=
 -\ov{\lam}\,
 \beta(t)\exp\left[\int_0^t(\ov{\beta(s)}-\beta(s))ds\right]\, .
 \stopeq
 Note that the function $B$ satisfies also the equation
 \begeq
 LB=-\ov{c(t)}\,\frac{B^2(t)}{\ov{B(t)}}
 \stopeq
 To each $\R$-valued function $u$, we associate the
 $\C$-valued function $w$  defined by
 \[
 w(r,t)=B(t)\ov{L} u(r,t)\, .
 \]
 We will refer to $w$ as the $L$-potential of $u$ with respect to  $P$.
 To the
operator $P$ we associate the first order operator $\cl$ defined by
 \begeq
 \cl w =Lw -c(t)\ov{w}\, ,
 \stopeq
where $c(t)$ is given by (9.4). We have the following proposition.

{\prop Suppose that $F$ is an $\R$-valued function and $u(r,t)$ is
$\R$-valued and solves the equation
 \begeq
 Pu(r,t)=F(r,t).
 \stopeq
 in the cylinder $A(0,R)$.
 Then its $L$-potential $w$ satisfies the equation
 \begeq
 \cl w(r,t)=B(t)F(r,t).
 \stopeq
Conversely, if $w$ is a solution of $(9.8)$ in $A(0,R)$, then there
is an $\R$-valued function $u$ defined in $A(0,R)$ that solves
$(9.7)$ and whose $L$-potential is $w$. More precisely, the function
$u$ can be defined by
 \begeq
 u(r,t)=\mathrm{Re}\int_{(r_0,t_0)}^{(r,t)}\frac{w(\rho,\ta)}{B(\ta)}
 \, \frac{d\zeta}{ia\zeta}
 \stopeq
 where $\zeta =\rho^\lam\ei{i\ta}$ and the integration is taken over
 any simple curve in $A(0,R)$ that joins the fixed point
 $(r_0,t_0)$ to the point $(r,t)$ }

 \vspace{.2cm}

\Proof{ Suppose that $u$ is $\R$-valued and solves (9.7). By using
(9.4) and (9.5), we see that its $L$-potential satisfies
 \[\begar{ll}
 Lw=L(B\ov{L}u) & =BL\ov{L}u+LB\ov{L}u =BF-\ov{\lam}\beta
 BLu\\
  & \dis =BF-\ov{\lam}\,\beta\frac{B}{\ov{B}}\, \ov{(B\ov{L}u)}=
  BF+c\ov{w}
 \stopar\]
 Thus $w$ solves (9.8). Conversely, suppose that $w$ solves
 (9.8). Let $(r_0,t_0)\in A(0,R)$ and consider the function
 $u(r,t)$ defined by (9.9). We need to verify that the integral is
 path independent. Let $U$ be a relatively compact subset of
 $A(0,R)$ whose boundary consists of simple closed curves.
 It follows from the proof of Green's identity (1.8) and
 (9.5) that
 \[\begar{ll}
 \dis\int_{\pa U}\!\!\frac{w(\zeta)}{B(\ta)}\,\frac{d\zeta}{ia\zeta} & =
 \dis\iint_U\!\! L\left(\frac{w}{B}\right)\frac{d\ov{\zeta}
 d\zeta}{2a^2|\zeta|^2} =\iint_{ U}\!\!
 \left[\frac{Lw}{B}-\frac{LB}{B^2}w\right]\frac{id\rho d\ta}{\rho}\\
 & =\dis\iint_U\!\! \left[F+\frac{c}{B}\,\ov{w}+
 \frac{\ov{c}}{\ov{B}}\, w\right]\frac{id\rho d\ta}{\rho}
 \stopar\]
 Since $F$ is $\R$-valued, then the real part of the above integral
 is zero and the function $u$ is well defined. That $u$ satisfies
 (9.7) follows easily by computing the derivatives $u_t$ and $u_r$ from
 (9.9) to obtain $\ov{L} u=w/B$ and then using (9.8) to get (9.7) }

 \section{The homogeneous equation $Pu=0$}
  We use the reduction given in Proposition 9.1 to obtain properties
  of the solutions of the equation $Pu=0$ from those of their
  $L$-potentials $w$. In
  particular, series representation  for $u$ in a cylinder is derived.
  Under an assumption on the spectrum of $\cl$, we prove a maximum
  principle for the equation $Pu=0$: The extreme values of $u$ can
  occur only on the distinguished boundary $\pa_0U$. It should be
  mentioned that many results in this section and the next are close
  to those
  obtained, in {\cite{Mez-CV4}}. There, the operator in $\R^2$
   has its principal part of
   the particular form $(x^2+y^2)\Delta$, where $\Delta$
  is the Laplacian. Such an operator, when written in polar
  coordinates has the form (9.2) with the vector field $L$ having
  the invariant $\lam =1$.

 \subsection{Some properties}
 The following simple properties for the solutions $u$ will be
 needed .
We start by considering the possibility of the existence of radial
solutions.

{\prop The equation $Pu=0$ has radial solutions $u=u(r)$ if and only
if the coefficient $\beta$ has the form
 \begeq
 \beta (t)=\frac{\lam}{a}p(t)-ik
 \stopeq
 where $k\in\Z$ and $p(t)$ is $\R$-valued and such that
 $\dis\int_0^{2\pi}\!\! p(t) dt=0$. In this case, the radial solutions
 have the form
 \begeq
 u(r) =\left\{\begar{ll}
 C_1\log r +C_2 &\quad\mathrm{if}\ \ k=0\\
 \dis C_1 r^{2ak} +C_2 &\quad\mathrm{if}\ \ k\ne 0
 \stopar\right.
 \stopeq
 where $C_1$, $C_2$ are arbitrary constants.
 The corresponding $L$-potentials are $w(r,t)=iC_1B(t)$, when
 $k=0$ and $w(r,t)=2iakC_1r^{2ak}B(t)$, when $k\ne 0$,  where
 \[
 B(t)=\ei{ikt}\exp\left(\frac{\ov{\lam}}{a}\int_0^tp(s)ds\right)
 \]
 Moreover, $w(r,t)$ is a basic solution of $\cl$ with character
 $(2ak,\ k)$.}

 \vspace{.2cm}

 \Proof{ If $u=u(r)$ solves $Pu=0$, then
 \[
 r^2u''(r)+(1-i(\ov{\lam}\,\beta(t)-\lam\,\ov{\beta(t)}))ru'(r)=0\, .
 \]
 Hence, $i(\ov{\lam}\,\beta(t)-\lam\,\ov{\beta(t)})$ is a real constant.
 If we set $\beta(t)=p(t)+iq(t)$ with $p$ and $q$ real-valued, then
 $aq(t)-bp(t)=M$,
 with $M\in\R$ constant. It follows from hypothesis (9.1) that
 $M=-ak$ with $k\in\Z$ and that the average of $p$ is zero. This gives
 $aq(t)=bp(t)-ak$ and consequently $\beta$ has the form (10.1).
 For such a coefficient $\beta$, the radial solutions are easily
 obtained from the differential equation}

 \Rem{10.1}{ Note that if $u=u(t)$ (independent on $r$)
 solves the equation $Pu=0$, then $u$ is necessarily constant.}

 The following lemma will be used in the proof of the next
 proposition.

{\lem Let $u(r,t)$ be a solution of $Pu=0$ in the cylinder $A(0,R)$
and let $w=B\ov{L}u$ be its $L$-potential. If
 \[
\mathrm{Re}\left[ \frac{\lam w(r,t)}{iaB(t)}\right] \equiv 0\, ,
 \]
 then $u$ is constant.}

 \vspace{.2cm}

\Proof{ Let $(r_0,0)$ be a fixed point in the cylinder $A(0,R)$.
 Let $u$ be as in the lemma and let
$w=B\ov{L}u$. The function
 \[
 v(r,t)=\mathrm{Re}\int_{\Gamma(r,t)}\frac{w(\zeta)}{B(\ta)}
 \frac{d\zeta}{i a\zeta}\,
 \]
where $\Gamma (r,t)$ is any piecewise smooth curve
 that joins the point $(r_0,0)$ to the point $(r,t)$,
solves $Pv=0$ (Proposition 9.1). We choose $\Gamma$ as
 $\Gamma =\Gamma_1\cup\Gamma_2$, where
 \[
 \Gamma_1=\{ (r_0,st),\ 0\le s\le 1\}\ \mathrm{and}\
 \Gamma_2=\{ ((1-s)r_0+sr,t),\ 0\le s\le 1\} .
 \]
With this choice of $\Gamma$ and with the hypothesis of the lemma on
the potential $w$, the integral over $\Gamma_2$ is 0 and the
expression for $v$ reduces to
\[
 v(r,t)=\frac{t}{a}\mathrm{Re}\int_0^1\frac{w(r_0,st)}{B(st)}ds\, .
\]
Hence, the function $v$ depends only on the variable $t$ and since
it solves $Pv=0$, then $v$ is constant (Remark 10.1). Consequently,
$w=B\ov{L}v=0$. This means $\ov{L}u=0$. Since $u$ is $\R$-valued,
then $Lu=0$ and so $u$ is constant }

{\prop Suppose that $u\in C^0(\ov{A(0,R)})$ solves $Pu=0$, then its
$L$-potential $w$ satisfies, $w\in C^0(A(0,R)\cup S_0)$ and
$w(0,t)\equiv 0$. Moreover, $u$ is constant along $S_0$.}

\vspace{.2cm}

 \Proof{ Since $P$ is elliptic for $r\ne 0$, then we
need only to verify the continuity of $w$ up to $S_0$ and its
vanishing there.
 As a solution of $\cl w=0$, the function $w$ has a Laurent series
 expansion (Theorem 6.2)
 \[
 w(r,t)=\sum_{j\in\Z}a_j^\pm w_j^\pm(r,t)
 \]
 where $w_j^\pm$ are the basic solutions of $\cl$. Let $\tau\in\R$
 be the order of $w$ along $S_0$ (that the order
 $\tau$ is finite is a
 consequence of the continuity of $u$ up to $S_0$). We are going to
 show that $\tau >0$. Let
 $w_1,\cdots , w_N$ be the collection of all basic solutions with
 order $\tau$ along $S_0$. That is,
 $w_m$ is a basic solution with $\charac(w_m)=(\sig_{j_m},j_m)$
 and such that the exponent satisfies
 $\mathrm{Re}(\sig_{j_m})=\tau$. We have then
 \[
 w(r,t)=\sum_{k=1}^Na_kw_k(r,t)+o(r^\tau)=w_\tau (r,t)+o(r^\tau)
 \]
 It follows from Lemma 10.1 that for each $k$, $\, \mathrm{Re}(\lam w_k
 /iaB)\not\equiv 0$. Let $t_0\in\R$ be such that
 \[
 \mathrm{Re}\left(\frac{\lam w_k(r,t_0)}{ia B(t_0)}\right) \ne
 0,\qquad k=1,\cdots ,N.
 \]
 Let $r_0<R$ be fixed. By using integration over the segment from
 $(r_0,t_0)$ to $(r,t_0)$, we find
 \[
 u(r,t_0)-u(r_0,t_0)=\int_0^1\mathrm{Re}
 \left[\frac{\lam w_\tau ((1-s)r_0+sr,t_0)}{iaB(t_0)} \right]
 \,\frac{(r-r_0)ds}{(1-s)r_0+sr}\, +\, o(r^\tau).
 \]
 Recall that each basic solution $w_1,\cdots ,w_n$ has an
 exponent $\sig_k=\tau +i\beta_k$ and so
 \[ \mathrm{Re}
 \int_0^1
 \frac{(r-r_0)\lam w_k ((1-s)r_0+sr,t_0)}{
 ((1-s)r_0+sr)iaB(t_0)}ds =\left\{
 \begar{ll}
 O(r^\tau) &\ \mathrm{if}\ \tau\ne 0\\
 O(\log r) &\ \mathrm{if}\ \tau = 0, \ \beta_k=0\\
 O(r^{i\beta_k}) &\ \mathrm{if}\ \tau =0,\ \beta_k \ne 0
 \stopar\right.
 \]
  From these estimates and the above integral, we deduce that
  in order for $u(r,t_0)-u(r_0,t_0)$ to have a limit as
  $r\to 0$, it is necessary that $\tau >0$. For such $\tau$,
  $w(0,t)=0$,
  and $u(0,t)$ is constant. }

  As a consequence of the proof of Proposition 10.2, we have the
following proposition.

{\prop Suppose that $\cl$ has no spectral values on $i\R^\ast$.
 If $u\in L^\infty (A(0,R))$ solves $Pu=0$, then $u$ is continuous
 up to the boundary $S_0$ and it is constant on $S_0$.}

\subsection{Main result about the homogeneous equation $Pu=0$}
 We use the basic solutions of the associated operator $\cl$ to
 construct $2\pi$-periodic functions $q_j^\pm(t)$ and establish a
 series expansion of the continuous solutions $u$.

 Let $\{\sigpmj\}_{j\in\Z}$ be the spectrum of
   the associated operator $\cl$ and $w_j^\pm$ be the corresponding
  basic solutions. Recall that if $\sigpmj\in\R$, then
  $w_j^\pm =r^{\sigpmj}f_j^\pm(t)$ with
  $\ind (f_j^\pm)=j$ and if $\sig_j^+\in\C\backslash\R$,
  then $\sig_j^-=\sig_j^+=\sig_j$ and
  \[
 w_j^+(r,t)=r^{\sig_j}\phi_j(t)+\ov{r^{\sig_j}\psi_j(t)}\, ,
 \qquad
 w_j^-(r,t)=i\left[r^{\sig_j}\phi_j(t)-\ov{r^{\sig_j}\psi_j(t)}\right]
  \]
  with $|\phi_j|>|\psi_j|$ and $\ind (\phi_j)=j$.
 Define the functions $q_j^\pm(t)$ as follows.
 For $\sigpmj\in\R^\ast$,
 \begeq
 q_j^\pm(t)=\frac{\lam}{ia\sigpmj}\frac{f_j^\pm(t)}{B(t)}
 \stopeq
 and for $\sig_j\in\C\backslash\R$,
 \begeq
 q_j^+(t)  =\dis\frac{1}{ia\sig_j}\left(\frac{\lam\phi_j(t)}{B(t)}
 -\frac{\ov{\lam}\,\psi_j(t)}{\ov{B(t)}}\right) ,\quad
 q_j^-(t)  =\dis\frac{1}{ a\sig_j}\left(\frac{\lam\phi_j(t)}{B(t)}
 +\frac{\ov{\lam}\,\psi_j(t)}{\ov{B(t)}}\right) .
 \stopeq
It follows from Theorem 4.1 that the asymptotic behaviors of
$q_j^\pm$ are
\[
 q_j^+(t)=\frac{\ei{ijt}}{iajB(t)}+O(j^{-2})
 \quad\mathrm{and}\quad
 q_j^-(t)=\frac{\ei{ijt}}{ajB(t)}+O(j^{-2}).
\]
We have the following representation theorem

{\theo If $u\in C^0(\ov{A(0,R)})$ is a solution of $Pu=0$, then
 $u$ is constant on $S_0$ and it has the series
 expansion
 \begeq
 u(r,t)=u_0+\!\!\sum_{\mathrm{Re}(\sigpmj) >0}
 u_j^\pm \mathrm{Re}\left[ r^{\sigpmj} q_j^\pm (t)\right]
 \stopeq
 where the functions $q_j^\pm$ are defined in $(10.3)$ and $(10.4)$,
 and where $u_j^\pm\in\R$. }

 \vspace{.2cm}

\Proof{ It follows from Proposition 10.2 that $u$ is constant on
$S_0$. Hence, by using integration over the segment from $(0,t)$ to
$(r,t)$, we obtain
 \[
 u(r,t)-u(0,0)=u(r,t)-u(0,t)=
 \mathrm{Re}\int_0^1\frac{\lam w(sr,t)}{iaB(t)}\frac{ds}{s}
 \]
 where $w$ is the $L$-potential of $u$. The function $w$, being a solution
 of $\cl w=0$,  has
 a series expansion
 \[
 w(r,t)=\sum_{\mathrm{Re}(\sigpmj)>0}\!\!\!c_j^\pm w_j^\pm (r,t).
 \]
 For the function $u$ we have then
 \[
 u(r,t)=u(0,0)+\!\!\sum_{\mathrm{Re}(\sigpmj)>0}\!\!\!c_j^\pm
 \mathrm{Re}
 \int_0^1\frac{\lam w_j^\pm (sr,t)}{iaB(t)}\frac{ds}{s}\, .
 \]
 Now for $\sig_j^\pm\in\R$, we have $w_j^\pm
 (r,t)=r^{\sigpmj}f_j^\pm (t)$ and
 \[
 \int_0^1\frac{\lam w_j^\pm (sr,t)}{iaB(t)}\frac{ds}{s}=
 r^{\sigpmj}\frac{\lam}{ia\sigpmj}\frac{f_j^\pm(t)}{B(t)}=
 r^{\sigpmj}q_j^\pm(t)\, .
 \]
 For $\sigpj=\sigmj=\sig_j \in\C\backslash\R$, we have
 \[\begar{ll}
 \dis\int_0^1\frac{\lam w_j^+ (sr,t)}{iaB(t)}\frac{ds}{s} & =
 \dis \int_0^1\left[(rs)^{\sig_j}\frac{\lam\phi_j(t)}{iaB(t)}
 +(rs)^{\ov{\sig_j}}\frac{\lam\ov{\psi_j(t)}}{iaB(t)}\right]
 \frac{ds}{s}\\
 & =\dis
 \frac{r^{\sig_j}}{\sig_j}\frac{\lam\phi_j(t)}{iaB(t)}+
 \frac{r^{\ov{\sig_j}}}{\ov{\sig_j}}
 \frac{\lam\ov{\psi_j(t)}}{iaB(t)}.
 \stopar\]
 From this and (10.4), we get
 \[
 \mathrm{Re}\int_0^1\frac{\lam w_j^+ (sr,t)}{iaB(t)}\frac{ds}{s}
 =\mathrm{Re}\left(r^{\sig_j}q_j^+(t)\right)
 \]
 A similar relation holds for the integral of $w_j^-$ and the series
 expansion (10.5) follows}

\Rem{10.2}{ A consequence of this theorem and of the asymptotic
expansion of the spectral values $\sig_j$, given in Theorem 4.1,
imply that the number $\lam$ is an invariant for the operator $P$ in
the following sense: Suppose that
 \[
 P_1=L_1\ov{L_1}+\ov{\lam_1}\beta_1(t)L_1+\lam_1\ov{\beta_1(t)}\,\ov{L_1}
 \]
 is generated by a vector field $L_1$ with invariant $\lam_1=a_1+ib_1\in \R^++i\R$
 and such that for every $k\in\Z^+$, there is a diffeomorphism, $\Phi^k$,
 in a neighborhood of the circle $S_0$ such that that $\Phi^k_{\ast}P$ is a multiple
 of $P_1$, then $\lam =\lam_1$}

\subsection{A maximum principle}
 We use the series representation of Theorem 10.1 to obtain a
 maximum principle when the spectrum satisfies a certain condition.

 Recall that the function $B(t)$ satisfies $\ind(B)=-k$, where $k\in\Z$
 is defined by (9.1).
We will say that the operator $P$ satisfies hypothesis $\mathcal{H}$
if the spectrum of $\cl$ satisfies the following conditions.
\begin{enumerate}
\item[$\mathcal{H}_1:$] $\dis\mathrm{Re}(\sigpmj)\le 0\,
\Longrightarrow\, j\le -k$.
\item[$\mathcal{H}_2:$]
$\dis\mathrm{Re}(\sigpmj)=\mathrm{Re}(\sig_m^\pm)\,
\Longrightarrow\, \sigpmj =\sig_m^\pm$
\end{enumerate}
Thus $P$ satisfies $\mathcal{H}$ means that the projection of
$\specl$ into $\R$ is injective and that the basic solutions $w$ of
$\cl$ with  positive orders have winding numbers $\ind(w)>k$.

{\theo Suppose that the operator $P$ satisfies $\mathcal{H}$. Let
$U\subset\R^+\times\cir$ be  open, bounded, and such that
$A(0,R)\subset U$ for some $R>0$. If $u\in C^0(\ov{U})$ satisfies
$Pu=0$, then the value of $u$ on $S_0$ is not an extreme value of
$u$. Thus the maximum and minimum of $u$ occur on $\pa U\backslash
S_0$.}

 \vspace{.2cm}

\Proof{ Let $\tau >0$ be the order along $S_0$ of the $L$-potential
of $u$.  It follows from Theorem 10.1 that
 \begeq
 u(r,t)-u(0,0)=\sum_{\mathrm{Re}(\sigpmj)=\tau}\!\!
  c_j^\pm \mathrm{Re}\left(r^{\sigpmj}q_j^\pm (t)\right)\,
  +o(r^\tau).
 \stopeq
 We consider two cases depending on whether $\tau$ is a spectral
 value of $\cl$ or $\tau$ is  only the real part of a spectral value.
   Note that it follows from $\mathcal{H}_2$ that
 the sum in (10.6) consists of either one term or two terms. It
 has one term, if $\tau$ is a spectral value with multiplicity one.
 It has two terms if $\tau$ is a spectral value with multiplicity
 two or if $\tau$ is not a spectral value.

 If $\tau$ is
 spectral value (say with multiplicity 2), then the corresponding
 basic solutions have the form $r^\tau f_{j}^\pm (t)$
 with winding number $j>-k$ (by condition $\mathcal{H}_1$).
 After replacing,
 in (10.6), the functions
 $q_j^\pm$ by their expressions given in (10.3), we find that
\[
u(r,t)-u(0,0)= r^\tau \mathrm{Re}\left(
 \frac{\lam}{ia\tau}\,\frac{c_j^+f_j^+(t)+c_j^-f_j^-(t)}{B(t)}\right)\,
 +o(r^\tau)
\]
 Recall that the functions $f_j^+$ and $f_j^-$ are $\R$-independent
 solutions of the differential equation (2.3). Thus,
 $c_j^+f_j^++c_j^-f_j^- $ has winding number $j$ and consequently
 \[
 \ind\left(
 \frac{\lam}{ia\tau}\,\frac{c_j^+f_j^+(t)+c_j^-f_j^-(t)}{B(t)}
 \right)=j+k \, >0
 \]
 (we have used the fact that $\ind(B)=-k$). Since the winding number
 is positive, the real part changes sign. This implies that
 $u(r,t)-u(0,0)$ changes sign (for $r$ small) and $u(0,0)$ is not an
 extreme value. The proof for the case when $\tau$ is a spectral value with
 multiplicity one is similar.

 If $\tau$ is not a spectral value, then there is a  unique
 spectral value $\sigma =\tau +i\mu$ with $\mu\ne 0$. The corresponding
 basic solutions $w_j^\pm$ have winding number $j>-k$. After substituting, in
 (10.6), the functions  $q_j^\pm$  by their
 expressions given in (10.4) we find
 \[
 u(r,t)-u(0,0) =r^\tau\mathrm{Re}\left[
 \frac{r^{i\mu}}{ia\sig}\left(
 D\lam\frac{\phi_j(t)}{B(t)}-\ov{D\lam}\frac{\psi_j(t)}{\ov{B(t)}}
 \right)\right] \, +o(r^\tau)\, ,
 \]
 where $D=c^++ic^-$ and where $\phi_j$ and $\psi_j$ are the
 components of the basic solutions. We have $|\phi_j|>|\psi_j|$
 and $\ind(\phi_j)=j$. The same argument as before shows that
 $u(r,t)-u(0,0)$ changes sign, as the real part of a function with
 winding number
 $j+k>0$ }

 \Rem{10.3}{If the condition $\mathcal{H}$ is not satisfied, then
 equation
 $Pu=0$ might have solutions with extreme values on $S_0$. Consider
 for example the case in which the function $\beta (t)$ is given by
 (10.1) with $k=1$ (In this example we have
  $\ind(B)=1$).
 The operator $P$ does not satisfy $\mathcal{H}_1$.
 Indeed, $2a$ is a spectral value,
 corresponding to the basic solution $r^{2a}B(t)$, with winding
 number $1$. The corresponding basic solution $r^{2a}$ has minimum
 value 0 and it is attained on $S_0$.}

 \section{The nonhomogeneous equation $Pu=F$}
 We construct here integral operators for the equation $Pu=F$.
 A similarity principle between the solutions of $Pu=0$ and those of
 a semilinear equation is then obtained through these
 operators.

 Let $\widehat{\Om}_1$ and $\widehat{\Om}_2$ be the functions
 given by (7.7). Define the function $S(z,\zeta)$ by
 \begeq
 S(z,\zeta)=\mathrm{Re}\left[
 \frac{-\lam}{2\pi aiB(t)}\int_0^1\!\!
 \left(\widehat{\Om}_1(sr,t,\zeta)B(\ta)+
 \ov{\widehat{\Om}_2(sr,t,\zeta)B(\ta)} \right)
 \frac{ds}{s}\right]
 \stopeq
and  the integral operator $\mathbb{K}$ by
 \begeq
 \mathbb{K}F(r,t)=\iint_{A(0,R)}\!\! S(r,t,\rho,\ta)F(\rho,\ta)
 \frac{d\rho d\ta}{\rho}
 \stopeq
We have the following theorem.

{\theo If $p>2$ and $R>0$, then there exist positive constants $C$
and $\delta$ such that
 $\mathbb{K}:\, L^p_a(A(0,R))\,\longrightarrow\, C^0(\ov{A(0,R)})$
 has the following properties
 \[
 P(\mathbb{K}F) =F,\quad \mathbb{K}F(0,t)=0,\quad
  \mathrm{and}\quad
 |\mathbb{K}F(r,t)|\,\le\, CR^\delta  ||F||_{p,a}\, .
 \]}
 \vspace{.2cm}

\Proof{
 For an $\R$-valued function $F\in L^p_a(A(0,R))$, with
$p>2$, consider
 \begeq
 \widehat{T}_R(B(t)F(z))=\frac{-1}{2\pi}
 \iint_{A(0,R)}\!\!\left(
 \widehat{\Om}_1(z,\zeta)B(\ta)+
 \ov{\widehat{\Om}_2(z,\zeta)B(\ta)}\right) F(\zeta)
 \frac{d\rho d\ta}{\rho}
 \stopeq
We know, from Theorem 7.4,  that $\widehat{T}_R(BF)\in
C^0(\ov{A(0,R)})$ satisfies
 \[\begar{c}
 \cl \widehat{T}_R(BF)=BF,\quad
  \widehat{T}_R(BF)(0,t)=0,\quad\mathrm{ and}\\
 |\widehat{T}_R(BF)(r,t)|\le CR^\delta ||BF||_{p,a}\le C||B||_0
 R^\delta  ||F||_{p,a}\quad\forall F\in L^p_a(A(0,R)).
 \stopar\]
 for some positive constants  $C_1$ and $\delta$.
 Furthermore, it follows from (11.1), (11.2), and (11.3) that
 \begeq
 \mathbb{K}(F)(r,t)=\mathrm{Re}\left(\frac{\lam}{ia}
 \int_0^1\frac{\widehat{T}_R(BF)(sr,t)}{B(t)}\,\frac{ds}{s}\right).
 \stopeq
  Then, from Proposition 9.1, we conclude that
  $\widehat{T}_R(BF)$ is the $L$-potential of
 $\mathbb{K}(F)$. The conclusion of the
 theorem follows from (11.4) and from the  properties of $\widehat{T}_R$ }

 To established a similarity principle between the solutions of
 $Pu=0$ and those of an associated semilinear equation, we need to
 use the modified kernels of section 5. For $j\in\Z$, let
  $\Om_{j,1}^\pm$ and  $\Om_{j,2}^\pm$ be the kernels given by (5.11)
  and (5.12). Define $S_{j}^\pm$ by
  \begeq
 S_{j}^\pm (z,\zeta)=\mathrm{Re}\left[
 \frac{-\lam}{2\pi aiB(t)}\int_0^1\!\!
 \left(\Om_{j,1}^\pm (sr,t,\zeta)B(\ta)+
 \ov{\Om_{j,2}^\pm (sr,t,\zeta)B(\ta)} \right)
 \frac{ds}{s}\right]
 \stopeq
and the operator $\mathbb{K}_{j}^\pm$ by
 \begeq
 \mathbb{K}_{j}^\pm F(r,t)=\iint_{A(0,R)}\!\!
 S_{j}^\pm (r,t,\rho ,\ta)F(\rho,\ta)\frac{d\rho d\ta}{\rho}.
 \stopeq
 The operators $T^\pm_j$, defined in (8.2), and $\mathbb{K}_j^\pm$ are related
 by
 \begeq
 \mathbb{K}_j^\pm F(r,t)=\mathrm{Re}\left(
 \frac{\lam}{ia}\int_0^1\frac{T_j^\pm (BF)(sr,t)}{B(t)}
 \,\frac{ds}{s}\right)
 \stopeq
 The operator $\mathbb{K}_j^\pm$ acts on the Banach
 space $r^{\sigpmj}L^p_a(A(0,R))$, defined in Section 8, and
 produces continuous functions that vanish along $S_0$.
 More precisely,  define the Banach space $r^{\sigpmj}\mathcal{E}(A(0,R))$
 to be the
 set of functions $v(r,t)$  that are in $C^1(A(0,R))$ such that
 $\dis\frac{v}{r^{\sigpmj}}$ and $\dis\frac{Lv}{r^{\sigpmj}}$ are
 bounded functions in $A(0,R)$. The norm of $v$ is
 \[
 \Arrowvert v\Arrowvert_{r^{\sigpmj}\mathcal{E}}=
 \big\Arrowvert\frac{v}{r^{\sigpmj}}\big\Arrowvert_0+
\big\Arrowvert\frac{Lv}{r^{\sigpmj}}\big\Arrowvert_0\, .
 \]
  The next theorem can be proved by using Theorem 8.2 and arguments
  similar to those used in the proof of Theorem 11.1.

 {\theo The operator
 \[
 \mathbb{K}_j^\pm :\,
 r^{\sigpmj}L^p_a(A(0,R))\,\longrightarrow\,
 r^{\sigpmj}\mathcal{E}(A(0,R))
 \]
 satisfies $P\mathbb{K}_j^\pm F=F$ and
 \[
 \Arrowvert \mathbb{K}_j^\pm F\Arrowvert_{r^{\sigpmj}\mathcal{E}}
 \,\le\, CR^\delta
 \Arrowvert F\Arrowvert_{p,a,\sigpmj}
 \]
 for some positive constants $C$ and $\delta$.}

 \vspace{.2cm}

 Let $f_0(r,t)$, $f_1(r,t)$, and $f_2(r,t)$ be bounded functions
 in $A(0,R)$ and let
 $g_1(r,t,u,w)$ and $g_2(r,t,u,w)$ be bounded functions in $A(0,R)\times
\R\times\C$.  Define the function $H$ by
 \begeq
 H(r,t,u,w)=uf_0+wf_1+\ov{w}f_2+|u|^{1+\alpha}g_1+|w|^{1+\alpha}g_2
 \stopeq
 with $\alpha >0$. For $\ep >0$, consider the semilinear equation
 \begeq
 Pu=r^\ep \mathrm{Re}\left(H(r,t,u,Lu)\right)\, .
 \stopeq
We have the following  similarity result between the solutions of
(11.9) and those of the equation $Pu=0$.

{\theo For a given function $H$ defined by $(11.8)$, there exists
$R>0$ such that, for every $u\in C^0(\ov{A(0,R)})$ satisfying $Pu=0$
and $u=0$ on $S_0$, there exists a function $m\in C^0(A(0,R))$
satisfying
  \[
  C_1\le m(r,t) \le C_2\qquad\forall (r,t)\in A(0,R)
 \]
 with $C_1$ and $C_2$ positive constants, such that the function
 $v=mu$ solves equation $(11.9)$.}

 \vspace{.2cm}

\Proof{ Let $u$ be a solution of $Pu =0$ with order $\tau >0$ along
$S_0$. Then there is $\sigpmj\in\specl$ such that $\tau
=\mathrm{Re}(\sigpmj)$. Hence,  $u\in
r^{\sigpmj}\mathcal{E}(A(0,R_0))$ for some $R_0>0$. Consider the
operator
 \[
 \mathcal{Q}:\  r^{\sigpmj}\mathcal{E}(A(0,R_0))\,\longrightarrow
 \, r^{\sigpmj}\mathcal{E}(A(0,R_0))
 \]
 given by
 $\mathcal{Q}v=\mathbb{K}_j^\pm(r^\ep\mathrm{Re}(H(r,t,v,Lv)))$.
 It follows from (11.8) that the function
 $r^\ep\mathrm{Re}(H(r,t,v,Lv))$
 is in the space $r^{\sigpmj}L^p_a$.  Now,  Theorems 11.2, 8.2, and
 relation (11.7), imply  that
 \[\begar{l}
 P\mathcal{Q}v(r,t)=r^\ep\mathrm{Re}(H(r,t,v,Lv)),\\
 L\mathcal{Q}v(r,t)=T_j^{\pm}\left[B(t)r^\ep\mathrm{Re}(H(r,t,v,Lv))\right].
 \stopar\]
 Consequently,
  $\Arrowvert \mathcal{Q} v\Arrowvert_{r^{\sigpmj}\mathcal{E}}
 \le CR_0^\delta
 \Arrowvert v\Arrowvert_{r^{\sigpmj}\mathcal{E}}$. If $R_0$ is small
 enough, we have $\Arrowvert \mathcal{Q} \Arrowvert <1$, and we can define
 the resolvent $\mathcal{F}=(I-\mathcal{Q})^{-1}$.  It is easily
 checked that for the solution $u$ of $Pu=0$ as above, the function
 $v=\mathcal{F}(u)$ solves equation (11.9) and $m=u/v$ is bounded away
 from 0 and $\infty$}

 \section{Normalization of a Class of Second Order Equations with a
 Singularity }
 This section deals with the normalization of a class of second
 order operators $\mathbb{D}$ in $\R^2$ whose coefficients vanish at
 a point. To such an operator, a complex number
 $\lam=a+ib\in\R^++i\R$ is invariantly associated. It is then shown that
 the operator
 $\mathbb{D}$ is conjugate, in a punctured neighborhood of the
 singularity, to a unique operator $\mathbb{P}$ given by (9.2). The
 properties of the solutions of the equations
 corresponding to $\mathbb{D}$ are, thus, inherited from
  the solutions of the equations for $\mathbb{P}$ studied in sections
   10 and 11.

 Let $\mathbb{D}$ be the second order operator  given in a
 neighborhood of $0\in\R^2$ by
 \begeq
 \mathbb{D}u=a_{11}u_{xx}+2a_{12}u_{xy}+a_{22}u_{yy}+
 a_1u_x+a_2u_y
 \stopeq
 where the coefficients $a_{11},\cdots ,\ a_2$ are $C^\infty$, real-valued
 functions vanishing at 0, with $a_{11}$ nonnegative, and
 \begeq
 C_1\,\le\,
 \frac{ a_{11}(x,y)a_{22}(x,y)-a_{12}(x,y)^2}{(x^2+y^2)^2}\,\le\,
 C_2
 \stopeq
 for some positive constants $C_1<C_2$.  It follows in particular
 that $a_{11}$ and $a_{22}$ vanish to second order at 0.
 Let $A$ and $B$ be the
 functions defined for $(x,y)\ne 0$ by
 \begeq\begar{ll}
 A(x,y) & = \dis
 \frac{(x^2+y^2)\sqrt{a_{11}a_{22}-a_{12}^2}}{a_{11}y^2-2a_{12}xy+a_{22}x^2}
 \\
 B(x,y) & = \dis
 \frac{(a_{22}-a_{11})xy+a_{12}(x^2-y^2)}{a_{11}y^2-2a_{12}xy+a_{22}x^2}
 \stopar\stopeq
  Note that it follows from (12.2) that these functions are bounded and
 $A$ is positive. Let
 \begeq
 \mu =\frac{1}{2\pi}
 \lim_{\rho\to 0^+}\int_{C_\rho}
 \frac{A(x,y)-iB(x,y)}{x^2+y^2}\, (xdy-ydx)\, ,
 \stopeq
where $C_\rho$ denotes the circle with radius $\rho$ and center 0.
 We will prove that  $\mu\in \R^++i\R$  is well defined
 and it is an invariant for the operator $\mathbb{D}$.

We will be using the
 following normalization theorem for a class of vector fields
 in a neighborhood of a characteristic curve.

{\theo  Let $X$ be a $C^\infty$ complex vector field in $\R^2$
satisfying the following conditions in a neighborhood of a smooth,
simple, closed curve $\Sigma$:
 \begin{enumerate}
 \item[$(\imath )$]
  $X_p\wedge\ov{X}_p\ne 0$ for every  $p\notin\Sigma$;
 \item[$(\imath\imath)$]
 $X_p\wedge\ov{X}_p$ vanishes to first order for $ p\in\Sigma$; and
 \item[$(\imath\imath\imath)$]
 $X$ is tangent to $\Sigma$.
 \end{enumerate}
 Then there exist an open tubular neighborhood $U$ of $\Sigma$, a positive
 number $R$,  a unique complex number $\lam\in\R^++i\R$, and a
 diffeomorphism
 \[
 \Phi :\, U\,\longrightarrow\, (-R,\ R)\times\cir
 \]
 such that
 \[
 \Phi_\ast X =m(r,t) \left[\lam\dd{}{t}-ir\dd{}{r}\right]
 \]
 where $m(r,t)$ is a nonvanishing function. Moreover, when
  $\lam\not\in\mathbb{Q}$, then
 for any given $k\in\Z^+$, the diffeomorphism $\Phi$
 and the function $m$ can be taken to be of class $C^k$.   }

 \vspace{.2cm}

 This normalization Theorem was proved in {\cite{Mez-JFA}}
 when $\lam\in\C\backslash\R$.
 When $\lam\in\R$, only a $C^1$-diffeomorphism $\Phi$ is achieved in
 {\cite{Mez-JFA}}.
 A generalization is obtained by Cordaro and Gong in
 {\cite{Cor-Gon}} to include $C^k$-smoothness of $\Phi$ when
 $\lam\in\R\backslash\mathbb{Q}$. It is also proved in {\cite{Cor-Gon}},
 that, in
 general, a $C^\infty$-normalization cannot be achieved.

 We will be using polar
 coordinates $x=\rho\cos\ta$, $y=\rho\sin\ta$ and we will denote
 this change of coordinates by $\Psi$. Thus,
 \[
 \Psi  :\R^2\backslash 0\, \longrightarrow\, \R^+\times\cir\,
 ,\qquad \Psi(x,y)=(\rho,\ta)\, .
 \]

 {\theo Let $\mathbb{D}$ be the second order operator
 given by $(12.1)$ whose coefficients  vanish at $0$
 and satisfy condition $(12.2)$. Then there is a neighborhood
 $U$ of the circle $\{0\}\times\cir$ in $[0,\ \infty)\times\cir$,
 a positive number $R$, a diffeomorphism
 \[
 \Phi :\, U\,\longrightarrow\,  [0,\ R)\times\cir
 \]
 sending  $\{0\}\times\cir$ onto itself,  such that
 \begeq
 (\Phi\circ\Psi)_\ast\mathbb{D}=
 m(r,t)\left[ L\ov{L}+\mathrm{Re}(\beta(r,t)L)\right]
 \stopeq
 where $m, \beta$ are differentiable functions with
 $m(r,t)\ne 0$ for every $(r,t)$ and
 \[
L=\lam\dd{}{t}-ir\dd{}{r}
 \]
 with $\lam=\dis\frac{1}{\mu}$ and $\mu$ given by $(12.4)$.
 Moreover, if the invariant
 $\mu\notin\mathbb{Q}$, then for every $k\in\Z^+$, the diffeomorphism
 $\Phi$, and the functions
 $m$, and $\beta$ can be chosen to be of class $C^k$.  }

 \vspace{.2cm}

 \Proof{ We start by rewriting $\mathbb{D}$ in polar coordinates:
 \begeq
 \mathbb{D}u=Pu_{\ta\ta}+2Nu_{\rho\ta}+Mu_{\rho\rho} +
 Qu_{\rho}+Tu_\ta
 \stopeq
where
 \[\begar{ll}
 P & = \dis\frac{1}{\rho^2}\left[
  a_{11}\sin^2\ta -2a_{12}\sin\ta\cos\ta +a_{22}\cos^2\ta\right]\\
  N & =\dis\frac{1}{\rho}\left[ -
  a_{11}\sin\ta\cos\ta +a_{12}(\cos^2\ta-\sin^2\ta) +
  a_{22}\cos\ta\sin\ta\right]\\
 M&=\dis a_{11}\cos^2\ta +2a_{12}\sin\ta\cos\ta +a_{22}\sin^2\ta\\
 Q & =\dis\frac{1}{\rho}\left[
  a_{11}\sin^2\ta -2a_{12}\sin\ta\cos\ta +a_{22}\cos^2\ta\right]
  +a_1\cos\ta +a_2\sin\ta\\
  T & = \dis\frac{1}{\rho^2}\left[
  a_{11}\sin\ta\cos\ta +a_{12}(\sin^2\ta-\cos^2\ta) -
  a_{22}\sin\ta\cos\ta\right]-\frac{1}{\rho}
  (a_1\sin\ta +a_2\cos\ta)
 \stopar\]
 Condition (12.2) implies that there is a constant $C_0>$ such
 that
 \[
 M(\rho,\ta)\ge C_0\rho^2\quad\mathrm{and}\quad
 P(\rho,\ta)\ge C_0\qquad\forall (\rho,\ta).
 \]
 We define the following $C^\infty$ functions (of $(\rho,\ta)$)
 \[
 N_1=\frac{N}{\rho P},\quad M_1=\frac{M}{\rho^2 P},
 \quad Q_1=\frac{Q}{\rho P},
 \quad T_1=\frac{T}{ P} .
 \]
 In terms of these function, (12.2) takes the form
 \begeq
 M_1(\rho,\ta)-N_1^2(\rho,\ta)\,\ge\,  C_2,\qquad\forall
 (\rho,\ta)\in [0,\ R_1]\times\cir\, ,
 \stopeq
 and  (12.6) becomes
 \begeq
 \frac{\mathbb{D}u}{P}=
 u_{\ta\ta}+2\rho N_1u_{\rho\ta}+\rho^2M_1u_{\rho\rho}+
 \rho Q_1u_\rho +T_1u_\ta
 \stopeq

  Let $X$ be the $C^\infty$ complex vector field defined by
 \begeq
 X=\dd{}{\ta}-\rho g(\rho,\ta)\dd{}{\rho}
 \stopeq
 with $g=N_1+i\sqrt{M_1-N_1^2}$.
 Although  we will use $X$ for $\rho\ge 0$, the vector field
 $X$  is defined
 in a neighborhood of $\{ 0\}\times\cir$ in $\R\times\cir$.
 By using $X$ and its  complex conjugate $\ov{X}$, we find that
 \begeq
 X\ov{X} u=u_{\ta\ta}+2\rho N_1u_{\ta\rho}+\rho^2M_1u_{\rho\rho}+
 \rho f u_\rho
 \stopeq
 where
 \[
 f=\frac{X(\rho\ov{g})}{\rho}=-|g|^2+X(\ov{g}).
 \]
  We also have
 \begeq
 \rho u_\rho =\frac{Xu-\ov{X} u}{r-\ov{g}}\quad\mathrm{and}\quad
 u_\ta =\frac{g\ov{X}u-\ov{g}Xu}{g-\ov{g}}\, .
 \stopeq
 It follows from (12.8), (12.10) and (12.11) that
 \begeq
 \frac{\mathbb{D} u}{P}=
 X\ov{X} u-\frac{f-Q-1+\ov{g}T_1}{g-\ov{g}}\, Xu +
 \frac{f-Q-1+gT_1}{g-\ov{g}}\, \ov{X}u
 \stopeq
 Since the coefficients of $\mathbb{D}$ and the function $u$
 are $\R$-valued,
 then the right hand side of (12.12) is real valued and can be
 written  as
 \begeq
 2\frac{\mathbb{D} u}{P}=X\ov{X} u +\ov{X}\, X u +
 B(\rho,\ta)Xu+\ov{B(\rho,\ta)}\, \ov{X} u
 \stopeq
 with
 \[
 B(\rho,\ta)=-\frac{f+\ov{f}-2Q_1+2\ov{g}T_1}{g-\ov{g}}\, .
 \]
 Now, for the vector field $X$, we have
  \[
 X\wedge\ov{X} =\rho(\ov{g}-g)\dd{}{\ta}\wedge\dd{}{\rho}
 =-2i\rho\sqrt{M_1-N_1^2}\, \dd{}{\ta}\wedge\dd{}{\rho},
 \]
  and so $X$ satisfies the conditions of Theorem 12.1 and
  therefore it can be
  normalized. In our setting, the
   invariant $\lam$ is given by (see {\cite{Mez-JFA}})
  $\lam =1/\widetilde{\mu}$ where
  \[
 \widetilde{\mu}=\frac{1}{2\pi}\int_0^{2\pi}\left(
 \sqrt{M_1(0,\ta)-N_1^2(0,\ta)}-
 iN_1(0,\ta)\right)d\ta\,=\mu ,
  \]
 and where $\mu$ is given by (12.4).
 Hence, there is a diffeomorphism $\Phi$ defined in a neighborhood
 of $\rho =0$ in $\R\times\cir$ onto a cylinder $(-R,\ R)\times \cir$
 such that $\Phi_\ast X=m(r,t)L$ with $L$ as in the Theorem and
 $m$ a nonvanishing function. Finally, it follows from this normalization
 of $X$ that, in the $(r,t)$ coordinates, expression (12.13) becomes
 \begeq
  2\frac{\mathbb{D} u}{P}=
  2|m|^2L\ov{L}u+(mB+\ov{m}\,\ov{L}m)\, Lu +
  (\ov{m}\,\ov{B}+mL\ov{m})\,\ov{L} u\, .
 \stopeq
 This completes the proof of the theorem }

\end{document}